\pdfoutput=1
\RequirePackage{ifpdf}
\ifpdf 
\documentclass[pdftex]{sigma}
\else
\documentclass{sigma}
\fi

\usepackage[all]{xy}

\numberwithin{equation}{section}

\newtheorem{theor}{Theorem}[section]
\newtheorem{prop}[theor]{Proposition}
\newtheorem{Corollary}[theor]{Corollary}
\newtheorem{Conjecture}[theor]{Conjecture}
\newtheorem{Lemma}[theor]{Lemma}

{\theoremstyle{definition}
\newtheorem{defi}[theor]{Definition}
\newtheorem{defis}[theor]{Definitions}

\newtheorem{Remark}[theor]{Remark}
\newtheorem{compari}[theor]{Comparison}

\newtheorem{exa}[theor]{Example}
\newtheorem{question}[theor]{Question}
\newtheorem{questions}[theor]{Questions}
}

\newcommand{\N}{\mathbb{N}}
\newcommand{\Z}{\mathbb{Z}}
\newcommand{\R}{\mathbb{R}}
\newcommand{\C}{\mathbb{C}}
\newcommand{\f}{\rightarrow}

\newcommand{\CDCD}{\operatorname{CD}}
\newcommand{\BG}{\operatorname{BG}}
\newcommand{\CAT}{\operatorname{CAT}}

\newcommand{\Ker}{\operatorname{Ker}}

\newcommand{\Vol}{\operatorname{Vol}}

\newcommand{\Cat}{\operatorname{Cat}}
\newcommand{\Ric}{\operatorname{Ric}}

\newcommand{\Ent}{\operatorname{Ent}}

\newcommand{\Max}{\operatorname{Max}}
\newcommand{\sys}{\operatorname{sys}}
\newcommand{\dias}{\operatorname{Dias}}
\newcommand{\diam}{\operatorname{diam}}
\newcommand{\inj}{\operatorname{inj}}
\newcommand{\g}{\operatorname{\gamma}}
\newcommand{\e}{\operatorname{\varepsilon}}
\def\Smin{\textrm{Scal}_\mathrm{min}}

\begin{document}

\newcommand{\arXivNumber}{2010.08207}

\renewcommand{\thefootnote}{}

\renewcommand{\PaperNumber}{046}

\FirstPageHeading

\ShortArticleName{On Scalar and Ricci Curvatures}

\ArticleName{On Scalar and Ricci Curvatures\footnote{This paper is a~contribution to the Special Issue on Scalar and Ricci Curvature in honor of Misha Gromov on his 75th Birthday. The~full collection is available at \href{https://www.emis.de/journals/SIGMA/Gromov.html}{https://www.emis.de/journals/SIGMA/Gromov.html}}}

\Author{Gerard BESSON and Sylvestre GALLOT}

\AuthorNameForHeading{G.~Besson and S.~Gallot}

\Address{CNRS-Universit\'e Grenoble Alpes, Institut Fourier, CS 40700, 38058 Grenoble cedex 09, France}
\Email{\href{mailto:g.besson@univ-grenoble-alpes.fr}{g.besson@univ-grenoble-alpes.fr}, \href{mailto:sylvestre.gallot@univ-grenoble-alpes.fr}{sylvestre.gallot@univ-grenoble-alpes.fr}}

\ArticleDates{Received October 19, 2020, in~final form April 05, 2021; Published online May 01, 2021}

\Abstract{The purpose of this report is to acknowledge the influence of M.~Gromov's vision of geometry on our own works. It~is two-fold: in the first part we aim at describing some results, in~dimension $3$, around the question: which open $3$-manifolds carry a complete Riemannian metric of positive or non negative scalar curvature? In the second part we look for weak forms of the notion of ``lower bounds of the Ricci curvature'' on non necessarily smooth metric measure spaces. We~describe recent results some of which are already posted in~[arXiv:1712.08386] where we proposed to use the volume entropy. We~also attempt to give a new \emph{synthetic version} of Ricci curvature bounded below using Bishop--Gromov's inequality.}


\Classification{51K10; 53C23; 53C21; 53E20; 57K30}

\renewcommand{\thefootnote}{\arabic{footnote}}
\setcounter{footnote}{0}

\begin{flushright}
\emph{Happy Birthday Misha}
\end{flushright}

\section{Introduction}

The purpose of this text is to acknowledge the influence of M.~Gromov's vision of geometry on our own works. It~started in the early 80's with the french version of the so-called ``green book''~\cite{Gr1}. We~then discovered some of his previous articles and followed his progression in the new realm that he was building.

The following text is two-fold. In~the first part (Section~\ref{chap:scalar}) we aim at describing some results, in~dimension $3$, around the question: which open $3$-manifolds carry a complete Riemannian metric of positive or non negative scalar curvature? The influential article is the joint work by~M.~Gromov and B.~Lawson~\cite{GromovLawson} and, in~particular, the beautiful Chapter~10 where the relation between a $3$-manifold carrying a metric with positive scalar curvature and the stable minimal surfaces that it contains is subtly exploited. We~focus on two types of~$3$-manifolds: the so-called decomposable $3$-manifolds and the contractible
$3$-manifolds. For the first family the results are extensions to open $3$-manifolds of a corollary of Perelman's works asserting that a closed $3$-manifold has a metric of positive scalar curvature if and only if it is
a connected sum of spherical manifolds (see Section~\ref{sec:decomposable} for the details). The second part of this section is about the family of contractible open $3$-manifolds, which is very rich, and found its origin
in early works by J.H.C.~Whitehead (see Section~\ref{intro:scalar}). The~goal here would be to show that only $\mathbb{R}^3$ carries a~complete metric of positive scalar curvature. It~is not yet achieved but the results are already quite striking and strongly influenced by the article~\cite{GromovLawson}. In~both parts the results culminate with works by Jian Wang.

So far we have been looking for smooth metrics on smooth Riemannian manifolds. In~view of~some recent developments we could search for non smooth metrics on these contractible spaces, for example. These raises quite a few exciting questions which are asked at the end of Section~\ref{chap:scalar}.

In the second part (Section~\ref{chap:BG}) we are indeed interested in non necessarily smooth metric measure spaces. The~underlying philosophy founds its origin as far as in the article~\cite{BCG-GAFA}. However, most of the results that are described are proved in the preprint~\cite{BCGS} and some others will be in a forthcoming article, by the same authors, to be posted soon. In~this Section~\ref{chap:BG} we look for weak forms of lower bounds of the Ricci curvature. This already exists, for example curvature-dimension conditions \textit{\`a la} Bakry--Emery, the so-called CD(K,N) condition where $K$ stands for a~lower bound of the Ricci curvature and $N$ for an upper bound on the dimension. The original version was about operators (or semi-groups generated by an operator), typically the Laplacian of a Riemannian manifold and it somehow used a weak version of Bochner formula as a definition of a lower bound on the Ricci curvature. It~was then extended, using measure transportation, to what is now called
\textit{synthetic Ricci curvature} by Lott--Villani and Sturm (see Section~\ref{BGexemples}($b$)). In~\cite{BCGS} we propose to use an upper bound on the volume entropy as a replacement to a lower bound on the Ricci curvature. The~volume entropy of a Riemannian manifold is related to the asymptotic growth of the volume of balls of radius $R$ around a point in its universal cover and its definition is given in~Definition~\ref{Entropies0}.

\looseness=1
By Bishop's comparison theorem it is easy to show that, if the Ricci curvature of a Riemannian metric $g$ satisfies that $\Ric (g)\geq -(n-1)\kappa g$ where $\kappa$ is a non negative real number, then the entropy is bounded above by $(n-1)\sqrt \kappa$. Hence an upper bound on the entropy could serve as a~(very) weak version of a lower bound on the Ricci curvature. Now, a~lower bound on the Ricci curvature, in~the Riemannian case, implies the Bishop--Gromov's comparison theorem which is an easy modification of Bishop's original result which has amazing consequences. This is the case also with (some of) the synthetic version of the Ricci curvature and we show in~\cite{BCGS} that it is also the case, at~least in a weaker form, with the upper bound on the volume entropy (together with other natural assumptions) applied to metric measure spaces. In~Section~\ref{chap:BG} we decided to use this consequence as a definition of the expression
``Ricci curvature bounded below'', yet another version in the spirit of Bakry--Emery $\CDCD(K,N)$ somehow but on the geometric side, where Bochner formula is replaced by Bishop--Gromov's comparison. We~develop both this
approach and the entropy approach applied to metric measure spaces in Section~\ref{chap:BG} and go as far as describing versions of the celebrated Margulis lemma and its consequences on proving finiteness and compactness results. The~reader is referred to this section for precise statements.

It would be clear by now that Misha's works has had a great influence on
ours and we are greatly indebted to him for the exchanges and discussions
that we had along these years. We~also owe a lot to our former advisor Marcel Berger who passed away on October~15, 2016 and his teaching remains in our memory. He admired Gromov, as~we do, and right after Misha settled
in France, at~any questions from us he would answer ``ask Gromov, he knows!''

\section[Scalar curvature: on some works of Schoen--Yau, Gromov--Lawson and Jian Wang]{Scalar curvature: on some works of Schoen--Yau,\\ Gromov--Lawson and Jian Wang}\label{chap:scalar}

\subsection{Introduction}\label{intro:scalar}

This section is a short account of some recent results concerning the existence of metrics with non negative scalar curvature on some $3$-manifolds. The~starting point of our interest is the ``proof'' of the Poincar\'e conjecture given in 1934 by J.H.C.~Whitehead in the article~\cite{White1}. Roughly speaking the scheme is the following (see \href{https://www.math.unl.edu/~mbrittenham2/ldt/poincare.html}{here}):
\begin{itemize}\itemsep=-.2pt
 \item Let $X$ be a simply connected closed $3$-manifolds then, for any point $p\in X$, $X\setminus\{p\}$ is contractible (by simple considerations left to the reader).
 \item The only contractible $3$-manifolds is $\mathbb{R}^3$.
 \item The one-point compactification of~$\mathbb{R}^3$ is $\mathbb{S}^3$.
\end{itemize}
In 1935 he realised that the second step was wrong (see~\cite{White2}) and constructed in~\cite{White3} the first contractible $3$-manifold which is not homeomorphic to $\mathbb{R}^3$, the now-called \textit{Whitehead manifold} which we shall denote by $\rm Wh$.

Let us dream of a version of the Poincar\'e conjecture for open manifolds. What could it~be? Well, an open manifold which has the same homotopy as $\mathbb{R}^3$, that is which is contractible, is~homeomorphic to $\mathbb{R}^3$! The Whitehead manifold is a counter-example to this statement and a~non compact version of the Poincar\'e conjecture is not true! However, the discovery of~$3$-mani\-folds which are contractible but not homeomorphic to $\mathbb{R}^3$, such as the Whitehead manifold, opened a wide playground for topologists and geometers. Mistakes could be useful, indeed a~similar experience happened to Poincar\'e too.

The Whitehead manifold $\rm Wh$ is an open subset of~$\mathbb{S}^3$ whose complement is a closed set, called the \textit{Whitehead continuum}, which looks locally like a product of an interval with a~Cantor
set. It~is easily seen, by construction, that $\rm Wh$ is contractible and the fact that it is not home\-o\-mor\-phic to~$\mathbb{R}^3$ follows from the fact that it is not \textit{simply-connected at infinity} (see a definition in Section~\ref{sec:contractible}). Such manifolds, not homeomorphic to~$\mathbb{R}^3$ but nevertheless contractible, turn out to be plentiful, as~was shown by McMillan (see~\cite{McMil}). In~fact there exist uncountably many such manifolds whereas there are only a countable family of closed topological manifolds. Many other nice properties can be proved and the reader is referred to the literature.

Let us come back to Poincar\'e's conjecture. In~the Clay--Perelman conference, in the lecture entitled ``What is a Manifold'', given at I.H.P. in 2011, Misha Gromov, talking about the Poincar\'e conjecture says ``and so what?'' (see \href{https://www.youtube.com/watch?v=u5DLpAqX4YA}{here}, at~time 13:30). This is typical of Misha's point of view, an information on one object, here being able to recognise the $3$-sphere only, is nice but does not contain enough knowledge. A~general picture is
more important, even if it is incomplete, or, in~other words, getting statistical informations is more useful than knowing a single value. This general picture was given by W.~Thurston in the 70's and is known as Thurston's geometrisation conjecture (see~\cite{Thur} and the post-Perelman literature). It~is a conjectural description of all closed $3$-manifolds and
the way to build or decompose them nicely, which yields the Poincar\'e
conjecture as a corollary. The~picture is so clear and so nice that the attempts to publish counter-examples to the Poincar\'e conjecture (almost) stopped. In~any case it took then less than 30 years to get G.~Perelman's proof of Thurston's conjecture (see~\cite{Per1, Per2, Per3}) developing a technique introduced by R.~Hamilton in~\cite{Ham}.

The beautiful idea of Thurston is to decompose the closed $3$-manifolds into pieces carrying each a very specific Riemannian geometry, like the case of closed surfaces which are either spherical, flat or hyperbolic but with more (in fact 8) possibilities. Now, what could be a~geometrisation conjecture for open $3$-manifolds? There is little hope to be able to state a~reasonable conjecture. For example, the starting point of Thurston's geometrisation conjecture are the Kneser--Milnor and the Jaco--Shalen and Johannson decompositions of closed $3$-manifolds none of which is true for open manifolds (see~\cite{S-T} and~\cite{Mai}). Therefore, the only sensible thing to do is to work out examples and this is the purpose of the next two sections. We~shall describe two extremely different families of open $3$-manifolds and explore their Riemannian geometry with an emphasis on
the scalar curvature.

\subsection{Decomposable 3-manifolds}\label{sec:decomposable}

This is a family of~$3$-manifolds that generalises connected sums. More precisely let $\mathcal X$ be a~class of {\it orientable }closed $3$-manifolds. A manifold $M$ is a (possibly infinite)
{\it connected sum} of members of~$\mathcal X$ if there exists a locally finite graph $G$ (or simply a locally finite tree $T$) and a map $v\mapsto X_v$ which associates to each vertex of~$G$ a copy of some manifold in $\mathcal X$ such that, by removing from each $X_v$ as many $3$-balls as edges incident to $v$ and gluing the thus punctured $X_v$'s to each other
along the edges of~$G$, one obtains a $3$-manifold diffeomorphic to $M$. Precise definitions and several useful comments may be found in~\cite{BBM}.

\begin{figure}[h!]
\centering
{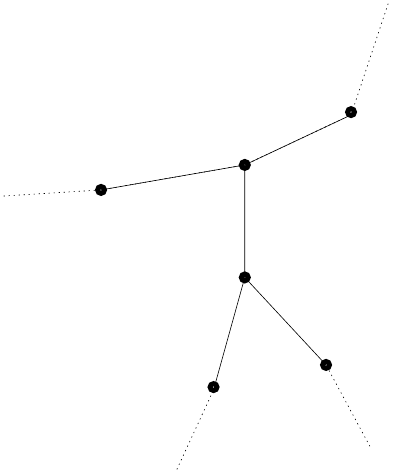\hspace{3cm}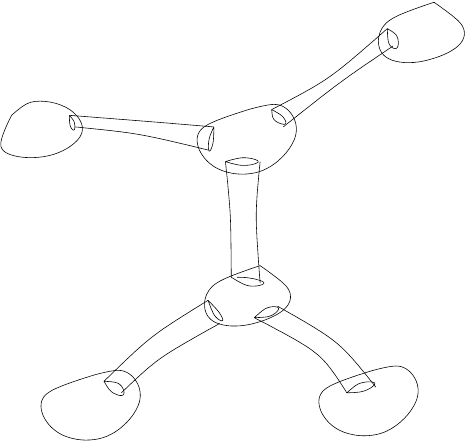}
\end{figure}

We now shall consider Riemannian manifolds with positive scalar curvature. This class has been extensively studied since the works of A.~Lichn\'erowicz, R.~Schoen, S.-T.~Yau, M.~Gromov, B.~Lawson and others (see, e.g.,~the survey articles~\cite{gromov:sign,rosenberg:report}).

Let $(M,g)$ be a Riemannian manifold. We~denote by $\Smin(g)$ the infimum of the scalar curvature of~$g$ and we say that $g$ has
{\it uniformly positive scalar curvature} if $\Smin(g)>0$. If~$M$ is compact then this amounts to saying that $g$ should have positive scalar curvature at each point of~$M$.

A $3$-manifold is {\it spherical} if it admits a metric of positive constant sectional curvature. M.~Gromov and B.~Lawson~\cite{gl:spin} have shown that any compact, orientable $3$-manifold which is a connected sum of orientable spherical manifolds and copies of~$S^2\times S^1$ carries a metric of positive scalar curvature. G.~Perelman~\cite{Per2}, completing pioneering work of Schoen--Yau~\cite{SchoenYau:incompressible} and Gromov--Lawson~\cite{GromovLawson}, proved the converse.

As mentioned before we are mostly interested in the non compact case. We~say that a~Riemannian metric $g$ on $M$ has {\it bounded geometry} if it has bounded sectional curvature
and injec\-tivity radius bounded away from zero. It~follows from the Gromov--Lawson construction that, if the $3$-manifold $M$ is a (possibly infinite) connected sum of spherical manifolds, taken in a~{\it finite} list (see the statement below), and of copies of~$S^2\times S^1$, then $M$ admits
a complete metric of bounded geometry and uniformly positive scalar curvature. The~following theorem by~L.~Bessières, G.~Besson and S.~Maillot shows that the converse holds, generalising Perelman's theorem:

\begin{theor}[\cite{BBM}]\label{thm:positive scalar general}
Let $M$ be a connected, orientable $3$--manifold which carries a complete
Riemannian metric of bounded
geometry and uniformly positive scalar curvature. Then there is a finite collection $\mathcal X$ of spherical manifolds such that $M$ is a connected sum of copies of~$S^2\times S^1$ and members of~$\mathcal X$.
\end{theor}

The proof relies on an extension of the Ricci flow with surgery, called {\it a surgical solution}, that is constructed in the same article~\cite{BBM}. It~is not the purpose of this text to give the details of this construction and the interested reader is referred to~\cite{BBM}. Let us though emphasise that the bounded geometry assumption is the necessary requirement in order to trigger this new Ricci flow and to ensure its existence for all time. The~finiteness of the class $\mathcal X$ is a consequence
of the combination of the lower bound on the injectivity radius and of the uniformly positive scalar curvature. More precisely, let us normalise the metric so that $\Smin = 6$. The~spherical manifolds are quotients of the $3$-sphere with a constant sectional curvature metric and with this
normalisation this constant is equal to $1$; there are infinitely many of
them, the lens spaces for example, and their injectivity radius goes to zero. Hence, a~lower bound on the injectivity radius combined with the lower bound on $\Smin$ reduces $\mathcal X$ to a finite set.

The article~\cite{cwy:taming} contains results closely related to Theorem~\ref{thm:positive scalar general}. Precisely, Theorem~5.1 in~\cite{cwy:taming} has the same conclusion than~\ref{thm:positive scalar general} without the assumption that the metric has bounded geometry but with the assumption that $M$ has a finitely generated fundamental group. Notice that, in~Theorem~\ref{thm:positive scalar general}, if $M$ is a connected sum of
an infinite number of copies of members of~$\mathcal X$ (different from $\mathbb{S}^3$), its fundamental group is not finitely generated. An interesting point is that the proof described in~\cite{cwy:taming} relies on $K$-theory and hence is completely different from the proof of~Theorem~\ref{thm:positive scalar general}. Somehow $K$-theory has to do with the Dirac operator and the relation between the square of the Dirac operator and
the scalar curvature is encoded in the Bochner--Lichn\'erowicz--Weitzenbock formula (see \cite[equation~(2.14), p.~112]{GromovLawson}).

Recently Jian Wang announced the proof of the following result which improves both Theorem~\ref{thm:positive scalar general} and the result in~\cite{cwy:taming}.

\begin{theor}[\cite{Wang3}]\label{Wang3}
A complete connected orientable $3$-manifold of uniformly positive scalar
curvature is homeomorphic to a~$($possibly infinite$)$ connected sum of spherical $3$-manifolds and copies of~$S^2\times S^1$.
\end{theor}

The progress here is that there is neither the assumption that $g$ has bounded geometry nor that the fundamental group of~$M$ is finitely generated. The~proof is yet another one; it uses the case of closed manifolds, proved by G.~Perelman~\cite{Per3}, but neither a version of the Ricci flow for open manifolds nor $K$-theory. The~theory of minimal surfaces and the special relations between stable minimal surfaces and the scalar curvature of the ambient metric in dimension $3$ is essential.

\subsection[Contractible 3-manifolds]
{Contractible $\mathbf{3}$-manifolds}\label{sec:contractible}

\begin{figure}[!h]
\centering\includegraphics[scale=0.8]{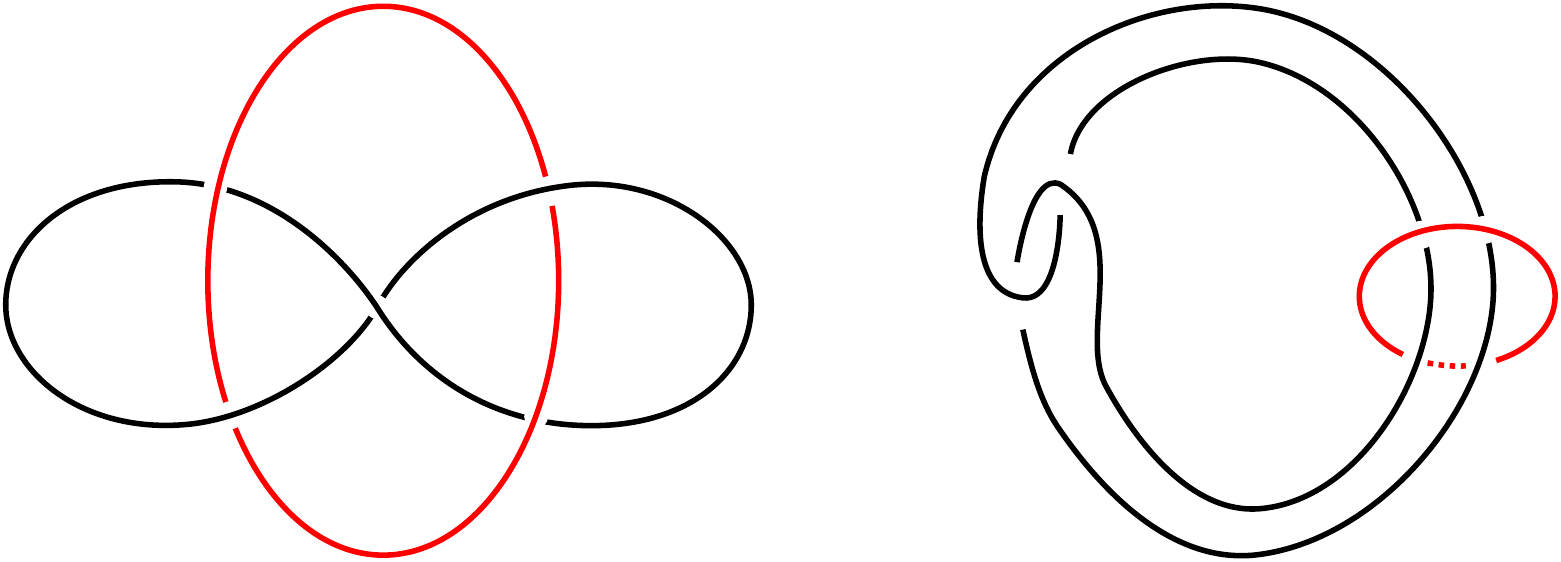}
\caption{Whitehead link.}
\label{Fig1}
\end{figure}

We now come to the second series of examples, namely the contractible open $3$-manifolds which are not homeomorphic to $\mathbb{R}^3$. They are somehow opposite to the class of decomposable manifolds in the sense that they are homotopically trivial and their topology is hidden in their structure at infinity.
Let us recall the construction of the main example, the Whitehead manifold. We~start with the Whitehead link which is a link with two components illustrated in Figure~\ref{Fig1} below in~two different ways. Notice that this
link is symmetric; this means that there is an isotopy of the ambient space which reverses the roles played by the black and red curves.

\begin{figure}[!h]
\centering
\includegraphics[scale=0.8]{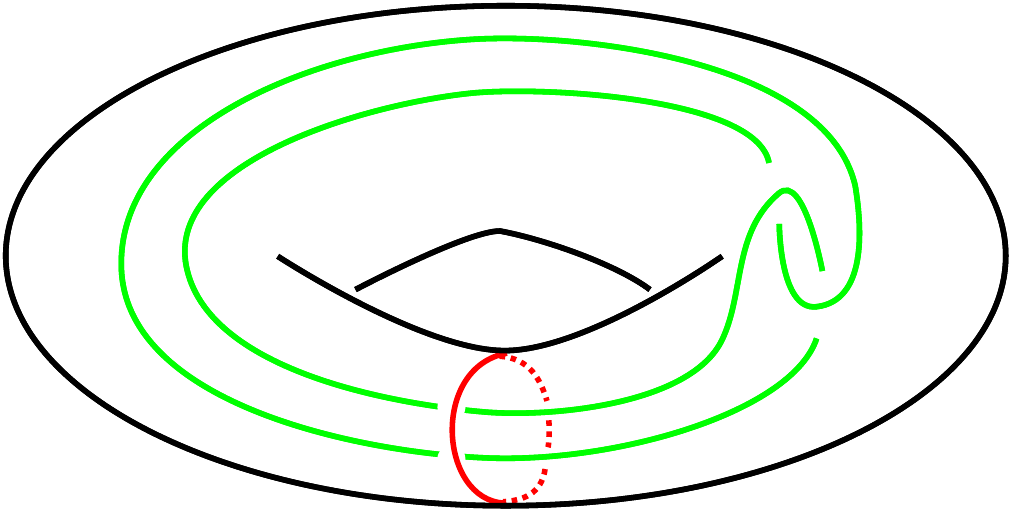}
\caption{$T_2\subset T_1$.}
\label{Fig2}
\end{figure}

For a closed solid torus $T$ we define the notion of a meridian curve. A meridian $\gamma\subset \partial T$ is an~embedded circle which is nullhomotopic in $T$ but not contractible in $\partial T$. On~Figure~\ref{Fig2} the red curve is a meridian.

 Now we choose a closed unknotted solid torus $T_{1}\subset\mathbb{S}^{3}$. It~is well known that the complement of~$T_1$ in~$\mathbb{S}^{3}$ is another solid torus (this time open). We~then embed a second solid torus $T_{2}$ inside~$T_{1}$ as a tubular neighbourhood of the green curve. The~green and red curves form a~Whitehead link. This is the basic pattern of the construction which we will repeat infinitely many times. Precisely, $T_{2}$ is also an unknotted solid torus in $\mathbb{S}^{3}$ and we then embed $T_{3}$ inside $T_{2}$ in the same way as $T_{2}$ lies in $T_{1}$. We~do this infinitely many times and define $T_{\infty}:=\cap_{k=1}^{\infty} T_{k}$ which is a~closed subset of~$\mathbb{S}^{3}$ called the {\it Whitehead continuum}. The~Whitehead manifold is then defined as $Wh:=\mathbb{S}^{3}\setminus T_{\infty}$ which is an open manifold (i.e., a~non compact manifold without boundary), an open subset of~$\mathbb{S}^{3}$ (hence of~$\mathbb{R}^{3}$ too).

The symmetry of the Whitehead link allows to give an alternative definition. Note that each~$T_{k}$ is an unknotted solid torus in $\mathbb{S}^{3}$, it is however knotted in $T_{k-1}$; we recall that this means that no isotopy of~$T_{k-1}$ can unknot $T_{k}$. Nevertheless $T_{k}$ is homotopically trivial in $T_{k-1}$. Being unknotted in $\mathbb{S}^{3}$, the complement of~$T_{k}$ is a (open) solid torus $N_{k}$. The~manifold $\rm Wh$ can then be defined as the increasing union of the open solid tori $\{N_{k}\}_{k\geq 1}$.

This second construction, more flexible, can support variations such as changing the knot at each step $k$, and yields a family of manifolds some of which are not embedded in $\mathbb{S}^{3}$. We~then say that an open $3$-manifold is genus one if it is the increasing union of solid tori $N_k$ so that, for each $k$, $N_k\subset \rm{Int}(N_{k+1})$ and such that a disc filling a meridian of~$N_{k+1}$ intersects the core of~$N_k$. They were introduced in~\cite{McMil} and there are uncountably many of them, some of which are subsets of~$\mathbb{S}^3$, some not (see~\cite{K-Mil}). Note that $\mathbb{R}^{3}$ is not genus one but we
could call it genus zero, since it is an increasing union of~$3$-balls. The~construction can also be made with handlebodies of higher genus~\cite{McMil} and the genus can also change at each stage of the construction, the worst case being when this genus goes to infinity. This produces an incredible zoology of contractible pairwise non homeomorphic $3$-manifolds!

Showing that $\rm Wh$ is contractible is easy with the second (ascending) description. Thanks to a theorem due (again) to Whitehead it suffices to show that all homotopy groups are trivial and this follows from the fact that $N_k$ is homotopically trivial in $N_{k+1}$. To prove that it is not homeomorphic to $\mathbb{R}^{3}$ we simply show that $\rm Wh$ is not simply-connected at infinity. Let us recall the definition of the notion of {\it simply-connectedness at infinity}.

\begin{defi}
A connected, locally compact and simply-connected topological space $M$ is \emph{simply-connected at infinity} if, for any compact set $K\subset M$, there exists a compact set $K'$ containing $K$ so that the induced map $\pi_{1}(M\setminus K')\rightarrow \pi_{1}(M\setminus K)$ is trivial.
\end{defi}

In other words, for any compact set $K$, the loops which are far away from $K$, say in the complement of~$K'$, can be contracted in the complement
of~$K$.

We now want to explore the Riemannian geometry of these spaces, the idea being that among all these contractible $3$-manifolds $\mathbb{R}^{3}$ should be special. The~starting point of our study is the article
\cite{SchoenYau:stable} by R.~Schoen and S.-T.~Yau in which they prove the following theorem.

\begin{theor}[{\cite[Theorem~3]{SchoenYau:stable}}]
Let $M$ be a complete non compact $3$-dimensional manifold with positive Ricci curvature. Then $M$ is diffeomorphic to $\mathbb{R}^{3}$.
\end{theor}

The key idea relies on showing that there are no stable complete minimal surface in $M$. \mbox{Indeed}, its Jacobi operator is related to the Ricci curvature of the ambient space whose positiveness would give a contradiction. Recently this result was extended by G.~Liu in~\cite{Liu} (see~Theorem~2)
who gets the conclusion that a contractible $3$-manifold cannot carry a complete Riemannian metric with non negative Ricci curvature unless it is $\mathbb{R}^{3}$. The~next step brings us to the amazing article by M.~Gromov and B.~Lawson,~\cite{GromovLawson}. Particularly Section~10 which starts by~Theorem~10.2 giving a nice relation between a compact stable minimal surface and the ambient scalar curvature. The~following result is Corollary~10.9 on p.~173.

\begin{theor}[{\cite[Corollary~10.9]{GromovLawson}}]
A complete $3$-manifold of uniformly positive scalar curvature and with finitely generated fundamental group is simply-connected at infinity.
\end{theor}

In particular contractible $3$-manifolds cannot carry a complete metric with uniformly positive scalar curvature unless they are diffeomorphic to $\mathbb{R}^{3}$; indeed, a~result by C.H.~Edwards~\cite{Edw} combined
with the proof of the Poincar\'e conjecture shows that a contractible $3$-manifold which is simply-connected at infinity is homeomorphic to $\mathbb{R}^3$. Let us recall that $\mathbb{R}^3$ does carry a~complete metric of uniformly positive scalar curvature.

Results in the same spirit then appear in~\cite{cwy:taming} where it is proved, in~particular, that

\begin{theor}[{\cite[Theorem~4.4]{cwy:taming}}]\label{cwy}
If a non compact contractible $3$-manifold $M$ has a complete Riemannian metric with uniformly positive scalar curvature outside a compact set, then it is homeomorphic to $\mathbb{R}^{3}$.
\end{theor}

Then, two striking results by Jian Wang gave a definitive answer for certain classes of contractible $3$-manifolds. The~first one is the following.

\begin{theor}[\cite{Wang1}]\label{Wang1}
No contractible genus one $3$-manifold admits a complete metric of non negative scalar curvature.
\end{theor}

The existence of complete metrics of non negative scalar curvature is related to the fundamental group at infinity which is the inverse limit of the fundamental groups of complements of compact subsets. It~is denoted by $\pi_1^{\infty}$. A well-known fact is that any $3$-manifold which is simply-connected at infinity has a trivial $\pi^\infty_1$. However the converse is not true, i.e., a $3$-manifold with trivial $\pi_{1}^{\infty}$ may be not simply-connected at infinity. For example, the Whitehead manifold has trivial $\pi_{1}^{\infty}$ but is not simply-connected at infinity, indeed, in~this case, the homotopy classes of non-trivial loops in
the complement of a compact set $K\subset X$ do not persist when we increase~$K$. We~then get,

\begin{theor}[\cite{Wang2}]\label{Wang2}
A contractible $3$-manifold with non negative scalar curvature and trivial~$\pi_1^{\infty}$ is homeomorphic to $\mathbb{R}^{3}$.
\end{theor}

In~\cite{Wang1, Wang2} the two previous statements are made with the assumption that the scalar curvature is positive. However a nice argument by J.~Kazdan then allows Jian Wang to state the results as above. The~method
of proof pertains to the same philosophy than in~\cite{SchoenYau:stable} and~\cite{GromovLawson}. Let us describe some of the ideas contained in the proof of~Theorem~\ref{Wang1}, in~the case where the scalar curvature is supposed to be positive.

 With the notations used to define genus one $3$-manifolds, let us consider a meridian $\gamma\subset\partial {\bar N}_k$ (where $\bar N_k$ is the
closure of~$N_k)$, this is an embedded curve which can be filled by a minimizing disk $D_k$. Let us assume, for a while, that this disk is included in $\bar N_k$. Jian Wang shows that the number of connected components of~$D_k\cap N_2$ intersecting $N_1$ goes to infinity with $k$. The~fact that $D_k$ is included in $\bar N_k$ is ensured when $\partial \bar N_k$ is mean convex and one can always deform the Riemannian metric in a small neighbourhood of~$\partial \bar N_k$ so that it becomes mean convex. The~disk $D_k$ is then included in $\bar N_k$ and is minimal for this new metric which coincides with the original one in, say, $N_{k-1}\subset N_k$ and this is sufficient for the rest of the argument. Now, let us assume that the sequence $D_k$ converges to a complete minimal surface $\Sigma\subset X$ which, in~that situation, is stable. By~the result of Schoen and Yau (see~\cite{SchoenYau:stable}) this surface is homeomorphic to $\mathbb{R}^2$ and the previous argument shows that the number of connected components of~$\Sigma\cap N_2$ intersecting~$N_1$ is infinite. By~a result of Meeks and Yau (see~\cite{Mee-Yau}), each of these components contains a definite amount of area. The~contradiction comes from a nice extrinsic version of Cohn--Vossen's inequality, proved by Jian Wang,
\[
\int_\Sigma \kappa (x)\,{\rm d} v(x) \leq 2\pi,
\]
where $\kappa (x)$ is the scalar curvature of~$X$ at~$x$ and $dv$ is the volume element of the induced metric on $\Sigma$. The~original Cohn--Vossen's inequality is the same, for a complete surface with positive curvature, where $\kappa$ is replaced by the intrinsic curvature. Since on $N_2$ the scalar curvature of~$X$ is bounded away from zero, by assumption, this inequality is in contradiction with the infinite area contained in~$N_2$. Now, this is by far too naive; indeed, in~general $D_k$ does not converge to a complete stable $\Sigma$ but, according to Colding and Minicozzi
(see~\cite{Col-Min}), to a lamination with complete stable minimal leaves. A variation of the above argument, much more involved, leads to the same contradiction with the extrinsic Cohn--Vossen's inequality. As~mentioned before, in~the case when the scalar curvature is non negative, a~trick
due to J.~Kazdan allows to deform it into a metric with positive scalar curvature.

The proof of Theorem~\ref{Wang2} is a variation on the one sketched above
with non-trivial modifications.

\subsection{Conclusions}

This short account of the Riemannian geometry of open $3$-manifolds is focused on the positive or non negative scalar curvature and on two very specific families of~$3$-dimensional spaces. Despite these restrictions it is quite difficult to prove results and the techniques involved are definitely sophisticated. If~one wants to go beyond Theorem~\ref{Wang3}, considering for example complete metrics with positive scalar curvature (not uniformly positive), one could expect to have more building blocks to consider, not only spherical manifolds. It~is however difficult to make any conjecture at this stage.

On the other hand, concerning Theorems~\ref{Wang1} and~\ref{Wang2} it is quite natural and hopefully safe to state the following conjecture.

\begin{Conjecture}
Let $M$ be a contractible open $3$-manifolds. If~$M$ admits a complete Riemannian metric of non negative scalar curvature then it is diffeomorphic
to $\mathbb{R}^{3}$.
\end{Conjecture}

That would be a beautiful geometric characterisation of~$\mathbb{R}^{3}$ among all contractible $3$-mani\-folds. Notice that, if we allow the manifold to have a non trivial fundamental group, examples carrying a complete hyperbolic metric do exist, see~\cite{Sou-Sto}. However, in~these examples by J.~Souto and M.~Stover the fundamental group is not finitely generated.

We did not address any issue concerning a metric whose curvature is bounded above. It~is clear, thanks to Cartan--Hadamard theorem, that any contractible $3$-manifold cannot carry a~complete metric with non positive sectional curvature unless it is diffeomorphic to $\mathbb{R}^{3}$. It~turns
out that this statement is true even if we relax the regularity. More precisely, D.~Rolfsen proved in~\cite{Rolf} that a complete open $\CAT(0)$ manifold of dimension $3$ is homeomorphic to $\mathbb{R}^{3}$. Hence, Whitehead's type manifolds cannot carry any complete $\CAT(0)$ metric. Notice
that the above result by D.~Rolfsen is not any more true in higher dimension (see~\cite{D-J}).

In recent years we have witnessed a huge activity around synthetic versions of the notion of Ricci curvature bounded below. New families of metric
spaces have been described, such as $\CDCD (\kappa, N)$ spaces (see~Section~\ref{BGexemples}($b$)), and lots of results show that they share plenty of
nice properties with manifolds whose (standard) Ricci curvature is bounded below by $\kappa$ and whose dimension is (bounded above by) $N$. These notions might not see the details of the local geometry, they rather focus on the geometry in the large; however, let us recall that the assumption in Theorem~\ref{cwy} is only at infinity, more precisely outside a compact set. Then, in~view of Liu's result (see~\cite{Liu}), we~are led to ask the following question.

\begin{question}
Does the Whitehead manifold (or any contractible open $3$-manifold not homeo\-morphic to $\mathbb R^3$) carry a geodesically complete $\CDCD (0, 3)$ metric?
\end{question}

This introduces the next section which is devoted to yet another version of synthetic Ricci curvature, even more flexible. It~relies on a weak version of Bishop--Gromov's inequality and serves as an alternative to a lower bound on the Ricci curvature. Using notations defined in the next section we may then ask the following question.

\begin{question}
Does the Whitehead manifold (or any contractible open $3$-manifold not homeo\-morphic to $\mathbb R^3$) carry a geodesically complete metric measure
structure satisfying condition $\BG (0, 3)$?
\end{question}

One difficulty is that most of these spaces do not have any quotient which is a manifold or even an orbifold. We~then loose all the tools that group actions could bring into play.

Now, going to dimension $4$, there is a family of open spaces which plays
a role comparable to the contractible open $3$-manifolds, that is the differentiable structures on $\mathbb R^4$. To our knowledge very little is known about their possible Riemannian geometries and we could dream of a result that would characterise the standard differentiable structure on $\mathbb R^4$ by some of its geometric properties.

What about a synthetic version of the non negative scalar curvature? An attempt to give a non differential definition of scalar curvature was made in~\cite{Vero}, which could lead to a notion of scalar curvature for metric spaces. Some recent works by M.~Gromov seem to pave the way towards such a notion; the interested reader is referred to~\cite{Misha2, Misha1}
and the preprints posted \href{https://www.ihes.fr/~gromov/category/positivescalarcurvature/}{here}. Other interesting articles along these lines are~\cite{ChaoLi2, ChaoLi1}, by Chao Li.

In~\cite{Gro-Large}, M.~Gromov raised the following problem which we state in a simple case. Let $(M,g)$ be a closed Riemannian manifold of dimension $n$ and let us look at its universal cover $\widetilde M$ endowed with the pulled-back metric $\tilde g$. If~one considers the growth function
of~$\big(\widetilde M, \tilde g\big)$, at~some point $x\in\widetilde M$, that is the volume of the ball of radius $R$ around~$x$, denoted by $v(x,R)$, then, when~$R$ goes to zero, this volume is asymptotic to the volume of a ball in $\mathbb{R}^n$, of radius $R$, with a~remainder which involves the scalar curvature at~$x$. Next, when~$R$ goes to $+\infty$ the volume of the ball of radius $R$ has a behaviour driven by the volume entropy of~$\big(\widetilde M, \tilde g\big)$ (see~Definition~\ref{Entropies0}), this raises the following
question:
\begin{question}\label{question:scal-ent}
Could we make precise the influence of the scalar curvature on the entropy?
\end{question}
In 1985, in~\cite[p.~108]{Gro-Large}, Gromov stated what he called ``a vague conjecture'' relating the function $\sup_{x\in\widetilde M}v(x,R)$ of a ``large'' manifold to the growth function of the Euclidean space of~the same dimension. He then described a series of situations in which he gave a more precise meaning to the term ``large''. The~statement of Question~\ref{question:scal-ent} is also vague since we do not
know what to expect. Indeed, the precise shape of the growth function is unknown in general, except for the two asymptotic behaviours mentioned above. Some progresses have been made though since the article~\cite{Gro-Large}, for example a proof of one of Gromov's version of the vague conjecture is given in~\cite{Guth-large}.

Question~\ref{question:scal-ent} is very challenging and exciting and has
always been in our mind.

\section[Bishop--Gromov inequality generalised: on a joint work by G.~Besson, G.~Courtois, S.~Gallot and~A.~Sambusetti]
{Bishop--Gromov inequality generalised:\\
on a joint work by G.~Besson, G.~Courtois, S.~Gallot \\and~A.~Sambusetti}\label{chap:BG}

In the sequel, in~a metric space $(X,d)$, when there is no ambiguity on the choice of the metric $d$, we~shall \textit{denote by} $B_X (x,r)$ (\textit{resp.\ by}
$\overline B_X (x,r)$) the open (resp.\ closed) ball of radius $r$ centred at a point $x \in X$. The~action by isometries of a group $\Gamma$ on $(X,d)$ is said
to be {\em proper} when, for every $R> 0$, the set $\{ \g \in \Gamma\colon d(
x, \g x) \le R\}$ is finite (this does not depend on the choice of~$x \in X$).

\subsection{Synthetic Ricci curvature \`a la Bishop--Gromov}

\subsubsection{The general framework}\label{generalframe}

The notion of ``Ricci curvature bounded from below'', valid on Riemannian manifolds, has many beautiful consequences and the comparison theorems are among them. It~was a revolution when M.~Gromov extended the classical Bishop's inequality to what is now called the Bishop--Gromov's inequality and interpreted Cheeger's finiteness theorem as a consequence of a compactness one. Extending this notion to metric measure spaces is a challenging question. It~is already present in M.~Gromov's works and was made effective by results of J.~Cheeger and T.~Colding~\cite{CC97, CC00a, CC00b}.
Using a compactification of the space of manifolds whose Ricci curvature is bounded from below, they proved that the limit-spaces still verify a Bishop--Gromov's inequality. In~this article, we~shall reverse this circle of
ideas and define a synthetic version of Ricci curvature bounded from below for metric measure spaces satisfying a Bishop--Gromov's inequality.

\medskip\noindent
{\bf (\emph{a}) Bishop--Gromov inequality revisited.}
The notion of ``synthetic Ricci curvature bounded below'', associated to Bishop--Gromov inequalities, is defined by the

\begin{defi}\label{Riccisynthetic0}
For every $N \in{} ]0, + \infty[$ and $K \in{} ]0, + \infty[$, a~metric measure space $(X,d,\mu)$ is said to satisfy condition $BG_x (K,N)$, or equivalently to
 have $N$-dimensional BG-synthetic Ricci curvature $\ge -K^2$, when centred at $x \in X$, if
\begin{equation}
\label{synthetic0}
\forall r \ge \frac{N}{K}\qquad 0 < \mu \big(B_X(x, r)\big) < +\infty \qquad \text{and} \qquad \frac{\mu \big( B_X(x, 2 r)\big)}{\mu \big( B_X(x,r)\big)} \le
2^{N}{\rm e}^{K r}.
\end{equation}
It is said to satisfy condition $BG (K,N)$, or equivalently to have $N$-dimensional BG-synthetic Ricci curvature $\ge -K^2$, if it satisfies condition
$BG_x (K,N)$ for every $x \in X$.
\end{defi}

Notice that the right-hand side of inequality~\eqref{synthetic0} goes to $2^N$ when $r\f 0$, and this is exactly the limit of~$\frac{\mu (B_X(x, 2 r))}{\mu ( B_X(x,r))}$ in the case where $(X,d)$ is a $N$-dimensional Riemannian manifold and $\mu$ its Riemannian measure. However, in~our setting, $N$ is generally not an integer.

The reader may observe that the assumption ``$N$-dimensional BG-synthetic Ricci curvature \mbox{$\ge - \kappa^2$}'' does not imply that the $N$-dimensional
BG-synthetic Ricci curvature is $\ge - K^2$ when $\kappa < K$, because the interval where the inequality~\eqref{synthetic0} is valid is smaller under the first
assumption. Indeed, even in the classical case, under the hypothesis ``Ricci curvature $\ge -K^2$'', the value of~$K$ does not give any geometric or topological information if
not rescaled by ano\-ther geo\-met\-ric invariant; in fact, for any good (synthetic or classical) notion of Ricci curvature bounded below, for any metric space
$(X,d)$, it is equivalent to say that the Ricci curvature of~$(X,d)$ is $\ge - K^2$ and to say that the Ricci curvature of~$(X, K d)$ is $\ge - 1$. When we aim at
bounding topological invariants or scales invariant geometric quantities,
what makes sense is thus to define the notion of~$N$-dimensional BG-synthetic
Ricci curvature $\ge - 1$, and then to define the notion of~$N$-dimensional BG-synthetic Ricci curvature $\ge - K^2$ by homothety.

On the contrary, if one wants to make the lower bound of the radii $r$ independent of the lower bound of the $N$-dimensional
BG-synthetic Ricci curvature, we~introduce the

\begin{defis}\label{generalframe0}
Given parameters $r_0> 0$, $ K \ge 0$ and $ C > 1$, a~metric measure space $(X,d,\mu)$ is said to verify
\begin{itemize}\itemsep=0pt
\item[$(i)$]
a weak Bishop--Gromov inequality at scale $r_0$, centred at $x \in X$, with factor $C$ and exponent $K$ if, for every $r \ge r_0$,
\begin{gather}\label{BGinequality}
 0 < \mu \big(B_X(x, r)\big) < +\infty \qquad \text{and} \qquad \frac{\mu \big( B_X(x, 2 r)\big)}{\mu \big( B_X(x,r)\big)} \le C {\rm e}^{K r},
\end{gather}
\item[$(ii)$]
a weak Bishop--Gromov inequality at scale $r_0$, with factor $C$ and exponent $K$ if it verifies the
corresponding weak Bishop--Gromov inequality at scale $r_0$ centred at every point \mbox{$x \in X$},
\item[$(iii)$]
a strong Bishop--Gromov inequality with factor $C$ and exponent $K$ if~\eqref{BGinequality} is verified for every $r > 0$ and every $x \in X$.
\end{itemize}
\end{defis}

These two notions are in some sense equivalent by the following trivial

\begin{Lemma}\label{Bishopvssynthetic}
A metric measure space whose $N$-dimensional BG-synthetic Ricci curvature
is $\ge -K^2$ when centred at some point $x$ verifies a weak Bishop--Gromov inequality,
centred at~$x$, at~scale $\frac{N}{K}$, with factor $2^{N}$ and exponent $K$.

Conversely a metric measure space, which verifies a weak Bishop--Gromov inequality, centred at some point $x$, at~scale $r_0$, with factor $C$ and
exponent $K$, has $N'$-dimensional BG-synthetic Ricci curvature $\ge -K^2$ when centred at~$x$, where
$N' = \Max \big( K r_0, \frac{\ln C}{\ln 2} \big)$.
\end{Lemma}

Lemma~\ref{Bishopvssynthetic} proves that the two notions: ``$N$-dimensional BG-synthetic Ricci curvature boun\-ded below'' and
``weak Bishop--Gromov inequality satisfied at some scale'' are equivalent only when this scale is sufficiently large. In~the applications
developed in Sections~\ref{systolebound1},~\ref{weakvsstrongBishop} and~\ref{finicompact}, we~have to use a universal upper bound of the quotient
$\frac{\mu ( B_X(x, R))}{\mu ( B_X(x,r))}$ when $r$ is small enough,~and the notion of weak Bishop--Gromov inequality at a
sufficiently small scale is well adapted to this purpose. On the contrary, the notion of~$N$-dimensional BG-synthetic Ricci curvature bounded below
is not adapted to the case where $r$ is small; indeed, on $n$-dimensional
Riemannian manifolds, extending to small balls inequality~\eqref{synthetic0},
would imply simultaneously $2^N \ge 2^n$ and $N < 1$, a~contradiction.

Referring to Definitions~\ref{generalframe0}, it is obvious that condition $(iii)$ implies condition $(ii)$, which implies condition $(i)$. Moreover condition $(iii)$ is much stronger than condition $(ii)$.
Indeed, a~weak Bishop--Gromov inequality at scale $r_0$ (resp.\ the condition
$BG (K, N)$) only concerns balls of radius $r \ge r_0$ (resp.\ of radius $r \ge \frac{N}{K}$) and gives no information on balls of
smaller radius, in~particular it gives no information on the local topology and geometry of the space which verifies these properties. A counter-example may
be constructed by taking the connected sum of a fixed Riemannian manifold
$(X,g)$, satisfying a strong Bishop--Gromov inequality, with any Riemannian manifold $(Y,h)$ of diameter $\le \frac{r_0}{20}$; then, for every choice of~$(Y,h)$, the corresponding manifolds
$(X\# Y, g\# h)$ all verify the same weak Bishop--Gromov inequality, with
the same scale~$r_0$, the same factor $C$ and the same exponent $K$, independently of the choice of~$(Y,h)$, even though the topology of~$Y$, and consequently the topology of a ball of radius $\le \frac{r_0}{10}$ in $(X\# Y, g\# h)$, could be anything (see the proof of the Corollary~3.17 of~\cite{BCGS} for clarifications).
On the contrary, a~strong Bishop--Gromov inequality implies restrictions on the local topology of the metric spaces under consideration
(see~\cite[Theorem~5.1]{CS4}).

Turning back to the difference between a weak Bishop--Gromov inequality centred at every point and a weak Bishop--Gromov inequality centred at only one point, we~have

\begin{prop}
Let $(X, d, \mu)$ be a metric measure space which admits a proper action,
by isometries preserving the measure, of a group $\Gamma$ such that
$\diam (\Gamma \backslash X) \le D$, if there exists a~point $x_0 \in X$ such that $(X, d, \mu)$ verifies a weak Bishop--Gromov inequality, centred
at~$x_0$, at~scale $r_0$, with factor $C$ and exponent $K$, then $(X, d, \mu)$ verifies a weak Bishop--Gromov inequality at~scale $r_0 + \frac{5}{2}\,
D$, with factor $C^2$ and exponent $2 K$ centred at every point $x \in X$.
\end{prop}

\begin{proof} For every point $x \in X$, there exists $\g \in \Gamma$ such that $d(x, \g x_0) \le D$. Hence, for every $r \ge r_0 + \frac{5}{2}\, D$, we~have
$B_X(x,2 r) \subset B_X(\g x_0, 2 r + D)$ and $ B_X(\g x_0, r - D) \subset B_X(x, r)$, and consequently, using the $\Gamma$- invariance of~$d$ and $\mu$, we~get $\mu \big( B_X(x, 2 r)\big) \le \mu \big( B_X( x_0, 2 r + D)\big)$
and $\mu \big( B_X(x, r)\big) \ge \mu \big( B_X( x_0, r - D)\big)$; this
and the fact
that $ r + \frac{D}{2} \le 2\, (r -D)$ yield
\begin{gather*}
\dfrac{\mu \big( B_X(x, 2 r)\big)}{\mu \big( B_X(x,r)\big)} \le \dfrac{\mu \big( B_X( x_0, 2 r + D)\big)}{\mu \big( B_X( x_0, r- D)\big)} \leq
\dfrac{\mu \big( B_X( x_0, 2 r + D)\big)}{\mu \big( B_X( x_0, r+ D/2)\big)} \, \dfrac{\mu \big( B_X( x_0, 2 (r - D))\big)}{\mu \big( B_X( x_0,r- D)\big)}
\\ \hphantom{\dfrac{\mu \big( B_X(x, 2 r)\big)}{\mu \big( B_X(x,r)\big)}}
{}\le C {\rm e}^{K ( r+ D/2)}\, C {\rm e}^{K ( r-D)} < C^2 {\rm e}^{2 K r}.
\tag*{\qed}
\end{gather*}
\renewcommand{\qed}{}
\end{proof}

This proposition does not prove that a weak Bishop--Gromov inequality centred at one point is equivalent to a weak Bishop--Gromov inequality centred
at
every point, not only because it~depends on an upper bound on the diameter but, above all, because, if the first property is valid at a small scale $r_0$, it
only implies that the second property is valid at a much larger scale $r_0 + \frac{5}{2}\, D$. We~shall see that many
of the applications are valid only if we are able to prove that the scale
may be chosen small with respect to the diameter of~$\Gamma \backslash X$.

These Bishop--Gromov inequalities provide the following ones, of a more classical form:

\begin{Lemma}\label{Bishopclassic}
Let $(X,d,\mu)$ be a metric measure space which satisfies a weak Bishop--Gromov ineq\-uality, centred at $x\in X$, at~scale $r_0$, with factor $C$
and exponent $K$ $($resp.\ whose $N$-dimensional BG-synthetic Ricci curvature is bounded below by $-K^2$ when centred at~$x)$ then, for all $r, R$
such that $ r_0 \le r < R$ $($resp.\ such that $\frac{N}{K} \le r < R)$, one has
\begin{gather*}
\dfrac{\mu \big( B_X(x, R)\big)}{\mu \big( B_X(x,r)\big)}
\le C \bigg(\dfrac{R}{r}\bigg)^{\frac{\ln C}{\ln 2}} {\rm e}^{K R}\qquad
\bigg(\text{resp.\ } \dfrac{\mu \big( B_X(x, R)\big)}{\mu \big( B_X(x,r)\big)} \le 2^N \bigg(\dfrac{R}{r}\bigg)^N {\rm e}^{KR}\bigg).
\end{gather*}
\end{Lemma}

\begin{proof} If $(X,d,\mu)$ satisfies a weak Bishop--Gromov inequality, centred at $x\in X$, at~scale $r_0$, with factor $C$ and exponent $K$,
let $p$ be the integer such that $2^p < \frac{R}{r} \le 2^{p+1}$. Inequality~\eqref{BGinequality} gives
\begin{gather*}
\dfrac{\mu \big( B_X(x, R)\big)}{\mu \big( B_X(x,r)\big)} \le \prod_{i=0}^{p-1} \dfrac{\mu \big( B_X(x, 2^{- i}R)\big)}{\mu \big( B_X(x, 2^{-(i+1)}R)\big)}
 \dfrac{\mu \big( B_X(x, 2^{- p}R)\big)}{\mu \big( B_X(x,r)\big)}
\le C^{p+1} {\rm e}^{KR (2^{-1}+ \cdots + 2^{- p})} {\rm e}^{K r}
\\ \hphantom{\dfrac{\mu \big( B_X(x, R)\big)}{\mu \big( B_X(x,r)\big)}}
{}\le C^{p+1} {\rm e}^{K R}.
\end{gather*}

The second inequality of Lemma~\ref{Bishopclassic} follows from the first part of Lemma~\ref{Bishopvssynthetic}.
\end{proof}

\medskip\noindent
{\bf (\emph{b}) Lower bound of the $\boldsymbol N$-dimensional BG-synthetic Ricci curvature versus entropy.}
\begin{defi}\label{Entropies0}
The entropy of a metric measure space $ (X,d, \mu) $ (denoted by $ \Ent (X,d, \mu) $) is the lower limit (when $ R \to +\infty $)
of~$ \frac{1}{R} \ln \big( \mu \big( B_X (x, R)\big) \big) $. It~does not depend on the choice of~$x$.
\end{defi}

This invariant, possibly infinite, gives an asymptotic, hence weak, information on the geometry of the metric space (see Section~\ref{comparaison}), nevertheless it becomes interesting
when there exists a group $ \Gamma $ acting properly by isometries on $ (X,d) $ (and possibly co-compactly) and when we restrict ourselves to Borel measures $ \mu $ which are invariant under this action. When the action of~$\Gamma$ is co-compact, a~particular role is played by the \emph{counting measure} $ \mu_{x}^{\Gamma} $ on the orbit $ \Gamma \cdot x $ of a point $ x $ defined by $ \mu^\Gamma_x = \sum_{\gamma \in \Gamma} \delta_{\gamma x}$, where $\delta_y$ is the Dirac measure at the point $y$.

Notice that, in~the co-compact case, the entropy does not depend on the chosen $ \Gamma$-invariant measure, as~shown by the

\begin{prop}\label{Entropies1}
Let $(X,d)$ be a non compact metric space and $ \Gamma $ be a group acting properly and co-compactly on $(X,d)$
by isometries. For every positive $\Gamma$-invariant measure $\mu $ on~$X$ such that all balls have finite measure, one has $\Ent (X,d, \mu) =\Ent \big(X,d, \mu_{x}^{\Gamma}\big)$ for every $x \in X$.

If, furthermore, $(X,d)$ is a length space, then $\Ent (X,d, \mu)$ is the limit $($when $R \to +\infty)$ of~$ \frac{1}{R} \ln \big(\mu \big(B_X (x, R)\big)\big)$.
\end{prop}

This proposition is classical, a~proof may be found in~\cite[Proposition~3.3]{BCGS}.
Pursuant to this proposition, we~shall sometimes use the notation $ \Ent (X,d) $ instead of~$ \Ent (X,d, \mu) $ for $\Gamma$-invariant measures.

The hypothesis ``entropy bounded from above by $K\,$'' may be viewed
as an asymptotic version of ``Ricci curvature bounded from below by
$- K^2\,$"; indeed it is much weaker than the hypothesis ``$N$-dimensional BG-synthetic Ricci curvature $\ge - K^2\,$'', as~proved by the

\begin{prop}\label{BishopvsEntropy}
A non compact metric measure space whose $N$-dimensional BG-synthetic Ricci curvature is $\ge -K^2$ when centred at one point $x$ $($see Definition~$\ref{Riccisynthetic0})$ has entropy $\le K$.

Conversely, when the metric space is a length space and admits a proper,
co-compact action of some group $\Gamma$ of isometries preserving the measure, then its entropy is the infimum of the values $K>0$ such that its $N$-dimensional BG-synthetic Ricci curvature is $\ge -K^2$ for at least one value of~$N$; equivalently, its entropy is also the infimum of the values $K>0$ such that the space satisfies a weak Bishop--Gromov inequality at
scale $r_0$, with factor $C$ and exponent $K$, for at least one value of~$r_0$ and of~$C$. Notice that these values $N$, $r_0$, $C$ then depend on $K$
and on the space.
\end{prop}

\begin{proof}
Denote by $ (X,d, \mu) $ the metric measure space under consideration. For the sake of~sim\-p\-li\-city, let us define $h_x (R) := \frac{1}{R}
\ln \big(\mu \big( B_X (x, R) \big)\big)$ and $\Ent := \Ent (X,d,\mu)$.

Let us first suppose that the $N$-dimensional BG-synthetic Ricci curvature is $\ge -K^2$ ($K> 0$) when centred at the point $x$. Then, by definition of the entropy as
$\liminf\limits_{R\to +\infty} h_x (R)$, for every $\e > 0$, there exists some $R'_0 = R'_0 (\e, X, d)$ such that, $ h_x(R) > \Ent - \e$ for every $R
\ge R'_0$
and there exists a sequence $(R_i)_{i \in \N^*}$ going to infinity such that $\lim\limits_{i \to +\infty} h_x(R_i) = \Ent $. Hence there exists some
$R_0 \ge R'_0$ such that, for every $R_i \ge R_0$, we~simultaneously have
$h_x (2 R_i) > \Ent - \e$ and $h_x(R_i) < \Ent + \e$.
When $R_i \ge \Max \big(R_0, \frac{N}{K} \big)$ this yields:
\begin{gather*}
\Ent - 3 \e < 2\, h_x(2 R_i) - h_x(R_i) = \dfrac{1}{R_i} \ln \Bigg(\frac{\mu \big(B_X (x, 2 R_i)\big)}{\mu \big( B_X (x, R_i) \big)} \Bigg) \le \dfrac{1}{R_i} \ln \big( 2^N {\rm e}^{K R_i}\big)
\\ \hphantom{\Ent - 3 \e }
{}\le K + \frac{N}{R_i} \ln 2.
\end{gather*}
Taking the limit when $R_i \f + \infty$ and when $\e \f 0$, we~get $\Ent
\le K$. It~follows that $\Ent $ is smaller or equal to the infimum
of the values $K$ such that the $N$-dimensional BG-synthetic Ricci curvature is $\ge -K^2$, when centred at~$x$, for at least one value of~$N$.

To prove the reverse implication, let $K$ be any value $> \Ent $, and choose $\e> 0$ such that
$ 2 (1+\e) (\Ent + \e) - (1-\e) (\Ent - \e) < K$. Fixing any point $x_0 \in X$,
Lemma~\ref{Entropies1} implies that $\Ent = \lim\limits_{R \to +\infty} h_{x_0}(R)$ and, if $\diam (\Gamma \backslash X) \le D$, that there exists
$R_0 = R_0 (\e, D, X, d)$ such that, for every $x \in X$, and every $R \ge R_0$,
\begin{gather*}
(1 - \e) (\Ent - \e) \le \bigg(1 - \frac{D}{R}\bigg) h_{x_0} (R-D) \le h_x (R) \le \bigg(1 + \frac{D}{R}\bigg) h_{x_0} (R+ D)
\\ \hphantom{(1 - \e) (\Ent - \e) }
{}\le (1 + \e)(\Ent +\e).
\end{gather*}
Hence, for every $R\ge R_0$ and every $ x \in X$, we~have
\begin{gather*}
\dfrac{1}{R} \ln\Bigg(\!\dfrac{\mu \big( B_X(x, 2 R)\big)}{\mu \big( B_X(x, R)\big)}\!\Bigg)\! = 2 h_x( 2 R) - h_x(R) \le 2 \big(1+\e\big) \big(\Ent + \e \big) - \big(1-\e\big) \big(\Ent - \e \big)\! < K
\end{gather*}
and, choosing $N := K R_0$, we~obtain that the $N$-dimensional BG-synthetic Ricci curvature is $\ge -K^2$.
This proves that the entropy is the infimum of the values $K$ such that the $N$-dimensional BG-synthetic Ricci curvature is $\ge -K^2$ for at least one value
of~$N$.

By Lemma~\ref{Bishopvssynthetic}, this is equivalent to saying that the entropy is the infimum of the values $K$ such that $ (X,d, \mu) $ satisfies a weak
Bishop--Gromov inequality at scale $r_0$, with factor $C$ and exponent~$K$, for at least one value of~$r_0$ and of~$C$.
\end{proof}

However, we~are not satisfied with the equivalence given by Proposition~\ref{BishopvsEntropy}; indeed, saying that all the elements of a family of metric
spaces have $N$-dimensional BG-synthetic Ricci curvature \mbox{$\ge -K^2$} needs
to control uniformly and simultaneously the two constants $N$ and $K$.
Proposition~\ref{BishopvsEntropy} gives an optimal control of the parameter $K$, at~the price that the dimension~$N$ generally goes to infinity. In~other words, in~order to keep the dimension $N$ bounded, we~may be forced to increase the value of the parameter $K$ with respect to the entropy.
Hence, the meaning of~Proposition~\ref{BishopvsEntropy} is that the assumption ``$N$-dimensional BG-synthetic Ricci curvature $\ge -K^2$'' is
stronger than the hypothesis ``Entropy $\le K$"; indeed, it is \emph{much} stronger, as~proved by Example~\ref{nonBishopbis}, which constructs a
sequence $(X_i,g_i)_{i \in \N^*}$ of non compact Riemannian manifolds whose entropy is bounded by some constant $K_1$ (independent on $i$) though,
for every choice of constants $r_0, C, K$, there are infinitely many of these manifolds $(X_i, g_i)$ which do not verify the weak Bishop--Gromov inequality
at scale $r_0$, with factor $C$ and exponent~$K$.

\medskip\noindent
{\bf (\emph{c}) Doubling conditions.}
In the literature, a~doubling condition (for all balls of every radius) is generally considered as a generalisation of the hypothesis ``nonnegative
Ricci curvature''. The~key notion here will be the notion of doubling condition for balls of radius $r \in \big[ \frac{1}{2}r_0, \frac{5}{2} r_0\big]$. If~the fact
that a doubling condition for balls of radius $r\leq\frac{5}{2} r_0$ may be considered as a generalisation of the hypothesis ``Ricci curvature bounded from below'' (see Section~\ref{BGexemples}($a$)), the fact that the condition concerns balls of radius $r \ge\frac{r_0}{2}$ means that
this doubling property gives no information about the local topology and geometry. A counter-example may be constructed by taking the connected
sum of any fixed Riemannian manifold whose Ricci-curvature is bounded from below with any Riemannian manifold of diameter $\le \frac{r_0}{20}$, the construction of this example is described after Lemma~\ref{Bishopvssynthetic}.

\begin{defis}
Given $r_0 > 0$, $C_0 > 1$, a~metric measure space $(X,d, \mu)$ is said to satisfy the $C_0$-doubling at scale $r_0$, centred at a given point $x \in X$, if
\begin{gather}\label{defdouble}
\forall r \in \bigg[\frac{1}{2} r_0, \frac{5}{2} r_0\bigg], \qquad
0 < \mu \big(B_X(x, r)\big) < +\infty \qquad
\text{and} \qquad\dfrac{\mu \big(B_X(x, 2 r)\big)}{\mu \big(B_X(x,r)\big)} \le C_0.
\end{gather}
If this condition~\eqref{defdouble} is satisfied for all $ x \in X $, we~say that $ (X,d, \mu) $ satisfies the $C_0$-doubling at~scale~$r_0$.
\end{defis}

\medskip\noindent
{\sloppy{\bf (\emph{d}) Doubling condition and weak Bishop--Gromov inequality are equivalent.}
In~each application, we~shall choose one hypothesis between the three
equivalent following assumptions: ``$(X,d, \mu)$ verifies a weak Bishop--Gromov inequality at scale $\frac{r_0}{2}$'' or
``$(X,d, \mu)$ has $N$-dimensional BG-synthetic Ricci curvature $\ge -K^2\,$'' (with $N \le \frac{1}{2} K r_0$) or ``$(X,d, \mu)$ satisfies a
$C_0$-doubling at~scale~$r_0$''. Our choice is determined by our will to get statements and applications as simple as possible. The~first equivalence is given by Lemma~\ref{Bishopvssynthetic}, the second equivalence is proved by~the

}

\begin{prop}\label{croissanceVol}
Let $(X,d, \mu)$ be a metric measure length space
\begin{itemize}\itemsep=0pt
\item[$(i)$]
if it verifies the weak Bishop--Gromov inequality $\frac{\mu (B_X(x, 2 r))}{\mu(B_X(x,r))} \le C\, {\rm e}^{K r}$ for every $r \ge \frac{r_0}{2}$ $($and every $x \in X)$, then it satisfies a $C_0$-doubling at scale $r_0$, where $C_0 = C\, {\rm e}^{\frac{5}{2} K r_0}$,
\item[$(ii)$]
if it satisfies the $C_0$-doubling at scale $r_0$ then, for every $x \in X$ and every $r \ge \frac{r_0}{2}$, $\frac{\mu (B_X (x,2r))}{\mu (B_X (x,r))} \le C_0^5 \, {\rm e}^{\frac{9}{2} \frac{ r}{r_0}\ln C_0} $.
\end{itemize}
\end{prop}

The proof of the part $(i)$ of this proposition is quite trivial. The~proof
of the part $(ii)$ is more involved (see~\cite{BCGS2}).

This result may be compared to~\cite[Lemma~14.6]{HK} (see also~\cite[inequality (21)]{HK}). Remark however that, at~the same scale $r_0$,
Lemma~14.6 of~\cite{HK} assumes the $C_0$-doubling condition to be valid for all balls of radius $r \in \big]0,\frac{5}{2} r_0\big]$, centred at some point $x$, while the main difficulty, in~the proof of our Proposition~\ref{croissanceVol}$(ii)$, is to conclude without assuming any doubling condition on~balls
of radius $r < \frac{r_0}{2}$.

Applying Proposition~\ref{croissanceVol}$(ii)$ (resp.\ Lemma~\ref{Bishopvssynthetic}), a~$C_0$-doubling condition (resp.\ the hypo\-thesis ``$N$-dimensional
BG-synthetic Ricci curvature bounded below'') only implies a weak Bishop--Gromov inequality where the quantity (scale$\,\times \,$exponent) is bounded below.
Hence the equivalence between these three notions given by these two results is not any more valid when this quantity (scale$\,\times \,$exponent)
is small. This is the reason why, when small radii are involved, we~shall
prefer the hypothesis
``$(X,d, \mu)$ verifies a weak Bishop--Gromov inequality at a given
scale''.

\medskip\noindent
{\bf (\emph{e}) Doubling of the counting measure vs packing condition vs Bishop--Gromov for any other measure.}
\begin{defis}\label{packings}
In a metric space $(X,d)$, for every $x \in X$, given $r, R$ such that $0 < r \le \frac{R}{2}$, we~call $r$-packing of~$B_X (x, R)$ any finite family
of disjoint balls of radius $r$ in $B_X (x, R)$ and, for any isometric proper action of a group $\Gamma$
on $(X,d)$, we~call $(\Gamma,r)$-packing of~$B_X (x, R)$ a $r$-packing whose centres belong to the orbit $\Gamma x$. We~denote by ${\rm Pack} (x, r, R)$ (resp.\ by ${\rm Pack}_\Gamma (x, r, R)$) the supremum (resp.\ the maximal) number of balls in a $ r$-packing (resp.\ in a
$(\Gamma,r)$-packing) of~$B_X (x, R)$. We~shall say that $(X,d)$ verifies the packing condition with bound $N_0$
at scale $r_0$ if ${\rm Pack} \big(x, \frac{r_0}{2}, 11\, r_0\big) \le N_0$ for every $x \in X$.
\end{defis}

When the metric space $(X,d)$ under consideration admits a proper action (by isometries) of some group $\Gamma$ then, among all possible $\Gamma$-invariant measures on~$X$, a~particular role is played by the counting measure $\mu_x^\Gamma$ of the orbit $\Gamma \cdot x$ of some point $x$.

Though not always explicitly said, many of our results are valid under the hypothesis that this counting measure verifies a $ C_0 $-doubling,
centred at $ x $, at~scale $ r_0 $. This last hypothesis is weaker
than the two equivalent hypotheses: ``$(X,d, \mu)$ verifies a weak
Bishop--Gromov inequality at scale $\frac{r_0}{2}$'' and
``$(X,d, \mu)$ satisfies a $C_0$-doubling at scale $r_0$'', it is also weaker than the packing hypothesis at scale $r_0$, as~proved by the following

\begin{prop}\label{comparedoublings} Consider a group $\Gamma$ which acts
properly by isometries on a metric space $ (X,d)$,
\begin{itemize}\itemsep=0pt
\item[$(i)$]
if $(X,d)$ satisfies the packing condition with bound $N_0$ at scale $r_0$ then, for every $x \in X$, its counting measure $\mu_{x}^{\Gamma}$
verifies a $N_0$-doubling centred at~$x$ at scale $2 r_0$,
\item[$(ii)$]
if there exists a $\Gamma $-invariant measure $\mu$ on~$X$ which verifies the $ C_0 $-doubling at scale $r_0$ $($centred at every point of~$X)$, then $ (X,d)$ satisfies the packing condition with bound $C_0^{109}$ at scale $r_0$, in~particular ${\rm Pack} (x, r, R) < +\infty$ when $r_0 \le r \le R/2$,
\item[$(iii)$]
for any $x \in X$, if there exists a $\Gamma $-invariant measure $\mu$ on~$X$ which verifies the weak Bishop--Gromov inequality $($centred at $x)$
at scale $r_0/2$ $($with factor $C$ and exponent $K)$, then
the counting measure $\mu_{x}^{\Gamma} $ verifies the following weak Bishop--Gromov inequality $($centred at $x){:}$
\begin{gather*}
\forall r \ge r_0\qquad
\dfrac{\mu_{x}^{\Gamma} \big(\overline B_X(x, 2 r)\big)}{\mu_{x}^{\Gamma} \big(B_X ( x, r)\big)} < C^{7/2} {\rm e}^{\frac{5}{2}K r}.
\end{gather*}
\end{itemize}
\end{prop}

The proof of this proposition is given by the inequalities:
\begin{gather}
\dfrac{\mu_{x}^{\Gamma} \big(\overline B_X(x, R- r)\big)}{\mu_{x}^{\Gamma}\big(B_X ( x, 2 r)\big)}
\le {\rm Pack}_\Gamma (x, r, R) \le \dfrac{\mu \big(B_X(x, R)\big)}{\mu\big(B_X(x,r)\big)},\nonumber
\\
\label{compardoublings1}
{\rm Pack} (x, r, R) \le \sup_{y \in X} \dfrac{\mu \big(B_X(y, 2R)\big)}{\mu\big(B_X(y, r)\big)},
\end{gather}
where the first inequality comes from a tool that we learned in M.~Gromov's book~\cite{Gr1} (see pp.~270 and 291--292, in~particular Exercise (b)), adapted here in order to prove that a maximal
$(\Gamma,r)$-packing $\big(B_X (\g_i x, r)\big)_{i \in I}$ of~$B_X (x, R)$ provides a concentric covering $\big(B_X (\g_i x, 2 r)\big)_{i \in I}$ of~$\overline B_X (x, R- r) \cap \Gamma \,x$; the second inequality in~\eqref{compardoublings1} is due to the fact that
$\mu \big(\cup_{i \in I} B_X (\g_i x, r)\big) = \# I \cdot \mu \left(B_X \big( x, r \big) \right)$. The~last inequality in~\eqref{compardoublings1} is a consequence
of the fact that, for any $r$-packing $\big(B_X(x_i,r)\big)_{i \in I}$ of~$B_X (x,R)$, there exists $i_0 \in I$ such that
$\# I \cdot \mu \big(B_X(x_{i_0}, r)\big) \le \mu \big(B_X(x, R)\big)$ and $\mu \big(B_X(x, R)\big) \le
\mu \big(B_X( x_{i_0}, 2 R)\big)$.

Another consequence of~\eqref{compardoublings1} is that the doubling of the counting measure $\mu_{x}^{\Gamma} $, centred at~$x$, at~scale $2 r_0$
is equivalent (roughly speaking) to the condition ${\rm Pack}_\Gamma \big(x, \frac{r_0}{2}, 11\, r_0\big) \le N_0$ (see~\cite[Lemma~3.11]{BCGS} for a complete proof, for more explanations,
see also~\cite[Lemmas 3.14 and 3.15]{BCGS}). At the same scale, this last
packing condition is generally much weaker than the
classical one, stated in Definitions~\ref{packings}; indeed there are clearly more $r$-packings than $(\Gamma,r)$-packings.

For this reason, we~can guess that, in~Proposition~\ref{comparedoublings}$(i)$, the packing condition at scale $r_0$ is generally \textit{strictly} stronger than the doubling of the counting measure $\mu_{x}^{\Gamma}$ at scale $2 r_0$ centred at~$x$. A proof of this fact is given by Examples~3.13(1), (2) and~(3) of~\cite{BCGS}, where are constructed sequences of pointed Riemannian manifolds $(M_k, x_k)_{k \in \N}$ such that all the $M_k$'s verify the $C_0$-doubling of the counting measure $\mu_{x_k}^{\Gamma}$ at scale $2 r_0$, centred at the point $x_k$, (for values of~$C_0$ and~$r_0$ independent on $k$) while the maximal number of disjoint balls of
radius $ \frac{r_0}{2}$ that can be packed inside a ball of radius $11\, r_0$ of~$M_k$ goes to infinity with $k$.

\subsubsection{Examples of spaces verifying a Bishop--Gromov inequality}\label{BGexemples}

Referring to Definitions~\ref{generalframe0}$(ii)$ and $(iii)$, in~the following list, the first two examples verify a~strong Bishop--Gromov inequality,
and thus a weak Bishop--Gromov inequality at any scale $r_0$; the last two
examples only satisfy a weak Bishop--Gromov inequality at a specified scale $r_0$.

\medskip\noindent
{\bf (\emph{a}) Riemannian manifolds with Ricci curvature bounded below.}
Let us denote by $b_{\kappa, n} (r)$ the volume of the ball of radius $r$ in the simply connected $n$-manifold~$X_{\kappa}^n$ of constant curvature~$\kappa$. The~celebrated original Bishop--Gromov inequality says

\begin{theor}[R.L.~Bishop, M.~Gromov]\label{BishopGromov}
Given $n \in \N^* \setminus \{1\}$ and $\kappa \in \R$, if $(X,g)$ is a~complete $n$-dimensional Riemannian manifold whose Ricci curvature verifies $\Ric \ge (n-1) \kappa \cdot g $ then, for every $r$, $R$ such
that $0 < r \le R < +\infty$, and every $x \in X$, one has $\frac{\Vol_g
 B_X(x, R)}{\Vol_g B_X(x, r)} \le \frac{b_{\kappa, n} (R)}{b_{\kappa, n} ( r)}$.
\end{theor}

Two strengths of this theorem are the facts that this inequality remains valid when the balls cross the cut-locus of the point $x$ and that the equality is
attained when $(X,g) = X_{\kappa}^n$.

When $\kappa \le 0$ and $\Ric \ge (n-1) \kappa \cdot g $, $(X,g)$ verifies a strong Bishop--Gromov inequality, with factor $C = 2^n$ and exponent $(n-1) \sqrt{|\kappa|}$; indeed, this is obvious when $\kappa = 0$ because $\frac{b_{0, n} (2 r)}{b_{ 0, n} ( r)} = 2^n$,
when $\kappa < 0$, it is a consequence of Theorem~\ref{BishopGromov}, which proves that
\begin{gather}\label{BGclassic}
\forall R>0, \quad \forall x \in X \qquad\dfrac{\Vol_g B_X(x, 2 R)}{\Vol_g B_X(x, R)} \le \dfrac{\int_0^{2 R} \sinh ^{n-1} \big(\sqrt {|\kappa|}\, t\big)\, {\rm d}t}{\int_0^R \sinh ^{n-1} \big(\sqrt {|\kappa|}\, t\big)\, {\rm d}t} \le 2^n {\rm e}^{(n-1) \sqrt{|\kappa|} R}.
\end{gather}
From~\eqref{BGclassic}, it is immediate that a complete $n$-dimensional
Riemannian manifold whose Ricci curvature verifies $\Ric \ge (n-1) \kappa \cdot g $
(with $\kappa \le 0$) automatically has $n$-dimensional BG-synthetic Ricci curvature $\ge (n-1)^2 \kappa$.

\medskip\noindent
{\bf (\emph{b}) Metric measure spaces verifying a synthetic condition of Ricci curvature boun\-ded below.}
For metric measure spaces verifying the synthetic notion of Ricci curvature bounded from below ${\rm CD}(K,N)$ of Lott--Sturm--Villani (see~\cite{LV} and~\cite{St}
for definitions and properties), this notion being inspired by the previous curvature-dimension condition of D.~Bakry and M.~Emery, one has the following Bishop--Gromov inequality, where $s_{\kappa} (t) = t$ when $\kappa = 0$, $s_{\kappa} (t) = \sinh (\sqrt{|\kappa|}\, t)$ when $\kappa < 0$ and, when $\kappa > 0$, $s_ {\kappa} (t) = \sin (\sqrt{\kappa}\, t)
$ if $\sqrt{\kappa}\, t \le \pi$ and $s_{\kappa} (t) = 0$ if $\sqrt{\kappa}\, t \ge \pi$,

\begin{theor}[{\cite[Theorem~2.3]{St}}]\label{Bishopsynth}
For every metric measure space $(X,d,\mu)$ which satisfies the curvature-dimension
condition ${\rm CD}(K,N)$ for some numbers $K, N \in \R$ $(N > 1)$ and for every $x$ lying in the support of~$\mu$, one has, for every $0 < r < R$
\begin{gather*}
\dfrac{\mu \big(B_X(x, R)\big)}{\mu \big(B_X(x, r)\big)} \le \dfrac{\int_0^{R} \big(s_{\frac{K}{N-1}} (t)\big)^{N-1}\, {\rm d}t}{\int_0^{r} \big(s_{\frac{K}{N-1}} (t)\big)^{N-1}\, {\rm d}t}.
\end{gather*}
\end{theor}

Under the hypothesis ${\rm CD}(K, N)$, when $K\le 0$, Theorem~\ref{Bishopsynth} and the computations made in~\eqref{BGclassic} imply that
$\frac{\mu (B_X(x, 2\,R))}{\mu (B_X(x, R))} \le 2^N{\rm e}^{\sqrt{(N-1) |K|} R}$, thus that $(X,d,\mu)$ verifies a strong Bishop--Gromov inequality, with factor $C = 2^N$ and exponent $ \sqrt{(N-1) |K|}$. Consequently, a~metric measure space $(X,d,\mu)$ which satisfies the curvature-dimension condition
${\rm CD}(K, N)$, with \mbox{$K < 0$}, automatically satisfies condition ${\rm BG}\big(\sqrt{(N-1) |K|},N\big)$, and its $N$-dimensional BG-synthetic Ricci
curvature is $\ge (N-1) K$.

\medskip\noindent
{\bf (\emph{c}) Metric measure spaces verifying a doubling condition.}
Let us recall (see Proposition~\ref{croissanceVol}) that a metric measure
length space which satisfies the $C_0$-doubling at scale $r_0$ verifies
a weak Bishop--Gromov inequality at scale $\frac{r_0}{2}$ with factor $C_0^5$ and exponent $\frac{9}{2} \frac{ \ln C_0}{r_0}$.

\medskip\noindent
{\bf (\emph{d}) Gromov-hyperbolic metric spaces admitting a cocompact action.}
The notion of~$\delta$-hyperbolic spaces (see Definition~\ref{hypdefinition0}) was introduced by M.~Gromov as a very weak metric notion of negative
curvature: for instance, for a metric space, being $\delta$-hyperbolic (with $\delta = \ln 3$) is a much weaker hypothesis than being ${\rm CAT}(-1)$
(see~\cite[Proposition~1.4.3, p.~12]{CDP}).
Indeed, even if the metric
space is a Riemannian manifold, being $\delta$-hyperbolic gives no information on the topology or on the geometry of balls of radius smaller than $\delta$:
a counter-example may be constructed by taking the connected sum of a $\frac{\delta}{4}$-hyperbolic Riemannian manifold with any Riemannian manifold
of diameter $\le \frac{\delta}{16}$.

It may seem paradoxical that an hypothesis which generalises the notion of ``curvature bounded from above'' would provide a Bishop--Gromov inequality,
while this inequality is usually the consequence of an hypothesis of the type ``(Ricci) curvature bounded from below''. Indeed, here, the hypothesis
``(Ricci) curvature bounded from below'' will be replaced by the assumption ``entropy bounded from above'' (see Definition~\ref{Entropies0}).

We have seen, in~Section~\ref{generalframe}($b$), that the condition ``entropy bounded above'' is much weaker than the condition of satisfying a weak Bishop--Gromov inequality at any scale $r_0$ with any factor~$C$ and any exponent~$K$.
Hence an open question is: {\em on which sets of metric spaces is it possible to prove that the
hypothesis ``entropy bounded above'' implies a Bishop--Gromov inequality at some specified scale $r_0$ and with specified values of~$K$ and $C$?}
A first answer is given by the

\begin{theor} [{\cite[Theorem~5.1]{BCGS}}]\label{cocompact2}
Let $(X, d) $ be any $ \delta$-hyperbolic metric space, for every proper
action
by isometries of a group $\Gamma$ on $(X,d)$ such that the diameter of~$\Gamma \backslash X$ and the
entropy of~$(X, d) $ are respectively bounded from above by $D$ and $K$,
then, for every $x \in X$
\begin{itemize}\itemsep=0pt
\item[$(i)$]
for every $\Gamma$-invariant measure $\mu$ on~$X$, for every $
R > r \ge \frac{5}{2} (7 D + 4 \delta)$,
\begin{gather*}
\dfrac{\mu \big(B_X(x, R)\big)}{\mu \big(B_X(x, r)\big)}
\le 3 {\rm e}^{K D} \bigg(\frac{R}{r} \bigg)^{25/4} \bigg(\frac{R}{r}\bigg)^{6 K D}
{\rm e}^{6 K (R - \frac{4}{5} r)},
\end{gather*}

\item[$(ii)$]
for every $R > r \ge 10\,(D + \delta)$, the counting measure
$\mu^\Gamma_x$ of the orbit $\Gamma x$ verifies the inequalities:
\begin{gather*}
\dfrac{\mu_{x}^{\Gamma} \big(B_X (x, 2r)\big)}{\mu_{x}^{\Gamma} \big(B_X (x, r)\big)}
\le 3^4 {\rm e}^{\frac{13}{2} \,K\, r}\quad
\text{and} \quad
\forall R \ge r \quad
\dfrac{\mu_{x}^{\Gamma}\big(B_X(x, R)\big)}{\mu_{x}^{\Gamma} \big(B_X(x, r)\big)}
< 3 \bigg(\frac{R}{r}\bigg)^{25/4} {\rm e}^{6 K (R- \frac{4}{5} r)}.
\end{gather*}
\end{itemize}
\end{theor}

This proves that, on a metric space $ (X,d)$ (endowed with a co-compact action of a group~$\Gamma$) which satisfies the hypotheses of Theorem~\ref{cocompact2},
every $\Gamma$-invariant measure verifies a weak Bishop--Gromov inequality at scale $\frac{5}{2} (7 D + 4 \delta)$, with factor $C = 2^{8} {\rm e}^{6 K D}$
and exponent $\frac{36}{5} K$; moreover, for every $x \in X$, the counting measure $\mu_{x}^{\Gamma}$ verifies a weak
Bishop--Gromov inequality (for balls centred at~$x$) at scale $10 (D + \delta)$, with factor $C = 3^{4}$ and exponent $\frac{13}{2} \,K$.

Remark that, in~Theorem~\ref{cocompact2}$(ii)$ the inequality itself only depends on the upper bound of the entropy, the upper bounds of the hyperbolicity
constant $\delta$ and of the diameter $D$ entering only in the computation of the scale. A by-product of Theorem~\ref{cocompact2}$(ii)$ and Proposition~\ref{BishopvsEntropy} is thus the following

\begin{Remark}\label{BGequivEnt}
Every Gromov-hyperbolic metric space $(X,d)$ (the value of the hyperbolicity constant being not specified), which admits a proper, co-compact action by
isometries of some group $\Gamma$, verifies the following property: there
exists some scale~$r_0 $ such that the counting measure~$ \mu^\Gamma_x $ of
the orbit of any point $x$ verifies (for balls centred at~$x$) a weak Bishop--Gromov inequality at scale $r_0$, with factor $C = 3^{4}$ and exponent $K$ satisfying
\begin{gather*}
\Ent (X,d) \le K \le \frac{13}{2} \, \Ent (X,d).
\end{gather*}
\end{Remark}

With respect to Proposition~\ref{BishopvsEntropy}, the new property here is that the value of the factor $C$ is fixed.

\begin{questions} In Theorem~\ref{cocompact2}, notice that the bound on the entropy is an absolutely necessary hypothesis, otherwise the Bishop--Gromov
inequality~\eqref{BGinequality} would not be verified for balls of large radius. What about the two other hypotheses? More precisely:
\begin{itemize}\itemsep=0pt
\item[(1)]
Can one get rid of the assumption ``Gromov-hyperbolic'' in Remark~\ref{BGequivEnt}?
\item[(2)]
Can one replace the hypothesis ``$\delta$-hyperbolic'' by
the hypothesis ``Gromov-hyperbolic'' in Theorem~\ref{cocompact2}$(ii)$, i.e., could the scale be estimated independently of
the value of the hyperbolicity constant $\delta$?
\item[(3)]
Is it possible to prove Theorem~\ref{cocompact2}$(i)$ or $(ii)$ in the non co-compact case?
\item[(4)]
Is it possible to prove Theorem~\ref{cocompact2}$(i)$ or $(ii)$ in the co-compact case, independently of~the value
of the upper bound $D$ of the diameter?
\item[(5)]
Is it possible to prove Theorem~\ref{cocompact2} in the ${\rm CAT}(-1)$ case independently of the value of~the upper bound $D$ of the diameter?
\end{itemize}
\end{questions}

The answers to Questions (3) and (4), as~they are written, are negative. Indeed Example~\ref{exemplebis} (resp.\ Examples~\ref{nonBishopbis} and~\ref{nonBishop}) construct $\delta_0$-hyperbolic spaces (resp.\ $\delta_0$-hyperbolic Riemannian manifolds), whose entropy is bounded above by some
constant $K_0$, which admit a co-compact proper action of a group of isometries and on which the counting measure of the orbit of~$\Gamma$ (resp.\ the Riemannian measure) does not verify any weak Bishop--Gromov inequality
at a scale $r_0 $ smaller than the diameter, whatever are the values of the factor
and of the exponent.
However is it possible to prove Theorem~\ref{cocompact2} when replacing the ``bounded diameter'' hypothesis by another one, such as an upper bound
of the Margulis constant (defined at the beginning of Section~\ref{MargLemme0})?

On the contrary, these examples do not answer negatively to Question (1) since, in~this question, the scale is allowed to change when the metric space
changes.

About Question (2), observe that Theorem~\ref{cocompact2} applies to all $\delta$-hyperbolic groups with bounded entropy (see definitions in Section~\ref{hypcase}), but it also applies to all subgroups of~such groups (see inequality~\eqref{actionsubgroup}), even if these subgroups are generally
not
$\delta$-hyperbolic.

\subsubsection{Hierarchy of the various possible hypotheses}\label{comparaison}
We start by a general comment: an upper bound of the entropy, as~well as the weak doubling condition for the counting measure on an orbit of the action of a group $ \Gamma $, makes sense on~gene\-ral metric spaces. On the
other hand, a~lower bound on the Ricci
curvature concerns only Riemannian manifolds.
However, even if we only consider Riemannian manifolds and the case where
$ \Gamma $ is the fundamental group of a closed manifold acting by deck transformations on its universal cover, the comparison between all these conditions is roughly summarised as follows:

\begin{compari}
An upper bound on the entropy is a condition that is strictly weaker than
the condition, for the counting measure of an orbit $ \Gamma x $, to verify a weak
Bishop--Gromov inequality for balls centred at the same point $x$ (see Proposition~\ref{BishopvsEntropy} and the example which follows this proposition), which
 is itself strictly weaker than a packing condition at a similar scale (see Proposition~\ref{comparedoublings}$(i)$ and
the discussion after this proposition), itself strictly weaker than the condition, for the Riemannian measure, to verify a weak Bishop--Gromov inequality for balls centred
at every point (see Proposition~\ref{comparedoublings}$(ii)$ and Example~\ref{nonBishopbis}), itself strictly weaker than the condition, for the Riemannian measure, to verify a strong Bishop--Gromov inequality for balls centred at every point (see the beginning of Section~\ref{weakvsstrongBishop} and Example~\ref{nonBishop}), itself strictly weaker than the curvature-dimension condition ${\rm CD}(K,N)$) (see Theorem~\ref{Bishopsynth}), itself strictly weaker than a lower bound on the Riemannian Ricci curvature.

If we furthermore restrict ourselves to comparing the various weak Bishop--Gromov conditions, the smaller the scale, the stronger the condition (obvious
by Definition~\ref{generalframe0}$(ii)$).
\end{compari}
Proofs, examples and further developments are given in~\cite[Section~3.3]{BCGS}.

\subsection{A first finiteness result}

\subsubsection{A gap result}
A generalisation of M.~Gromov's almost flat theorem (see~\cite{BK}), based on a result of E. Breuillard, B. Green, T. Tao (see~\cite[Corollary~11.2]{BGT} and Theorem~\ref{BGT}) is the

\begin{theor}
Given $N \in \, ]0, +\infty [$, there exists $ \e(N) > 0$ such that the following holds: for every path-connected metric measure length space $(X,d, \mu)$, if there exists
some point $x \in X$ such that the $N$-dimensional BG-synthetic Ricci curvature of~$(X,d, \mu)$, centred at~$x$, is $\ge - 1$, every group $\Gamma$,
acting properly, by isometries preserving the measure, on $(X,d, \mu)$ and verifying $\diam (\Gamma \backslash X) \le \e(N)$, is virtually nilpotent.
\end{theor}

\begin{proof}
Define $C_N := {\rm e}^{\frac{15}{2} N} $, $\e(N) := \frac{N}{\nu (C_N)}$
(where $\nu (\cdot)$ is the function introduced in Theorem~\ref{BGT})
and $R_N := 2 N$.
Let $ S := \{\g \in \Gamma\colon d(x, \g x) \le 2 \e (N) \}$; a result of M.~Gromov
(see~\cite[Proposition~3.22]{Gr1}), whose proof is written for Riemannian manifolds but is still valid on path-connected metric spaces, proves that $S$ is a generating set of~$\Gamma$.
As, by hypothesis, $(X,d, \mu)$ verifies the weak Bishop--Gromov inequality (centred at~$x$) at scale $\frac{R_N}{2}$, with factor $2^N$ and exponent
$1$, Proposition~\ref{comparedoublings}$(iii)$ guarantees that
\begin{gather*}
\dfrac{\mu^\Gamma_x \big(\overline B_X \big(x, 2 R_N\big)\big)}{\mu^\Gamma_x \big(B_X \big(x, R_N\big)\big)} < 2^{\frac{7}{2} N} {\rm e}^{\frac{5}{2}R_N} = 2^{\frac{7}{2} N} {\rm e}^{5 N} < {\rm e}^{\frac{15}{2} N} = C_N.
\end{gather*}
Define $A := \{\g \in \Gamma\colon d(x, \g x) \le R_N\}$, this last Bishop--Gromov inequality implies that
\begin{gather*}
\# (A\cdot A) \le \mu^\Gamma_x \big( \overline B_X \big( x, 2 R_N\big)\big) < C_N \cdot \mu^\Gamma_x \big(B_X \big(x, R_N\big)\big) \le C_N \cdot \# (A),
\end{gather*}
where $A\cdot B := \{\g g \colon \g \in A \text{ and } g \in B\}$. Defining
by iteration $S^p$ as the subset $S^{p-1}\cdot S$ of~$\Gamma$, we~have
$2 \e(N) \cdot \nu (C_N) = 2 N = R_N$, and thus $S^{\nu (C_N)} \subset A$ by the triangle inequality.
Applying Theorem~\ref{BGT}, we~deduce that $S$ generates a virtually nilpotent group, thus that $\Gamma$ is virtually nilpotent.
\end{proof}

\subsubsection{The general case}
M.~Gromov and S.~Gallot proved that, on closed Riemannian manifolds $(M^n, g)$ satisfying Ricci curvature $\ge - (n-1) K^2$ and diameter $\le D$,
the first Betti number $b_1 (M^n)$ is bounded above by~some universal constant, depending only on $n$, $K$ and $D$. The~proof of S.~Gallot being analytic, we~shall follow here the viewpoint
of M.~Gromov, using action of groups, which fits better to our purposes.

The following proposition is essentially due to M.~Gromov (see \cite[Lemma~5.19]{Gr1}).

\begin{prop}\label{gromovarg}
Let $(X,d)$ be a connected length space and $\Gamma$ a group acting properly, by isometries, on $(X,d)$ such that $\diam (\Gamma \backslash X) \le
D$ then, for
each point $x_0 \in X$ and for every $R> 0$, there exists a finite family
$\{\g_i\}_{i \in I}$ of elements of~$\Gamma$ satisfying the following properties:
\begin{itemize}\itemsep=0pt

\item[$(i)$] $ \forall i \in I$ $d(x_0, \g_i x_0) \le 2 D + R$,

\item[$(ii)$] $ \forall i, j \in I$ such that $i\ne j$, $d(\g_i x_0, \g_j x_0) \ge R$,

\item[$(iii)$] the subgroup $\Gamma'$ generated by $\{\g_i\colon i \in I\}$ has
finite index in $\Gamma$.

\end{itemize}
\end {prop}

The proof is similar to M.~Gromov's one and there are no extra difficulties (see~\cite[proof of~Lemma~5.19]{Gr1}).
Let us now settle the main theorem:

\begin{theor}\label{finitelygenerated}
Given any connected length space $(X, d)$ and a group $\Gamma$ acting properly by isometries on $(X, d)$ such that $\diam ( \Gamma \backslash X)\le D$. If~there exists a point $x \in X$ and a $\Gamma$-invariant measure $\mu$ on~$X$ such that $(X,d,\mu)$ has $N$-dimensional BG-synthetic Ricci curvature
$\ge -K^2$, when centred at~$x$, then $\Gamma$ contains a subgroup $\Gamma'$ with finite index, generated by $ N_0 (N, K, D)$ elements, where
$N_0 (N, K, D) = \big[121^N {\rm e}^{3 K D}\big]$.
\end{theor}

\begin{proof} Proposition~\ref{gromovarg} guarantees, for this $x \in X$,
the existence of a finite family
$\{\g_i\}_{i \in I_0}$ of~elements of~$\Gamma$ which generates a subgroup
$\Gamma'$ of finite index in $\Gamma$ and verifies, for every $r_0 > 0$:
\begin{gather*}
\forall i \in I_0\quad d(x, \g_i x) \le 2 D + 2 r_0 \quad \text{and} \quad \forall i, j \in I_0 \quad\text{such that}\quad i\ne j,\quad d(\g_i x, \g_j x) \ge 2 r_0.
\end{gather*}
In this case, $\big(B_X ( \g_i x, r_0)\big)_{i \in I_0}$ is a family of disjoint balls included in $B_X (x, 2 D + 3 r_0)$. We~deduce that
\begin{gather*}
\mu \big( B_X (x, 2 D + 3 r_0)\big) \ge \sum_{i \in I_0} \mu \big(B_X ( \g_i x, r_0)\big) = \# I_0 \cdot \mu \big(B_X ( x, r_0)\big),
\end{gather*}
which implies, choosing $r_0 := \frac{N}{K}$, using the fact that $(X,d,\mu)$ has $N$-dimensional BG-synthetic Ricci curvature $\ge -K^2$
when centred at~$x$ and Lemma~\ref{Bishopclassic}:
\begin{gather*}
\# I_0 \le \dfrac{\mu \left(B_X \big(x, 2 D + 3 \frac{N}{K}\big)\right)}{\mu \left(B_X \big(x, \frac{N}{K}\big)\right)} \le 2^N
\left(\dfrac{2 D + 3 \frac{N}{K}}{\frac{N}{K}}\right)^{N} {\rm e}^{K \left(2 D + 3 \frac{N}{K}\right)} < 121^N {\rm e}^{3 K D}.
\end{gather*}
We end the proof by recalling that $S := \{\g_i \colon i \in I_0\}$ is a generating set of~$\Gamma'$.
\end{proof}

Theorem~\ref{finitelygenerated} (and mutatis mutandis Corollary~\ref{Bettinumbers0}) can be rewritten by replacing the hypothesis
``$(X,d,\mu)$ has $N$-dimensional BG-synthetic Ricci curvature $\ge
-K^2$ when centred at~$x$'' by the hypothesis ``$(X,d,\mu)$ satisfies a weak
Bishop--Gromov inequality, centred at $x\in X$, at~scale $r_0$, with factor $C$ and exponent $K$'', these two hypotheses being equivalent by Lemma~\ref{Bishopvssynthetic}.

\medskip

Now, given a compact, arcwise connected, length space $(X, d)$, which admits a universal covering $\pi\colon \widetilde X \f X$, we~can define the length of a path $c$ in $\widetilde X$ as the length of
the path $\pi \circ c$ in~$(X,d)$ and define the pull-back distance $\tilde d$ on
$\widetilde X$ as the associated length distance. Given a~Borel measure $\mu$ on~$X$, we~define the pull-back measure $\tilde \mu$ on $\widetilde X$ by
\begin{gather*}
\int_{\widetilde X} f (y)\, {\rm d} \tilde \mu (y) := \int_{X} \bigg(\sum_{y \in \pi^{-1}(x)} f (y) \bigg) {\rm d} \mu (x),
\end{gather*}
if the fundamental group $\Gamma$ acts properly, $\Gamma$ is countable and $\pi^{-1}(x)$ also, hence the sum makes sense (eventually infinite)
for every $f\ge 0$. Notice that $\tilde \mu$ is $\Gamma$-invariant by construction.

\begin{Corollary}\label{Bettinumbers0}
Let $(X,d, \mu)$ be any arcwise connected, compact measure length space with diameter $\le D$, which admits a simply connected
universal covering
$ \pi\colon \big(\widetilde X, \tilde d, \tilde \mu\big) \f (X, d, \mu)$, where $\tilde d $ and $\tilde \mu$ are the pull-back distance and measure on $\widetilde X$. If~there exists a point $\tilde x \in \widetilde X$ such that
$\big(\widetilde X, \tilde d, \tilde \mu\big)$ has $N$-dimensional BG-synthetic Ricci curvature $\ge -K^2$ when centred at $\tilde x$, then
$\dim \big( H_1 (X, \R) \big) \le \big[121^N {\rm e}^{3 K D}\big]$.
\end{Corollary}

\begin{proof}
The proof follows the arguments developed by M.~Gromov in~\cite[Section~5.20]{Gr1}. Let $N_0 := \big[121^N {\rm e}^{3 K D}\big]$
for the sake of simplicity; by Hurewicz's theorem, the quotient of~$\Gamma \simeq \pi_1 (X,x)$ by its subgroup of commutators is isomorphic to $H_1 (X, \Z)$. If~$S \subset \Gamma$ generates a~subgroup
of finite index in $\Gamma $, its image under this quotient map generates a subgroup of finite index in~$H_1 (X, \Z)$ and its image in $H_1 (X, \R)$ generates
the real vector space $H_1 (X, \R)$; this implies that $\dim \big( H_1 (X, \R) \big) \le \# S$. As~Theorem~\ref{finitelygenerated}
provides
a subset $S$ of~$\Gamma$ which generates a subgroup of finite index in $\Gamma $ such that $\# S \le N_0$, it yields $\dim \big( H_1 (X, \R) \big)
\le N_0$.
\end{proof}

Theorem~\ref{finitelygenerated} and Corollary~\ref{Bettinumbers0} applies
in particular to the various cases listed in Section~\ref{BGexemples}, i.e., to Riemannian manifolds
with Ricci-curvature bounded below and diameter bounded above (this is the original result of M.~Gromov), to metric measure spaces (with dia\-me\-ter bounded above) verifying the synthetic condition $CD (K, N)$ of Ricci curvature bounded below, to metric measure spaces (with diameter bounded above) whose universal cover verify a~doubling condition at some scale $r_0$, and to quotients (with bounded diameter) of~$\delta$-hyperbolic spaces
with bounded
entropy. This last application is developed in the following subsection.

\subsubsection[The delta-hyperbolic case]
{The $\boldsymbol\delta$-hyperbolic case}\label{hypcase}

The following results are corollaries of the Bishop--Gromov versions of Theorem~\ref{finitelygenerated} and Corollary~\ref{Bettinumbers0} and of the
 Bishop--Gromov inequality established in Theorem~\ref{cocompact2}$(ii)$ and proved in~\cite[Theorem~5.1]{BCGS}:

\begin{Corollary}\label{Bettinumbers}
Let $(X,d)$ be any arcwise connected length space with diameter $\le D$, if it admits a simply connected universal cover $ \big(\widetilde X, \tilde d\big)$, which is $\delta$-hyperbolic
and has entropy bounded above by $K$, then $\dim \big( H_1 (X, \R) \big) \le 3^6 {\rm e}^{13 K (16 D + 15 \delta)}$.
\end{Corollary}
Even if we limit our focus to the set ${\cal R} (\delta, K, D)$ of closed connected Riemannian manifolds with diameter $\le D$ which are quotients, by a discrete
group of isometries, of a complete $\delta$-hyperbolic Riemannian manifold with entropy $\le K$, Example~\ref{infinitetopologies} proves, when $\delta$, $K$, $D$
are not too small,\footnote{The case where $K < \frac{\ln 2}{27 \delta + 10 D}$ is studied in~\cite[Section~5.2.2]{BCGS}.} that the local topology of such a
manifold can be anything and that there is an infinite number of topologies in ${\cal R} (\delta, K, D)$, as~those obtained by connected sum of an element of~${\cal R} (\delta, K, D)$ with any other closed manifold of
the same dimension. We~wish to emphasize that this example is not in contradiction with Corollary~\ref{Bettinumbers}, thanks to the extra hypothesis, assumed in this
corollary, that the covering space with bounded entropy is simply connected.

{\sloppy
Remark that, under the hypotheses of Corollary~\ref{Bettinumbers}, one cannot expect to bound the other Betti numbers, as~proved by Example~\ref{infiniteBetti}. Indeed, for every integer $k \in [2, n-2]$, this example constructs a sequence $(M_i, g_i)_{i \in \N} $ of Riemannian
manifolds such that $\dim \big( H_k (M_i, \R) \big) \f + \infty$ when $i \f +\infty$, though there exists constants $\delta, D, K$ (independent on $i$) such that each $(M_i, g_i)$ has diameter bounded by $ D$ and $\delta$-hyperbolic universal cover $\big(\widetilde M_i, \tilde g_i\big) $ satisfying
\mbox{$ \Ent \big(\widetilde M_i, \tilde g_i\big) \le K$}.

}

This validates the following
\begin{question}
In Corollaries~\ref{Bettinumbers0} and~\ref{Bettinumbers}, what extra hypothesis can be added in order to obtain upper bounds of
the other Betti numbers\,? to obtain a finite number of homotopy types\,?
\end{question}
First answers to this question are given in~\cite{BCGS2} and in Section~\ref{finicompact}.

Given a \emph{marked group} $(\Gamma, \Sigma)$, i.e. a finitely generated group $\Gamma$ and a finite generating set~$\Sigma$ of~$\Gamma$, its Cayley
graph ${\cal G} (\Gamma, \Sigma)$ is endowed with the distance $(g, \g) \mapsto d_\Sigma (g, \g)$, which is the minimal length of a path (in the graph) connecting
$g$ to $\g$ (all the edges of the graph being of~length $1$).
Classically, $ (\Gamma, \Sigma)$ is said to be $\delta$-hyperbolic if $\big({\cal G} (\Gamma, \Sigma), d_\Sigma \big)$ is a $\delta$-hyperbolic metric space, the entropy of
this metric space is called algebraic entropy of~$ (\Gamma, \Sigma)$ and is denoted by $\Ent (\Gamma, \Sigma)$.

\begin{Corollary}
Let $\Gamma$ be a group which admits a finite generating set $\Sigma$ such that $ (\Gamma, \Sigma)$ is $\delta$-hy\-per\-bolic with entropy $\le K$, then
$\Gamma$ is commensurable to a subgroup generated by $\big[ 3^6 {\rm e}^{13 K
(15 \delta + 16)}\big]$ elements.
\end{Corollary}

Let ${\rm Hyp} (\delta, K)$ be the set of groups $\Gamma$ which admit a finite generating set, say $\Sigma$, such that $ (\Gamma, \Sigma)$ is
$\delta$-hyperbolic with entropy $\le K$; one cannot hope the set ${\rm Hyp} (\delta, K)$ to contain only a~finite number of isomorphisms classes: a trivial counter-example is obtained by considering the products $(\Gamma \times G_i)_{i \in \N}$ of any group $\Gamma \in {\rm Hyp} (\delta, K)$ with a sequence $(G_i)_{i \in \N}$ of non isomorphic finite groups whose Cayley graphs all have bounded diameters.\footnote{For example, every
finite group admits a generating
set such that the corresponding Cayley graph has diameter~$1$: it is sufficient to
take the whole group as a generating set.} It~is then easy to verify the
existence of some~$\delta'$ (independent on~$i$) such that all the $\Gamma \times G_i$ belong to ${\rm Hyp} (\delta', K)$.

To our knowledge, the following question is still open:
\begin{question}
Does there exist a finite set $\{\Gamma_1, \dots, \Gamma_k \}$ of groups such that every element of~${\rm Hyp} (\delta, K)$ is commensurable to one
of the $\Gamma_i$'s\,?
\end{question}

\subsection{Margulis lemmas under weak Bishop--Gromov condition}\label{MargLemme0}

For any proper isometric action of a group $\Gamma$ on a metric space $(X,d)$, for any $x \in X$ and $r \ge 0$ one defines the set
$\Sigma_r (x)$ of the elements $\g$ of~$\Gamma$ such that $d(x, \g x) \le r$ and the subgroup $\Gamma_r (x)$ generated by $\Sigma_r (x)$. We~call ``\emph{Margulis constant}'' of this action, denoted by ${\rm Marg}_\Gamma (X,d)$, the supremum of the values $r \ge 0$ such that $\Gamma_r (x)$
is virtually nilpotent for every $x \in X$; when $(M,g)$ is a Riemannian manifold whose Riemannian universal covering is $\big(\widetilde M, \tilde g\big)$, we~define the
``\emph{Margulis constant}'' ${\rm Marg}(M,g)$ of~$(M,g)$ as the Margulis constant ${\rm Marg}_{\Gamma_M} \big(\widetilde M, d_{\tilde g}\big)$ of the
canonical action of the fundamental group $\Gamma_M$ of~$M$ on $\widetilde M$, endowed with the Riemannian metric~$d_{\tilde g}$.

\subsubsection{A first Margulis lemma}\label{MargLemme}

{\bf (\emph{a}) When a weak Bishop--Gromov inequality is assumed:}
The celebrated Margulis lemma~(\cite{Mar} and~\cite[Section 37.3]{BZ}) proves (in particular) that there exists a universal constant $\e_0 (n, K)$ such that
${\rm Marg}(M,g) \ge \e_0 (n, K)$ for every $n$-dimensional Riemannian manifold $(M,g)$ whose sectional curvature $\sigma$ verifies $-K^2 \le \sigma \le 0$.

The origin of this problem goes back to Bieberbach's theorem~\cite{Bi} (when $\Gamma$ is a discrete group of isometries of the Euclidean space).
Intermediate results in the same direction were obtained by H.~Zassenhaus~\cite{Zass} (in the case where $\Gamma$ is a discrete subgroup of a Lie group),
 by D.~Kazhdan and G.~Margulis~\cite{KM} (in the case where $\Gamma$ acts on a symmetric space) and by E.~Heintze~\cite{zbMATH01825238} in the
negatively curved case. For a more complete study of the background of this problem, see~\cite{Cou}.

{\sloppy
After a short sketch of proof by J.~Cheeger and T.~Colding (see~\cite{CC96}), V.~Kapovitch and~B.~Wilking~\cite{KW} recently proved the existence of
a universal constant $\e_0 (n, K)$ such that \mbox{${\rm Marg}(M,g) \ge \e_0 (n, K)$} for every $n$-dimensional Riemannian manifold $(M,g)$ whose Ricci curvature
is bounded from below by $-(n-1) K^2$ (see also~\cite{Cou}).

}

We generalised these results,\footnote{However, the generalisation is not complete; indeed several previous results also give an upper bound of the
index of the nilpotent subgroup of the group $\Gamma_r (x)$ when this one
is virtually nilpotent.} proving the following one, which improves a previous one by E.~Breuillard, B.~Green and T.~Tao (see~\cite{BGT}).

\begin{prop}[{\cite[Corollary~5.20]{BCGS}}]\label{Marg1}
For every $r_0 >0$, $N, C > 1 $ and $K\ge 0$, there exist constants
$N_1 (N)$ and $N_0 (C, Kr_0)$ such that the following property holds: given
a~metric space $(X,d)$ and a proper action by isometries of a~group $\Gamma$ on this space,
\begin{itemize}\itemsep=0pt
\item[$(i)$] if there exists some $\Gamma$-invariant measure $\mu$ on~$X$ such that $(X,d,\mu)$ has $N$-dimensional BG-synthetic Ricci curvature $\ge -K^2$, when centred at~$x$, then $\Gamma_r (x)$ is virtually nilpotent for every $x \in X$ and every $r \le \frac{N}{N_1 (N)}\, \frac{1}{K}$,
\item[$(ii)$] if there exists some $\Gamma$-invariant measure $\mu$ on~$X$ which verifies a weak Bishop--Gromov inequality, centred at $x\in X$, at~scale
$r_0$ with factor $C$ and exponent $K$, then $\Gamma_{r} (x)$ is virtually nilpotent for every $r \le \frac{r_0}{ N_0 (C, Kr_0)}$.
\end{itemize}
\end{prop}

We obtained a first partial version of this result as a corollary of a result of M.~Gromov, which appeared as a short remark at p.~71 of his celebrated paper~\cite{Gr3}
about the equivalence between virtual nilpotency and polynomial growth: indeed, this remark provided a quantitative version of this equivalence which
proved that, if the constants involved in the polynomial growth are controlled for all balls of radius $\le R_0$, then the group is virtually nilpotent.

Several authors wrote quantitative versions of this theorem of M.~Gromov (among them: B.~Kleiner, Y.~Shalom and T.~Tao, E.~Hrushovski). Among them, we~have chosen the following result by E.~Breuillard, B.~Greene and T.~Tao, which is convenient to our purposes:

\begin{theor}[{\cite[Corollary~11.2]{BGT}}]\label{BGT}
For every $ C > 1$, there exists $ \nu = \nu (C) \in \N^*$ such that the following
property holds for every group $ G$ and any finite symmetric\footnote{A generating set $S$ is said to be symmetric if, for each $s \in S$,
one also has $s^{-1} \in S$.} generating set $S$ of~$G$: if there exists
some $A \subset G$ which contains $S^{\nu (C)}$ and satisfies $\# (A\cdot
A) \le C \,\# (A)$, then
$G$ is virtually nilpotent.
\end{theor}

This theorem is the main tool in the proof of Proposition~\ref{Marg1}, it
also gives the estimate $N_1 (N) = \frac{1}{2} \,\nu \left(300^{3N}\right)$
and $N_0 (C, K r_0) = \frac{1}{2} \,\nu \left(C^3 {\rm e}^{15 K r_0} + 1\right)$ of the
universal constants involved in Proposition~\ref{Marg1}; unfortunately the function $C \mapsto \nu (C)$ is not explicit. It~is an open problem to decide in
which cases this function can be made explicit.

\medskip\noindent
{\bf (\emph{b}) When the entropy is bounded.}
With respect to Proposition~\ref{Marg1}, the idea here is to~weaken the hypothesis verified by the metric space $(X,d)$ and by the action of the group~$\Gamma$ on~this space, i.e., to replace the hypothesis ``there exists some $\Gamma$-invariant measure $\mu$ on~$X$ which verifies a weak Bishop--Gromov inequality, at~scale $r_0$ (with factor $C$ and exponent $K$)"
by the (much weaker) hypothesis ``there exists some
$\Gamma$-invariant measure $\mu$ on~$X$ such that the entropy of~$(X,d, \mu)$ is bounded above''.

However, this weakened assumption on the metric space is paid by more restrictive hypotheses on the group $\Gamma$ itself, for example we can suppose
that the group $\Gamma$ belongs to the set $\text{\rm Hyp}_{\rm sub} (\delta_0, H_0)$ defined in the following way:

\begin{defi}
Given any real parameters $\delta_0, H_0 > 0$, we~denote by $\text{\rm Hyp} (\delta_0, H_0)$ the set of non virtually cyclic groups $G$ which admit a
finite generating set $S$ such that $(G, S)$ is $\delta_0$-hyperbolic and verifies $\Ent (G, S) \le H_0 $; we denote by $\text{\rm Hyp}_{\rm sub}
(\delta_0, H_0)$ the set of all non virtually cyclic subgroups $\Gamma$ of groups $G$ belonging to $\text{\rm Hyp}(\delta_0, H_0)$.
\end{defi}

\begin{theor}[{\cite[Theorem~6.13(ii)]{BCGS}}]\label{transportnil}
Given $\delta_0, H_0 > 0$, a~proper action by isometries of a group $\Gamma \in \text{\rm Hyp}_{\rm sub} (\delta_0, H_0)$ on
a metric space $ (X, d) $ verifies $\text{\rm Marg}_\Gamma (X, d) \cdot
\Ent (X,d, \mu) \ge \alpha_0 (\delta_0, H_0) > 0$ for every $\Gamma$-invariant
measure $\mu$, where $\alpha (\delta_0, H_0)$ is a positive universal constant which only depends on $\delta_0$ and $H_0$.
\end{theor}
Notice that there is no restriction on the metric space on which this property is valid and that the entropy of~$(X,d, \mu)$ is only used to rescale the Margulis
constant in order that it becomes invariant by the homotheties.
See~\cite[Theorem~6.13]{BCGS} for proofs and other statements.

\subsubsection{A lower bound of the diastole}\label{diastolebound}

{\bf (\emph{a}) The classes of groups which are considered here:}
\begin{defi}
We denote by $\text{\rm Hyp}_{\rm act}$ the set of non virtually nilpotent groups $\Gamma$ which admit
a proper action by isometries on some connected, non elementary Gromov-hyperbolic metric space.
\end{defi}

\begin{defi}
We denote by $\text{\rm Hyp}_{\rm sub}^*$ the set of non virtually cyclic
groups which are isomorphic to some subgroup of some finitely
generated Gromov-hyperbolic group.
\end{defi}

Notice that $\text{\rm Hyp}_{\rm sub}^* \subset \text{\rm Hyp}_{\rm act}$
because, if $G$ is any Gromov-hyperbolic group (endowed with a generating
set
$S$) and $\Gamma$ any subgroup of~$G$ one can consider the action of~$\Gamma$ on the Cayley graph of~$(G,S)$, which is a Gromov-hyperbolic metric space. We~moreover use the fact that, as~the action of~$G$ on its Cayley graph is co-compact, any virtually nilpotent subgroup of~$G$ is virtually
cyclic (this is a corollary of~\cite[Section~8, Th\'eor\`eme 37, p.~157]{GH}).

\medskip\noindent
{\bf (\emph{b}) Systole, diastole and $\boldsymbol r$-thin subsets.}

\begin{defis}
For every proper action by isometries of a group $\Gamma$ on a metric space $(X,d)$,

\begin{itemize}\itemsep=0pt
\item at any point $x \in X$, $\sys_\Gamma (x)$ (resp.\ $\sys^{\diamond}_\Gamma (x)$) is the minimum of~$d(x, \g x)$ when $\g$ runs through the elements of~$\Gamma^*$ (resp.\ through the torsion-free elements of~$\Gamma^*$),

\item the diastole $\dias (\Gamma \backslash X)$ (resp.\ the torsion-free diastole $\dias^{\diamond} (\Gamma \backslash X)$) of~$\Gamma \backslash X$ is the supremum\footnote{The diastole and the torsion-free diastole are both correctly defined on $\Gamma \backslash X$ for $\sys_\Gamma (\cdot)$ and
$\sys^{\diamond}_\Gamma (\cdot)$ are $\Gamma$-invariant functions.} of~$\sys_\Gamma (x)$ (resp.\ of~$\sys^{\diamond}_\Gamma (x)$)
when $x$ runs in $X$,

\item the systole (resp.\ the torsion-free systole) of~$\Gamma \backslash X$ is the infimum of~$\sys_\Gamma (x)$ (resp.\ of~$\sys^{\diamond}_\Gamma (x)$)
when $x$ runs in $X$,

\item the $ r$-thin subset (resp.\ the torsion-free $ r$-thin subset) of~$X$ is the open set $ X_r$ (resp.\ $ X^{\diamond}_r $) of the points $x \in
X$ such that
$\sys_\Gamma (x) < r$ (resp.\ such that $\sys^{\diamond}_\Gamma (x) < r$).
\end{itemize}
\end{defis}

\medskip\noindent
{\bf (\emph{c}) Lower bounds of the diastole.}

\begin{theor}\label{transyst0}
Given a metric space $(X,d)$ and a proper action by isometries of a group $\Gamma \in \text{\rm Hyp}_{\rm act}$ on this space
such that ${\rm Marg}_\Gamma (X,d) > 0$, then $X^{\diamond}_{r}$ is empty or disconnected for every $r < {\rm Marg}_\Gamma (X,d)$.
Consequently, the torsion-free diastole verifies $\dias^{\diamond} (\Gamma \backslash X) \ge {\rm Marg}_\Gamma (X,d)$.
\end{theor}

As $N_0 (C, Kr_0)= \frac{1}{2} \,\nu \big(C^3 {\rm e}^{15 K r_0} + 1\big)$ is the universal constant occurring in Proposition~\ref{Marg1}$(ii)$, we~have the

\begin{Corollary}\label{transyst1}
For every $r_0 >0$, $C > 1 $ and $K\ge 0$, given a metric space $(X,d)$ and a~proper action by isometries of a group
$\Gamma \in \text{\rm Hyp}_{\rm act}$ on this space, if there exists some
$\Gamma$-invariant measure~$\mu$ on~$X$ which verifies a weak Bishop--Gromov inequality at scale $r_0$ with factor~$C$ and exponent~$K$, then $X^{\diamond}_{r}$ is empty or disconnected for every $r \le \frac{r_0}{N_0 (C, Kr_0)}$. Hence
the torsion-free diastole verifies $\dias^{\diamond} (\Gamma \backslash X) \ge \frac{r_0}{N_0 (C, Kr_0)}$.
\end{Corollary}

\begin{proof} It follows from Proposition~\ref{Marg1}$(ii)$ that, for every $x \in X$ and every $r \le \frac{r_0}{N_0 (C, Kr_0)}$, $\Gamma_{r} (x)$ is virtually
nilpotent; Theorem~\ref{transyst0} then implies that $X^{\diamond}_{r}$ is disconnected and that the torsion-free diastole of~$(X,d)$ is larger than
$\frac{r_0}{N_0 (C, Kr_0)} $.
\end{proof}

Weakening the assumption on the metric space and strengthening the hypotheses on the group in the same spirit as in Section~\ref{MargLemme}$(b)$, we~obtain:

\begin{Corollary}[{\cite[Theorem~6.15(ii)]{BCGS}}]\label{transyst2}
Given $\delta_0, H_0 > 0$, for any proper action by isometries of a group
$\Gamma \in \text{\rm Hyp}_{\rm sub} (\delta_0, H_0)$ on
a metric space $ (X, d) $, its torsion-free diastole verifies $\dias^{\diamond} (\Gamma \backslash X) \Ent (X, d, \mu) \ge \alpha_0 (\delta_0, H_0)$
for every $\Gamma$-invariant measure $\mu$, where $\alpha (\delta_0, H_0)$ is the constant of Theorem~$\ref{transportnil}$.
Moreover $X^{\diamond}_{r}$ is empty or disconnected for every $r \le \frac{\alpha_0 (\delta_0, H_0)}{\Ent (Y,d, \mu)}$.
\end{Corollary}
This corollary is an immediate consequence of Theorems~\ref{transyst0} and~\ref{transportnil}; see~\cite[Section~6.3]{BCGS} for other statements
and proofs.

\subsubsection{A first lower bound on the systole}\label{systolebound1}

\begin{theor}\label{minorsystole}
For every $r_0, D >0$, $C > 1 $ and $K\ge 0$, given a metric space $(X,d)$ and a~proper co-compact action by isometries of a~group
$\Gamma $ on this space such that the diameter of~$\Gamma \backslash X$ is bounded above by $D$, if there exists some point $y\in X$ such that $\sys^{\diamond}_\Gamma (y) \ge r_0$ and if there exists some
$\Gamma$-invariant measure $\mu$ on~$X$ which verifies a weak Bishop--Gromov inequality at scale $\frac{r_0}{2}$, with factor $C$ and exponent $K$,
then every torsion-free element $\g$ of~$\Gamma$ verifies
\begin{gather*}
\inf_{x \in X} d(x, \g x) \ge \inf_{x \in X} \sys^{\diamond}_\Gamma (x)
\ge \frac{D}{C} \bigg( 1 + 6\, \frac{D}{r_0}\bigg)^{- \frac{\ln C}{\ln 2}} {\rm e}^{- K\left (3 D + \frac{r_0}{2}\right)}.
\end{gather*}
\end{theor}

\begin{proof}[A sketch of the proof]
As the measure $\mu$ is $G$-invariant for any subgroup $G$ of~$\Gamma$, inequalities~\eqref{compardoublings1}, applied to G-invariant measures, yield the following comparison with the counting measure $\mu_{y}^{G} := \sum_{\g \in G}\, \delta_{g y}$, where $\delta_y$ is the Dirac measure at
the point $y$:
\begin{gather}\label{actionsubgroup}
\forall r > 0\quad \forall R\ge r \quad \forall y \in X \qquad
\dfrac{\mu_{y}^{G} \big(\overline B_X(y, R)\big)}{\mu_{y}^{G}\big(B_X(y, r)\big)}
\le \frac{\mu \big(B_X \big(y, R + \frac{r}{2}\big)\big)}{\mu \big(B_X \big(y, \frac{r}{2}\big)\big)}.
\end{gather}
At a point $y$ where $\sys^{\diamond}_\Gamma (y) \ge r_0$, any torsion-free $\g \in \Gamma^*$ verifies
$d\big(g y, \g^k g y\big) = d \big(y, g^{-1} \g^k g y\big) \ge r_0$ for every $k \in \Z^*$ and every $g \in \Gamma$. If~$G$ is the subgroup generated by $\g$, this
implies that $\mu_{g y}^{G} \big(B_X ( g y, r_0)\big) = 1$ and as, for every $x \in X$, there exists $g \in \Gamma$ such that $d(x, g y) \le D$,
$\mu_{gy}^{G} \big(\overline B_X(g y, 3 D)\big) \ge \mu_{x}^{G} \big(\overline B_X(x, D )\big) \ge 2\, \big[\frac{D}{d(x, \g x)}\big] + 1$, this and inequality~\eqref{actionsubgroup} yield
\begin{gather}\label{minorsystoleeq}
2\bigg[\dfrac{D}{d(x, \g x)}\bigg] \!+ 1 \le \dfrac{\mu_{gy}^{G} \big(\overline B_X(g y,3 D\big))} {\mu_{g y}^{G} \big(B_X (g y, r_0)\big)}
\le \dfrac{\mu \big(B_X \big( g y,3 D\! +\! \frac{r_0}{2}\big)\big)}{\mu \big(B_X(gy, \frac{r_0}{2})\big)} =\dfrac{\mu \big(B_X \big( y, 3 D \!+\! \frac{r_0}{2}\big)\big)}{\mu \big(B_X (y, \frac{r_0}{2})\big)}.\!\!
\end{gather}

We may now apply the weak Bishop--Gromov inequality at scale $\frac{r_0}{2}$ and Lemma~\ref{Bishopclassic} gives
\begin{gather*}
2 \bigg[ \dfrac{D}{d(x, \g x)}\bigg] + 1 \le \dfrac{\mu \big(B_X \big(y, 3D + \frac{r_0}{2}\big)\big)}{\mu \big( B_X ( y, \frac{r_0}{2})\big)}
\le C\bigg(1 + 6\,\frac{D}{r_0}\bigg)^{\frac{\ln C}{\ln 2}} {\rm e}^{K\left(3D + \frac{r_0}{2}\right)}.
\end{gather*}
A direct computation ends the proof.
\end{proof}

A direct corollary of Theorem~\ref{minorsystole} and of the lower bound of the torsion-free diastole given by Corollary~\ref{transyst1} is the following corollary,
where $ \nu (\cdot)$ is the function defined in Theorem~\ref{BGT}:

\begin{Corollary}\label{minorsystole1}
For every $D >0$, $C > 1 $ and $K\ge 0$ we define $\e_1 = \e_1 (K, C, D) := \frac{D}{\nu \left(C^3 {\rm e}^{15 K D} + 1\right)}$. Given a metric space $(X,d)$ and a proper co-compact action by isometries of a group $\Gamma \in \text{\rm Hyp}_{\rm act}$
on this space such that the diameter of~$\Gamma \backslash X$ is bounded above by $D$, if there exists some $\Gamma$-invariant measure $\mu$ on~$X$ which verifies a weak Bishop--Gromov inequality at scale $\e_1$, with factor $C$ and exponent $K$, then every torsion-free element $\g$ of~$\Gamma$ verifies
\begin{gather*}
\inf_{x \in X} d(x, \g x) \ge \inf_{x \in X} \sys^{\diamond}_\Gamma (x)
\ge \frac{D}{C} \bigg(1 + 3\, \frac{D}{\e_1}\bigg)^{- \frac{\ln C}{\ln 2}} {\rm e}^{- K (3 D + \e_1)}.
\end{gather*}
\end{Corollary}

In most of the cases, when we prove a weak Bishop--Gromov inequality, its scale is rather large with respect to the diameter (see Sections~\ref{BGexemples}$(c)$ and $(d)$). This is sufficient to obtain a~bound from below of the torsion-free diastole (see Corollaries~\ref{transyst1} and~\ref{transyst2}), but not to obtain a bound from below of the torsion-free systole; indeed such a lower bound of the systole (given by Corollary~\ref{minorsystole1}) is proved only
if the weak Bishop--Gromov inequality is valid at a scale $\e_1$ small with respect to the diameter.

A question is thus: what are the extra hypotheses which allow to deduce a
weak Bishop--Gromov inequality at a small scale from a weak Bishop--Gromov inequality
at a large scale? A~first answer to this question is given in the following section.

\subsection{From weak to stronger Bishop--Gromov inequalities}\label{weakvsstrongBishop}

It is clear, when comparing Definitions~\ref{generalframe0}$(ii)$ and $(iii)$ that a strong Bishop--Gromov inequality implies the weak Bishop--Gromov inequality at
any scale and with the same values of the factor~$C$ and of the exponent $K$. Example~\ref{nonBishop} proves that the strong condition is in fact much
stronger than the weak one; furthermore Example~\ref{infinitetopologies} proves that, contrarily to the strong condition, the weak one implies no restriction
on the topology of the spaces under consideration.

In order to deduce a weak Bishop--Gromov inequality at a small scale from
a weak Bishop--Gromov
inequality at a large scale, we~introduce the notion of Busemann space (i.e., metric spaces whose distance is convex), whose definition and first
properties
are given in Section~\ref{Busemannspaces}. Notice that a Busemann space is always contractible, thus simply connected.

\subsubsection{Busemann spaces and strong Bishop--Gromov inequality}\label{strengthenBishop0}

\begin{theor}\label{strengthenBishop}
Given $r_0, D > 0$ and $C > 1$, on any Busemann space $(X,d)$ satisfying
the property of extension of local geodesics $($see Definitions~$\ref{geodextension})$, any proper action by isometries of a group $\Gamma$ such that the diameter of~$\Gamma \backslash X$ is $\le D$ satisfies the following property:
if there exists a $\Gamma$-invariant measure on~$X$ which satisfies a weak Bishop--Gromov inequality at scale $r_0$, with factor $C$
and exponent $K$ then, for every $x \in X$ and every $r, R> 0$ such that $r \le R$, the counting measure $\mu_{x}^{\Gamma}$ verifies:
\begin{itemize}\itemsep=0pt
\item[$(i)$]
if $r \ge 2 r_0 $, then
\[
\dfrac{\mu_{x}^{\Gamma} \big(B_X( x, R)\big)}{\mu_{x}^{\Gamma}
\big(B_X (x, r)\big)} \le C
\bigg(1 +2\,\dfrac{R}{r} \bigg)^{\frac{\ln C}{\ln 2}}{\rm e}^{K \left(R + \frac{r}{2}\right)},
\]
\item[$(ii)$]
if $0 < r <2 r_0$, then
\[
\dfrac{\mu_{x}^{\Gamma} \big(\overline B_X(x, R)\big)}{\mu_{x}^{\Gamma} \big(B_X (x, r)\big)} \le C \bigg(\bigg(1 + \dfrac{D}{r_0}\bigg)\bigg(1 +2\,\dfrac{R}{r}\bigg)\bigg)^{\frac{\ln C}{\ln 2}}
{\rm e}^{K (D+ r_0) \left(1 +2\,\frac{R}{r}\right)},
\]
\item[$(iii)$]
for every $r$, $R$ such that $0 < r < R$ and $r \le r_0$,
\[
{\rm Pack} (x, r, R) \le C \bigg(\bigg(1 + \dfrac{D}{r_0}\bigg) \frac{R}{r} \bigg)^{\frac{\ln C}{\ln 2}} {\rm e}^{K (D+ r_0) \frac{R}{r}}.
\]
\end{itemize}
\end{theor}

 Remark that Example~\ref{exemplebis} (resp.\ Example~\ref{nonBishop}) proves that, without the hypothesis of extension of local geodesics (resp.\ without
the hypothesis ``Busemann"), it is generally impossible to deduce a
weak Bishop--Gromov inequality at scale $< D$ from a weak Bishop--Gromov inequality
at scale $r_0 > D$, even if all the other hypotheses of Theorem~\ref{strengthenBishop} are verified.

\medskip

A question is thus: \textit{can we replace the Busemann hypothesis by a weaker one?}

\medskip

A first alternative hypothesis is to suppose that only the balls of radius $(D+ r_0) \big(1 +2\,\frac{R}{r}\big)$ are Busemann.
A second alternative hypothesis is to replace the Busemann condition of convexity of the distance by a condition of convexity\footnote{A geodesic metric space $(X, d)$ is said to verify the condition of convexity modulo
the defect $\e_0 > 0$ if,
for every pair of geodesics $c_0$ and $c_1\colon [0, 1] \f X $, with common origin $c_0 (0) = c_1 (0)$, endowed with their natural parametrization (see Definition~\ref{naturalparameter}),
one has $ d \big(c_0(t), c_1 (t) \big) \le t\cdot d \big(c_0 (1), c_1 (1) \big) + \e_0$.} modulo some convexity defect $\e_0 > 0$ (see~\cite{BCGS2} where these ideas are developed). It~is worth noticing that, in~a $\delta$-hyperbolic space, the distance verifies a condition of convexity modulo the defect $\delta$.

We shall only sketch here the proof of Theorem~\ref{strengthenBishop} (see~\cite{BCGS2} for complete proofs of this theorem and of the two following
preliminary lemmas).

Let $(X, d) $ be any Busemann metric space and $x$ be any point of~$X$, for every $\lambda \in{} ] 0, 1]$, let us consider the map
$\varphi_{x, \lambda}\colon X \f X$ defined as follows: for every point $y $,
consider the geodesic segment $c_y$ from $x$ to $y$, endowed with its natural parametrization (i.e., $c_y (1) = y$, see Defi\-nition~\ref{naturalparameter}) and define $\varphi_{x, \lambda} (y) := c_y \left( \lambda\right)$. The~following lemma was suggested to us by~\cite[Lemma~4.5]{CS4}, established under ${\rm CAT}(K)$ hypotheses:

\begin{Lemma}\label{contractionbis} If $(X, d) $ is a Busemann space, then $d\big( \varphi_{x, \lambda} (y), \varphi_{x, \lambda} (y')\big) \le
\lambda \, d (y,y')$ for every $x, y, y' \in X$ and every
$\lambda \in{} ] 0, 1]$. If~moreover $(X, d) $ satisfies the property of extension of local geodesics $($see Definitions~$\ref{geodextension})$ then, for every $x \in X$ and every
$ \e, r, R$ such that $ 0 < \e < r \le R$, one has $ {\rm Pack} (x, \e, r) \le {\rm Pack} \big(x, \frac{R}{r} \e,\frac{R}{r} r\big)$.
\end{Lemma}

Another useful (and almost trivial) argument is the

{\sloppy
\begin{Lemma}\label{CovervsPack0} For every proper cocompact action, by isometries preserving the measure, of~a~group~$\Gamma$ on a~metric measure space $(X,d, \mu)$, if $D$ is an upper bound of the diameter of~$ \Gamma \backslash X$ then, for every $r, R$ such that $D < r < R$, one has
${\rm Pack}(x,r,R)\le {\rm Pack}_\Gamma (x,r - D,R) \le \frac{\mu (B_X ( x, R))}{\mu (B_X (x, r- D))}$.
\end{Lemma}

}

\begin{proof}[Sketch of the proof of Theorem~\ref{strengthenBishop}]
If $r \ge 2 r_0$, applying first inequality~\eqref{compardoublings1} and afterwards Lemma~\ref{Bishopclassic}, the weak Bishop--Gromov hypothesis yields:
\begin{gather*}
\dfrac{\mu_{x}^{\Gamma} \big(B_X(x, R)\big)}{\mu_{x}^{\Gamma} \big(B_X (x, r)\big)}
\le \dfrac{\mu \big(B_X\big(x, R + \frac{r}{2}\big)\big)}{\mu \big(B_X\big(x,\frac{r}{2}\big)\big)}
\le C \bigg(1 +2 \, \frac{R}{r}\bigg)^{\frac{\ln C}{\ln 2}} {\rm e}^{K \left(R + \frac{r}{2}\right)},
\end{gather*}
and this proves $(i)$.

Let us now suppose that $r < 2 r_0$. For every $\alpha > 1$, using first~\eqref{compardoublings1} and then Lemma~\ref{contractionbis}, one has
\begin{gather*}
\dfrac{\mu_{x}^{\Gamma} \big(B_X (x, R\big)\big)}{\mu_{x}^{\Gamma}\big(B_X (x, r)\big)}
\le {\rm Pack}_\Gamma \bigg(\!x, \frac{r}{2}, R + \frac{r}{2}\bigg)
\le {\rm Pack} \bigg(\!x, \frac{r}{2}, R + \frac{r}{2}\bigg)
\le {\rm Pack} \bigg(\!x,\alpha\,\frac{r}{2}, \alpha \bigg(R + \frac{r}{2} \bigg)\!\bigg).
\end{gather*}
Choosing now $\alpha := 2\, \frac{D + r_0}{r}$, applying first Lemma~\ref{CovervsPack0}, and then the weak Bishop--Gromov inequality at scale $r_0$
and Lemma~\ref{Bishopclassic}, we~obtain:
\begin{gather*}
\dfrac{\mu_{x}^{\Gamma} \big(B_X(x, R)\big)}{\mu_{x}^{\Gamma} \big(B_X (x,r)\big)}
\le {\rm Pack}_\Gamma \bigg(x,r_0, (D+ r_0)\bigg(1 +2\,\frac{R}{r}\bigg)\bigg)
\le \dfrac{\mu \big(B_X \big(x, (D+ r_0) \big(1 +2\,\frac{R}{r}\big)\big)\big)}{\mu \big(B_X(x,r_0)\big)}
\\ \hphantom{\dfrac{\mu_{x}^{\Gamma} \big(B_X(x, R)\big)}{\mu_{x}^{\Gamma} \big(B_X (x,r)\big)}}
{}\le C \bigg(\bigg(1 + \frac{D}{r_0}\bigg)\bigg(1 +2\,\frac{R}{r}\bigg)\bigg)^{\frac{\ln C}{\ln 2}} {\rm e}^{K (D+ r_0) ( 1 +2 \, \frac{R}{r})}.
\end{gather*}
This proves $(ii)$. The~proof of $(iii)$ is analogous to the proof of $(ii)$.
\end{proof}

\subsubsection{A lower bound on the systole (revisited)}

In Theorem~\ref{minorsystole} and Corollary~\ref{minorsystole1}, we~gave a lower bound of the global systole when the metric measure space verifies a
Bishop--Gromov inequality at a small scale, equal to the lower bound of the diastole given by Theorem~\ref{transyst0} and Corollary~\ref{transyst1}, thus
much smaller than the diameter~$D$. The~problem is that we generally only
establish these weak Bishop--Gromov ineq\-ualities at a~scale larger than the
diameter (see Example~\ref{nonBishop}), as~in Section~\ref{BGexemples}$(c)$ and Theorem~\ref{cocompact2}. As Theo\-rem~\ref{strengthenBishop} allows to deduce a Bishop--Gromov inequality at a small scale from
a~Bishop--Gromov inequality at a large scale (under the extra-hypothesis that the distance verifies some convexity condition), we~obtain the following estimate
as a corollary of inequality~\eqref{minorsystoleeq}, of Theo\-rem~\ref{strengthenBishop}$(ii)$ and of the fact that (by Corollary~\ref{transyst1}) the torsion-free diastole is bounded from below by $\frac{r_0}{N_0}$,
where $N_0 = N_0 (C, Kr_0)= \frac{1}{2} \,\nu \left(C^3 {\rm e}^{15 K r_0} + 1\right)$ (for a complete proof, see~\cite{BCGS2}):

\begin{theor}\label{minorsystole5}
{\sloppy
Given $r_0, D >0$, $C > 1 $ and $K\ge 0$, given a Busemann space $(X,d)$
and a~proper co-compact action by isometries of a group
$\Gamma \in \text{\rm Hyp}_{\rm act}$ on this space such that $\diam (\Gamma \backslash X) \le D$, if there exists some
$\Gamma$-invariant measure $\mu$ on~$X$ which verifies a weak Bishop--Gromov inequality at scale $r_0$, with factor $C$ and exponent $K$, then every torsion-free element $\g$ of~$\Gamma$ verifies
\begin{gather*}
\inf_{x \in X} d(x, \g x) \ge \inf_{x \in X} \sys^{\diamond}_\Gamma (x)
\\ \hphantom{\inf_{x \in X} d(x, \g x) }
{}\ge C^{-13/5} D \bigg(\bigg(1 + \frac{D}{r_0}\bigg) \bigg(1 + 4 N_0\,\frac{D}{r_0}\bigg)\bigg)^{-\frac{\ln C}{\ln 2}} {\rm e}^{ - 3 K (D+ r_0) ( 1 + 4 N_0 \frac{D}{r_0})}.
\end{gather*}}
\end{theor}

\subsection{Finiteness and compactness results}\label{finicompact}

The proofs of the results of this section are only sketched, their statements have been sometimes weakened for the sake of simplicity, for sharper
statements and complete proofs see~\cite{BCGS2}. For example, for the sake of simplicity, we~shall consider here only torsion-free groups $\Gamma$, for a~version of the following results considering non torsion-free groups $\Gamma$, see~\cite{BCGS2}.

\subsubsection{Bounds for the number of groups modulo isomorphisms}\label{boundgenerators0}

A reference for the results of this section, though assuming different hypotheses, is the
following result of M.~Anderson

\begin{theor}[{\cite[Theorem~2.3]{An}}]
Given $K, D, V > 0$ and
$n \in \N \setminus \{0,1\}$, in~the class of closed $n$-dimensional
Riemannian manifolds $M$ such that $\Ric_M \ge -(n-1) K^2\cdot g$, $\diam_M \le D$ and $\Vol_M \ge V$, there are finitely many isomorphism classes
of~$\pi_1 (M)$.
\end{theor}

We first obtained an analogous to this theorem when we replace the assumption ``Ricci curvature bounded from below'' by a strong Bishop--Gromov
inequality (see~\cite{BCGS2}); but we now follow a different approach.
Recalling that $\Sigma_{r} (x) := \{ \sigma \in \Gamma \colon d(x, \sigma x) \le r\}$, we~have the

\begin{prop}\label{majorgroupes1}
Given $\e'_0, \,D >0$, $C > 1 $ and $K\ge 0$, there exist
$N'_1 = N'_1 \big(C, \frac{D}{\e'_0}, K D \big)$ and a~list of groups $G_1,\dots, G_{N'_1}$
such that any arcwise connected and simply connected measure length space
$(X,d,\mu)$ which satisfies a weak Bishop--Gromov inequality at scale
\smash{$\frac{\e'_0}{2}$}, with factor~$C$ and exponent $K$, has the following property: every torsion-free group $\Gamma$ acting properly on $(X,d,\mu)$,
by isometries
preserving the measure, with $\diam (\Gamma \backslash X) \le D$ and such
that there exists a point $x \in X$ where $\sys_\Gamma (x) > \e'_0$, is isomorphic
to one of the $G_k$'s.
\end{prop}

The proof of this proposition and of the following ones will be sketched in the following Sections~\ref{boundgenerators0}$(a)$ and $(b)$.
$N_0 (C, Kr_0)= \frac{1}{2} \,\nu \big(C^3 {\rm e}^{15 K r_0} + 1\big)$ being the universal constant which occurs in Proposition~\ref{Marg1}$(ii)$, we~deduce the

\begin{Corollary}\label{majorgroupes2}
Given $r_0, D >0$, $C > 1 $ and $K\ge 0$; there then exist $N_2 = N_2 \big(C, K r_0, \frac{D}{r_0}\big)$ and a list of groups $G_1,\dots, G_{N_2} $
such that any arcwise connected and simply connected measure length space
$(X,d,\mu)$ which satisfies a weak Bishop--Gromov inequality at scale
$\frac{r_0}{3\, N_0 (C, Kr_0)}$, with factor $C$ and exponent $K$, has the following property: every torsion-free group $\Gamma \in {\rm Hyp}_{\rm act}$ acting properly on $(X,d,\mu)$, by isometries preserving the measure, such that $\diam (\Gamma \backslash X) \le D$, is isomorphic to one of the $G_k$'s.
\end{Corollary}

In this last corollary, the weak point is that we assume a weak Bishop--Gromov inequality at a rather small scale; however we can deduce this
from a weak Bishop--Gromov inequality at a~large scale $r_0$ by using the results of Section~\ref{strengthenBishop0}. An example of corollary is the

\begin{Corollary}\label{majorgroupes3}
Given $r_0$, $D >0$, $C > 1 $ and $K\ge 0$, there exist $N_3 = N_3 \big(C, K r_0, \frac{D}{r_0}\big)$ and a list of groups $G_1,\dots, G_{N_3}$
such that any Busemann metric measure space $(X,d,\mu)$ which verifies the property of extension of local
geodesics and which satisfies a~weak Bishop--Gromov inequality at scale $r_0$, with factor $C$ and exponent $K$, has the following property: every torsion-free group
$\Gamma \in \text{\rm Hyp}_{\rm act}$ acting properly on $(X,d,\mu)$, by isometries preserving the measure, with $\diam (\Gamma \backslash X) \le D$, is isomorphic to one of the $G_k$'s.
\end{Corollary}

\begin{theor}\label{majorgroupes5}
Given $\delta, K, D\! >0$, there exist $N_4\! =\! N_4 \big(\delta, K, D \big)$ and a list of groups $G_1,\dots, G_{N_4}$ with the following property:
every torsion-free group, which admits a proper action by isomet\-ries on some $\delta$-hyperbolic Busemann space $(X,d)$, with
the property of extension of local geodesics and entropy $\le K$, such that $\diam (\Gamma \backslash X) \le D$, is isomorphic to one of the $G_k$'s.
\end{theor}

In the sequel, $\Gamma$ is any torsion-free group acting properly, co-compactly, by isometries, on~a~path-connected metric space $(X,d)$ such that
$\diam (\Gamma \backslash X) \le D$ and we fix a point $x \in X$ where the diastole is attained.

\medskip\noindent
{\bf (\emph{a}) A bound of the number of generators.}
Let us recall that, for every $x \in X$, $\Sigma_{2 D} (x)$ is the set of
the elements $\g$ of~$\Gamma$ such that $d(x, \g x) \le 2 D$. A result of M.~Gromov~\cite[Proposition~3.22]{Gr1} proves that $\Sigma_{2 D} (x)$ is a symmetric generating set of~$\Gamma $, the properness of
the action implying the finiteness of~$\Sigma_{2 D} (x) $. One gets the following bound for the number of these generators:

\begin{prop}
For every data $\e'_0$, $D >0$, $C > 1 $ and $K\ge 0$, given a length space
$(X,d)$ and a proper co-compact action by isometries of a torsion-free group
$\Gamma $ on this space such that the diameter of~$\Gamma \backslash X$ is bounded above by $D$, if the diastole of this action
on~$(X,d)$ is $ > \e'_0$ and if there exists some
$\Gamma$-invariant measure $\mu$ on~$X$ which verifies a weak Bishop--Gromov inequality at scale \smash{$\frac{\e'_0}{2}$}, with factor $C$ and exponent $K$, then
there exists $x \in X$ such that
\begin{gather*}
 \# \big( \Sigma_{2 D} (x) \big) \le C \bigg(1 + \dfrac{4 D }{\e'_0}\bigg)^{\frac{\ln C}{\ln 2}} {\rm e}^{K (2 D + \e'_0/2)}.
 \end{gather*}
\end{prop}

\begin{proof}
Let $x$ be a point of~$X$, where $\sys^{\diamond}_\Gamma (x) > \e'_0$ then, by inequality~\eqref{compardoublings1} and Lemma~\ref{Bishopclassic},
\begin{gather*}
\# \big( \Sigma_{2 D} (x) \big) = \dfrac{\mu_{x}^{\Gamma} \big( \overline B_X (x, 2 D) \big)}{\mu_{x}^{\Gamma} \big( \overline B_X (x, \e'_0) \big)} \le \dfrac{\mu \big( B_X (x, 2 D + \e'_0/2) \big)}{\mu \big( B_X (x, \e'_0/2) \big)}
\\ \hphantom{\# \big( \Sigma_{2 D} (x) \big)}
{}\le C \bigg(1 + \dfrac{4 D }{\e'_0}\bigg)^{\frac{\ln C}{\ln 2}} {\rm e}^{K (2 D + \e'_0/2)}.
\tag*{\qed}
\end{gather*}
\renewcommand{\qed}{}
\end{proof}

Using the lower bounds of the diastole given in Section~\ref{diastolebound}$(c)$, we~deduce upper bounds of~$\# \big(\Sigma_{2 D} (x) \big) $ in
the case where a strong Bishop--Gromov inequality is assumed.

In the case where only a weak Bishop--Gromov inequality, at~scale $r_0 $ larger than the diastole, is assumed, we~deduce, from this first Bishop--Gromov inequality, a~weak Bishop--Gromov inequality, at~scale $\e'_0 $ smaller than the diastole, in~the case where $(X,d)$ is a Busemann space satisfying the property of extension of local geodesics (see Definitions~\ref{geodextension}). Indeed, we~may then apply Theorem~\ref{strengthenBishop} and carry on with
the proof as before.

\medskip\noindent
{\bf (\emph{b}) Bounding the number of groups modulo isomorphisms.}
In the sequel, as~in part~$(a)$, we~denote by $\Sigma := \Sigma_{2 D} (x) $ the canonical generating set of~$\Gamma$, by $\mathbb F (\Sigma)$
the free group generated by $\Sigma$ and by $\varphi_{\Sigma} \colon \mathbb F
(\Sigma) \f \Gamma$ the canonical epimorphism which maps $\Sigma$ onto
$\Sigma$. Obviously $\varphi_{\Sigma}$ induces an isomorphism: $ \mathbb
F (\Sigma) / {\Ker \varphi_{\Sigma}} \to \Gamma $.

Once the number of elements of~$\Sigma := \Sigma_{2 D} (x) $ is bounded
by some integer $N$ (see part~$(a)$), the problem is to associate,
to each group $\Gamma$, one of its presentation $\langle \Sigma, R \rangle$ by generators and relations\footnote{Recall that $\langle \Sigma, R \rangle$ is a presentation of~$\Gamma$ by generators and relations if the kernel of~$\varphi_{\Sigma}$ is the smallest normal subgroup of~$\mathbb F (\Sigma)$ which contains $R$.}
such that the word-length\footnote{The word-length of an element $g \in \mathbb F (\Sigma)$ is the minimal number of factors when one writes
$g$ as a product of elements of~$\Sigma$.} of all the relations belonging
to $R$ is bounded above independently of~$\Gamma$ and $\Sigma$. If~such a bound, say $p$, exists, the number of possible presentations $\langle \Sigma, R \rangle$ with $\# \Sigma = k$ is smaller than the number of subsets $R$
of the ball of radius $p$ in $\mathbb F (\Sigma)$, thus smaller than $2^{(2 k)^p}$. It~follows that the number of possible presentations $\langle \Sigma, R \rangle$ (and consequently the number of possible groups with less than $N$ generators modulo isometric isomorphisms) is bounded above by $\sum_{k = 2}^N 2^{(2k)^p}$.

There are two cases where we know how to bound the word-length of all the
elements of~a~gene\-ra\-ting set of relations, and thus the number of groups with
less than $N$ generators:

\begin{itemize}\itemsep=0pt
\item When $(X,d)$ is simply connected, each group $\Gamma$ admits a presentation $\langle \Sigma_{2 D} (x), R \rangle$
such that the word-length of all the elements of~$R$ is $\le 3$. This is an immediate consequence of the following lemma, which is due to J.-P.~Serre~\cite[p.~30]{Se} and M.~Gromov~\cite{Gr1}:

\begin{theor}
On every arcwise connected and simply connected topological space $X$, every action by
homeomorphisms of a group $\Gamma$ has the following properties: for every arcwise connected open set $U \subset X$ such that
$ \cup_{\g \in \Gamma} \, \g (U) = X$, the set $\Sigma := \{ \g \in \Gamma \colon \g(U) \cap U \ne \varnothing\} $ is a symmetric generating set of~$\Gamma$ and
there exists a generating set
$R$ of the relations\footnote{The generating set $R$ is defined as the set of~$ s\cdot \sigma \cdot t \in \mathbb F (\Sigma)$ such that $U \cap \varphi_\Sigma (s).U \cap \varphi_\Sigma (s \sigma).U \ne \varnothing$ and $ \varphi_\Sigma (t) = \varphi_\Sigma \big( (s \sigma)^{-1}\big)$.} between elements of~$\Sigma$ such that each element of~$R$ can
be written $ s\cdot \sigma \cdot t$, where $s,\sigma, t \in \Sigma$.
\end{theor}
To conclude in this case, it is sufficient to make $U := B_X (x, r)$ in
this theorem (for any $r > D$) and to recall that $(X,d)$ is a length space; this ends the sketch of the proof of~Proposition~\ref{majorgroupes1} and Corollaries~\ref{majorgroupes2} and~\ref{majorgroupes3}.

\item When $(X,d)$ is $\delta$-hyperbolic and not any more assumed to be simply connected, it follows from the works of M.~Gromov on $\delta$-hyperbolic groups
(see~\cite[Theorems~3.22 and~5.12, p.~88]{GH}) that every group acting properly and cocompactly by isometries on a
Gromov $\delta$-hyperbolic space is a finitely generated $\delta'$-hyperbolic group (with $\delta' := \delta' (\delta, D)$) when endowed with the generating set $\Sigma_{k D} (x))$ (with $k \ge 3$). For example, choosing the gene\-ra\-ting set $ \Sigma_{6 D} (x))$, we~obtain that $\delta' (\delta, D)
= 8 \big( 5\,\frac{\delta}{D} + 4\big)$ (see~\cite{BCGS2}).

One can find the proof of the following theorem of M.~Gromov in several books, here we have followed the quantified result given in~\cite[Chapitre~5, proof of
Theorem~2.3]{CDP} (see also~\cite[Chapitre~4, Proposition~17]{GH}), whose result may be rewritten:

\begin{theor}
Let $(\Gamma, \Sigma)$ ($\#\Sigma <+\infty$) be a $\delta'$-hyperbolic marked group, and $d_\Sigma$ the associated algebraic distance,
then $S = S_\Sigma := \{\g \in \Gamma^*\colon d_\Sigma (e, \g) < 4 \delta' + 2 \}$ is a symmetric generating set of~$\Gamma$ with the following properties: let $\mathbb B_S (3)$ be the set of elements of word-length $2$ or $ 3$ in the free group $\mathbb F(S)$
and choose $R := \Ker \varphi_S \cap \mathbb B_S (3)$, then $\langle S, R\rangle$ is a~pre\-sentation of~$\Gamma$.
\end{theor}

As the numbers of elements of~$S$ and $\mathbb B_S (3)$ are bounded in terms of~$\delta'$ and of the upper bound $N$ of~$\# \Sigma$, this concludes the sketch of the proof of Theorem~\ref{majorgroupes5}.
\end{itemize}

\subsubsection{Bounding the number of spaces modulo homotopy equivalences}

Given $r_0, D >0$, $C > 1$ and $K\ge 0$, define ${\cal M} (r_0, C, K,
D)$ as the set of arcwise connected length spaces $(X,d)$, with diameter $\le D$, whose
fundamental group $\Gamma_X$ is torsion free and belongs to~${\rm Hyp}_{\rm act}$, and which admit a metric universal cover $\big(\widetilde X, \tilde d\big)$
satisfying the following properties: $\big(\widetilde X, \tilde d\big)$ is Busemann with the property of extension of local geodesics and admits a~$\Gamma_X$-invariant
measure which verifies a Bishop--Gromov inequality at scale $r_0$, with factor $C$ and exponent $K$. The~homotopy class of a space $(X,d)$ being the subset
of all the spaces $(Y,d) \in {\cal M} (r_0, C, K, D)$ which are homotopically equivalent to $(X,d)$ and
$N_3 \big(C, K r_0, \frac{D}{r_0}\big)$ being the universal constant introduced in Corollary~\ref{majorgroupes3}, we~obtain the

\begin{theor}\label{boundhomotopies1}
For every $r_0$, $D >0$, $C > 1$ and $K\ge 0$, there are less than $N_3 \big(C, K r_0, \frac{D}{r_0}\big)$ homotopy classes of~spaces in ${\cal M} (r_0, C, K, D)$.
\end{theor}

Let us now define ${\cal H} (\delta, K, D)$ as the set of arcwise connected length spaces $(X,d)$, with diameter $\le D$, whose
fundamental group $\Gamma_X$ is torsion free, and which admit a metric universal cover $\big(\widetilde X, \tilde d\big)$
satisfying the following properties: $\big(\widetilde X, \tilde d\big)$ is Busemann, $\delta$-hyperbolic, with the property of extension of local geodesics and
entropy $\le K$. Recalling that $N_4 \left(\delta, K, D \right)$ is the universal constant introduced in Theorem~\ref{majorgroupes5}, we~obtain the

\begin{theor}\label{boundhomotopies2}
For every $ \delta, K, D >0$, there are less than $N_4 \left(\delta, K,
D \right)$ homotopy classes of
spaces in ${\cal H} (\delta, K, D)$.
\end{theor}

The proof of Theorem~\ref{boundhomotopies1} (resp.\ of Theorem~\ref{boundhomotopies2}) is a consequence of the fact that, modulo isomorphisms,
the number of fundamental groups of spaces belonging to ${\cal M} (r_0, C, K, D)$ (resp.\ to ${\cal H} (\delta, K, D)$) is less than $N_3 \big(C, K r_0, \frac{D}{r_0}\big)$ (resp.\ than
$N_4 (\delta, K, D)$) by Corollary~\ref{majorgroupes3} (resp.\ by Theorem~\ref{majorgroupes5}), and of the contractibility of the universal cover
of every $(X,d)$ in ${\cal M} (r_0, C, K, D)$ (resp.\ in ${\cal H} (\delta, K, D)$) along the geodesics issued from one point.

\medskip

The problem of trying to prove that, on some families of metric spaces which are moreover compact for the Gromov--Hausdorff distance, a~homotopical
finiteness result implies a finiteness theorem for topologies has a quite
long history, it would thus be interesting to prove that
the family ${\cal M} (r_0, C, K, D)$ (resp.\ ${\cal H} (\delta, K, D)$) is compact with respect to the Gromov--Hausdorff distance. A section of~\cite{BCGS2} is devoted to this problem, however we already proved the

\begin{theor}\label{GHcompactness}
The family ${\cal H}_0 (\delta, K, D)$ of the elements of~${\cal H} (\delta, K, D)$ which are locally ${\rm CAT}(0)$ is compact with respect to the Gromov--Hausdorff distance.
Moreover there exists a universal constant $\e_0 = \e_0 (\delta, K, D)
> 0$ such that every ball of radius $\e_0$, in~any $(X,d) \in {\cal H} (\delta, K, D)$
which is not a circle, is contractible.
\end{theor}

In addition to the compactness property, this last theorem proves the existence of a lower bound of the contractibility radius; this is an important tool in
proving that homotopical finiteness implies topological finiteness. This lower bound of the contractibility radius is a consequence of the lower bound of the systole
given by Theorems~\ref{minorsystole5} and~\ref{cocompact2}.

\begin{question}
Does $\,{\cal H}_{0} (\delta,H,D)\,$ and ${\cal M} (r_0, C, K, D)$ only contain a finite number of topo\-logies?
\end{question}

In~\cite{BCGS2}, we~give a partial answer to this question (this is only a work in progress):
the subset of the elements of~${\cal H}_0 (\delta, K, D)$
(resp.\ of~${\cal M} (r_0, C, K, D)$) which are topological manifolds only contains a finite number of topologies

A (very rough) sketch of the proof is as follows:
\begin{itemize}\itemsep=0pt
\item
In dimension $\ne 3$, it is a consequence of the following result of S.~Ferry~\cite[Theorem~1]{Fe}: in the set of compact metric spaces
(endowed with the Gromov--Hausdorff distance), any precompact subset, whose elements all have contractibility radius uniformly bounded from below, contains a finite number of topologies. Ferry's result concludes a series
of results of L.~Siebenmann~\cite{Sie}, M.~Gromov~\cite{Gro-Large}, F.~Farrell--L.~Jones~\cite{FJ}, K.~Grove--P.~Petersen--J.~Wu~\cite{GPW} (corrected in~\cite{GrPW}), R.~Greene--P.~Petersen~\cite{GPe}, \dots\ this list being
non exhaustive. We~thus have to prove first that
${\cal H}_0 (\delta, K, D)$ (resp.\ ${\cal M} (r_0, C, K, D)$) is precompact, secondly that the contractibility radii of all the elements of~${\cal H}_0 (\delta, K, D)$ (resp.\ of~${\cal M} (r_0, C, K, D)$) are bounded from below. The~precompactness is a corollary of
M.~Gromov's precompactness theorem (see~\cite[Proposition~5.2]{Gr1}) and of Theorem~\ref{strengthenBishop}$(iii)$, which bounds the number of~$\e$-balls
in a packing of every element of~${\cal H}_0 (\delta, K, D)$ (resp.\ of~${\cal M} (r_0, C, K, D)$). The~lower bound of the contractibility radius is given, in~the Busemann case, by Theorem~\ref{GHcompactness}.

\item In dimension 3, the same finiteness result is a corollary of G. Perelman's proof of Thurston's geometrisation conjecture and Poincar\'e's conjecture, with the addition of a result of M.~Kreck and W.~L\"uck~\cite[Theorem~0.7]{KL}.
\end{itemize}

\subsection{Appendices}

\subsubsection{Geodesics}\label{geodesics}
The following definitions are classical (see for example~\cite[D\'efinitions~I.1.3, p.~4]{BH}):
\begin{defis}
In any metric space $ (X,d)$
\begin{itemize}\itemsep=0pt
\item a (normal) ``geodesic'' is a map $c$ from some interval $I \subset \R$ to $X$ verifying the property: for every $ t, \, t' \in I$
$ d(c(t), c(t')) = |t - t'|$,

\item when $I$ is a closed interval (resp.\ $]{-}\infty, + \infty [$) the geodesic is called a geodesic segment
(resp.\ a geodesic line), the image of a geodesic segment $ c $ with origin $x$ and endpoint $y$ is often denoted by $[x, y]$
(though this does not suppose that this geodesic segment is unique),

\item a (normal) ``local geodesic'' is a map $c$ from some interval
$I \subset \R$ to $X$ verifying the property: for every $ t \in I$,
there exists $\e > 0$ such that $ d(c(t'), c(t'')) = |t' - t''|$ for every $ t', \, t'' \in \, ]t- \e, t+ \e[\,\cap I$,

\item a metric space $(X,d)$ is geodesic if any two points can be joined by at least one geodesic.
\end{itemize}
\end{defis}

\begin{defi}\label{naturalparameter}
Given a geodesic segment $[x_0, x_1]$, the ``natural parametrization'' of this segment is the map $t \mapsto x_t$ from $[0, 1]$ to $[x_0, x_1]$,
defined by $d(x_0, x_t) = t \, d(x_0, x_1)$.
\end{defi}

\begin{defis}[{cf.~\cite[Definition~II.5.7, p.~208]{BH}}]\label{geodextension}
A geodesic metric space is said to verify the ``property of extension of local geodesics'' if, for every local geodesic $c \colon [a,b] \to X$ ($a
< b$),
there~exists $\e > 0$ and a local (eventually not unique) geodesic $ c' \colon
[a,b+ \e ] \to X $ which extends~$c$ $\big($i.e.,~$ c'_{|_{[a,b]}} = c\big)$.

This space is said to be ``geodesically complete'' if every local geodesic $c\colon [a,b] \to X$ ($a < b$), can be extended as a local geodesic
$ \bar c\colon ]{-}\infty, +\infty[ {} \to X$.
\end{defis}
Notice that every complete or open Riemannian manifold verifies the property of extension of local geodesics.

\begin{Lemma}[{cf.~\cite[Lemma~II.5.8(1), p.~208]{BH}}]
A \emph{complete} geodesic metric space verifies the property of extension of local geodesics if and only if it is geodesically complete.
\end{Lemma}

\subsubsection{About Gromov-hyperbolic spaces}

Given three nonnegative numbers $\alpha$, $\beta$, $\gamma$, we~define
the tripod $T := T(\alpha, \beta, \gamma)$ as the metric simplicial tree
with $3$ vertices
$x'$, $y'$, $z'$ of valence $1$ (the ``endpoints''), one vertex $c$
of valence $3$ $($the ``branching point''$)$, and $3$ edges $[cx']$,
$[cy']$, $[cz']$
of respective lengths $\alpha$, $\beta$, $\gamma$ (the ``branches''). We~denote by $d_T (u,v) $ the distance on this tree between two points $u,v \in
T$, i.e., the minimal length of a path contained in $T$ and joining $u$ to
$v$.

\medskip

For the sake of simplicity, we~only consider geodesic metric spaces (see definition in Section~\ref{geodesics}). In~such a space a geodesic triangle
$\Delta = [x, y, z]$ is the union of three geodesics $[x, y ]$, $[y, z ]$ and $[z, x ]$. Given three points $x,y,z$ in a geodesic metric space, there
exists at least one geodesic triangle $\Delta = [x, y, z]$ whose sides have respective lengths $d(x,y),\, d(y,z)$ and $d(x,z)$.

\begin{Lemma}\label{prodist}
To any geodesic triangle $\Delta$ corresponds a metric tripod $(T_\Delta, d_T)$ and a surjective map $f_\Delta \colon \Delta \f T_\Delta$ $($called the
``approximation of~$\Delta$ by a tripod"$)$ such that, in~restriction
to each side of~$\Delta$, $f_\Delta$ is an isometry,
\end{Lemma}

Indeed, $T_\Delta$ is constructed as the tripod $ T(\alpha, \beta, \gamma)$, where (by the triangle inequality) $(\alpha, \beta, \gamma)$ is the
unique element of~$[0, +\infty[^3$ such that $d(x, y) = \alpha + \beta$, $d(x, z) = \alpha + \gamma$ and $d(y, z) = \beta + \gamma$. This choice of~$(\alpha, \beta, \gamma)$ implies the existence of the map $f_\Delta \colon \Delta \f T_\Delta$ as asserted in Lemma~\ref{prodist}.

\begin{defis}\label{hypdefinition0}
A geodesic triangle $\Delta$ of~$(X,d)$ is said to be $\delta$-thin if, for every $u \in T$ and every $x, y \in f_\Delta^{-1} (\{u\})$, one has $d(x,y) \le
\delta$.

In this article, a~metric space is said to be $\delta$-hyperbolic if it is geodesic, proper, and if all its geodesic triangles are $\delta$-thin.
\end{defis}

For more informations about Gromov-hyperbolic metric spaces see the original publication~\cite{Gr4} and, for more explanations,~\cite{CDP,GH} and~\cite{BH}.

\subsubsection{Busemann spaces}\label{Busemannspaces}

\begin{defi}
A Busemann space is a geodesic and proper metric space (not reduced to one point) whose distance between geodesics is convex, i.e. for every pair
 of geodesic segments $c_1$, $c_2 $ such that $c_1 (0) = c_2 (0)$ (with their natural parametrization, see Definition~\ref{naturalparameter}),
the
function $t \mapsto d(c_1(t), c_2 (t))$ is convex.
\end{defi}

Examples of Busemann spaces are Hadamard spaces and, more generally, $\Cat (0)$-spaces.

\medskip

The following results are classical and based on the facts that every Cauchy sequence is bounded and that every Busemann space is proper:

\begin{Remark}
Every Busemann space is a complete metric space.
\end{Remark}

\begin{Remark}
On every Busemann space, any two points are joined by a single geodesic.
\end{Remark}

\begin{Lemma}
On any Busemann space, every local geodesic is a $($minimizing$)$ geodesic.
\end{Lemma}

\begin{Lemma}
Every Busemann space $(X,d)$ satisfying the property of extension of local geo\-de\-sics is non compact and every geo\-de\-sic extends as
a $($minimizing$)$ geodesic line $ ]{-} \infty, +\infty [ {}\to X$.
\end{Lemma}

\subsubsection{Examples}

\begin{exa}\label{exemplebis} A family $(X_{\e}, d_{\e})_{\e > 0}$ of~$\delta_0$-hyperbolic spaces with entropy $ \le H_0$, which admit a proper isometric action of a group $\Gamma$ such that $\diam (\Gamma \backslash X_{\e}) \le D$ ($\delta_0, H_0, D$ being independent of~$\e$) such that
$(X_{\e}, d_{\e})$ is Busemann, though there exists $x \in X$ such that $\big(X_{\e}, d_{\e}, \mu^\Gamma_x\big)$ do not verify any weak Bishop--Gromov inequality at any scale $r_0 < 2 D$ when $\e$ is small enough.
\end{exa}

\begin{proof}Choose any $r_0, D$ such that $0 < r_0 < 2 D$; for every $ \e \le 2 D - r_0$, $(X_{\e}, d_{\e})$ is obtained by gluing to $\R$ an interval $I_k = \big[0, \frac{r_0}{2} \big]$ at each point $\e k \in
\e \mathbb Z$,
identifying the origin of~$I_k$ with the point $\e k$. Obviously, the group $\Gamma := \e \cdot \Z$ acts on $(X_{\e}, d_{\e})$ by translations sending $I_k$
onto $I_{p+k}$ for every $ \e\cdot p \in \e \cdot \Z$.
$(X_{\e}, d_{\e})$ is then $0$-hyperbolic (thus Busemann), has trivial Entropy, verifies $\diam (\Gamma \backslash X) = \frac{r_0 + \e}{2} \le D $. Let $x$ be the
endpoint of the interval $I_0$, let $r := r_0 + \e$, it is obvious that
$\mu_{x}^{\Gamma}\big(B_X (x, r) \big) = 1$ and that
$\mu_{x}^{\Gamma}\big(B_X (x, 2 r) \big) \ge 2 \big[ \frac{r_0}{\e}\big] + 3 $. Hence, for every choice of~$C, K$, the inequality
$\frac{\mu_{x}^{\Gamma}(B_X (x, 2 r))}{\mu_{x}^{\Gamma}(B_X(x, r))} \le C\, {\rm e}^{Kr}$ is not verified when $\e$ is small enough.
\end{proof}

In the following examples, saying that almost all the elements of a sequence $\big(\overline M_i, \bar g_i\big)_{i \in \N^*}$ of~Riemannian manifolds
do not verify any weak
Bishop--Gromov inequality at a given scale $r_0$ means that, for every $C, D > 0$, there exists $i_0 $ such that, for every $i \ge i_0$, the Riemannian
measure of~$\big(\overline M_i, \bar g_i\big)$ do not verify the weak Bishop--Gromov inequality at scale $ r_0$, with factor $C$ and exponent $K$.
Similarly, saying that almost all the Riemannian manifolds of this sequence do not verify any weak Bishop--Gromov inequality means that,
for every $r, C, D > 0$, there exists $i_0 $ such that, for every $i \ge i_0$, the Riemannian measure of~$\big(\overline M_i, \bar g_i\big)$
do not verify the weak Bishop--Gromov inequality at scale $ r$, with factor $C$ and exponent $K$.

\begin{exa}\label{nonBishopbis} There exist positive constants $\delta$ and $K_1$ and sequences of~$\delta$-hyperbolic
Riemannian manifolds $\big(\overline M_i, \bar g_i\big)_{i \in \N^*}$, whose entropy is $ \le K_1$ and such that each $\big(\overline M_i, \bar g_i\big)$
admits a compact quotient by a discrete group of isometries, though almost all the $\big(\overline M_i, \bar g_i\big)$'s do not verify any weak
Bishop--Gromov inequality.
\end{exa}

\begin{exa}\label{nonBishop} There exist positive constants $\delta$, $K_1$ and $D$ and a sequence of~$\delta$-hyperbolic
Riemannian manifolds $\big(\overline M_i, \bar g_i\big)_{i \in \N^*}$, such that each $\big(\overline M_i, \bar g_i\big)$ has entropy $\le K_1$,
admits a quotient of diameter $\le D$ by a discrete group of isometries, and satisfies the following properties: on each $\big(\overline M_i, \bar g_i\big)$, the
Riemannian measure verifies a weak Bishop--Gromov inequality at scale
$\frac{5}{2} (7 D + 4 \delta)$, with exponent $\frac{36}{5} K_1$ and factor $C_0 (K_1, D)$, though almost all the $\big(\overline M_i, \bar g_i\big)$'s do not verify
any weak Bishop--Gromov inequality at scale $ \frac{3}{4} D$.
\end{exa}

The factor $C_0 (K_1, D)$ is computed at the end of the following proof.

\begin{proof}[Construction and proof of Examples~\ref{nonBishopbis} and~\ref{nonBishop}]
We consider a sequence of closed Riemannian manifolds of dimension $n$, $(N_i, h_i)_{i \in \N^*}$, with Ricci curvature $\ge - \frac{K_0^2}{n-1}$
and sectional curvature $\le - \frac{K_0^2}{100 (n-1)^2}$, such that $\diam (N_i, h_i) = D_i $, and, on each $N_i$, we~choose a pair of points~$x_i$,~$y_i$ at distance $D_i$. On the other hand, for each $i \in \N^*$, let $(Y_i, k_i)$ be any closed $n$-dimensional Riemannian manifold, with small diameter $\le \frac{1}{i}$ and big volume $\Vol (Y_i, k_i) \ge C_0 \big( {\rm e}^{3 i^2 K_0}\big)$ (such metrics exist on every closed
manifold). We~construct a new sequence of closed $n$-dimensional Riemannian manifolds
$(M_i, g_i)_{i \in \N^*}$, each $(M_i, g_i)$ being obtained as the connected sum of~$(N_i, h_i)$ with $(Y_i, k_i)$; the construction is made by gluing $Y_i$ to
$N_i$ inside
the ball $B_{N_i} \big(y_i, \frac{1}{i} \big)$. More precisely, defining
$r_i := \frac{1}{2} {\rm e}^{-i}\cdot \min \big( \inj (Y_i, k_i), \inj (N_i, h_i)\big)$, we~choose
a point $y'_i \in Y_i$, we~excise one ball $B_{Y_i} (y'_i, r_i)$ (resp.\ $B_{N_i} (y_i, r_i)$) from $Y_i$ (resp.\ from $N_i$) and glue
$N_i \setminus B_{N_i} (y_i, r_i)$ and $Y_i \setminus B_{Y_i} (y'_i, r_i)$ by identification of their boundaries.\footnote{To be isometric, this identification may
require to slighly modify the two metrics $h_i$ and $k_i$ on $ B_{N_i} (y_i, 2 r_i)$ and $ B_{Y_i} (y'_i, 2 r_i)$ respectively, in~order that they become flat
on these two balls.} As~$N_i \setminus B_{N_i} \big(y_i, \frac{1}{i}\big)$ is included in both spaces
$N_i $ and $M_i$,
one can construct a continuous map $f_i\colon M_i \f N_i$, such that $f_i^{-1} (y_i) = Y_i \setminus B_{Y_i} (y'_i, r_i) $), which coincides with
the identity map on $N_i \setminus B_{N_i} \big(y_i, \frac{1}{i}\big)$, the new metric $g_i$ coinciding with
$h_i$ on $N_i \setminus B_{N_i} \big(y_i, \frac{1}{i}\big)$. We~call $(f_i)_* $ the induced morphism from $\pi_1 (M_i, x_i)$ to $\pi_1 (N_i, x_i)$.

Let us now define $\pi_i \colon (\widetilde N_i, \tilde h_i) \f (N_i, h_i)$ as
the Riemannian universal covering of~$(N_i, h_i)$ and choose $\tilde x_i \in \pi_i^{-1} (x_i)$
and $\tilde y_i \in \pi_i^{-1} (y_i)$ such that $d_{\tilde h_i} (\tilde x_i, \tilde y_i) = D_i $. Let now $ p_i \colon \big(\overline M_i, \bar g_i\big) \f
(M_i, g_i)$ be the Riemannian
covering of~$(M_i, g_i)$ such that $(p_i)_* \big(\pi_1 \big(\overline M_i, \bar x_i\big) \big) = \Ker (f_i)_* $ for some point $\bar x_i \in p_i^{-1} (x_i)$ (see~\cite[Corollary~6.9]{GrH}). As~$\Ker (f_i)_*$ is normal in $\pi_1 (M_i, x_i)$, $p_i$ is a Galois covering, and the group ${\rm Aut} (p_i)$ of automorphisms of this covering is the
group $\pi_1 (M_i, x_i) / \Ker (f_i)_* $, hence $(f_i)_* $ induces an isomorphism $f_i^*$ from $\pi_1 (M_i, x_i) / \Ker (f_i)_*$ to the fundamental
group $\Gamma_i$ of~$N_i$, viewed as the group of automorphisms of the universal covering $\pi_i$.
Classicaly (see~\cite[Theorem~(6.1)]{GrH}), there then exists a map $\bar
f_i \colon \overline M_i \f \widetilde N_i$ such that $\pi_i \circ \bar f_i =
f_i \circ p_i$
and $ \bar f_i (\g\cdot \bar x) = f_i^* (\g) \cdot \bar f_i (\bar x)$ for every $\bar x \in \overline M_i $ and every $\g \in {\rm Aut} (p_i)$.

As these two actions of~$\Gamma_i$ commute with $\bar f_i $ and as $f_i$ is the identity map from $ N_i \setminus B_{N_i} \big(y_i, \frac{1}{i}\big)$ onto itself,
direct computations give first that $p_i^{-1} \big( f_i^{-1} (y_i)\big)$ is the union of the left-translates $ \g \bar f_i^{-1} ( \tilde y_i) $ of~$\bar f_i^{-1} ( \tilde y_i) $, all isometric to $ f_i^{-1} (y_i) = Y_i
\setminus B_{Y_i} (y'_i, r_i) $ via the map $p_i$, and secondly that~$\bar f_i$ is isometric from $p_i^{-1} \big( N_i \setminus B_{N_i} \big(y_i, \frac{1}{i}\big)\big) \subset \overline M_i$ to $\widetilde N_i \setminus \bigcup_{\g \in \Gamma_i} B_{\widetilde N_i} \big(\g \tilde y_i, \frac{1}{i}\big)$.

This shows that, from a geometric point of view, $\big(\overline M_i, \bar g_i\big)$ is the connected sum of~$\big(\widetilde N_i, \tilde h_i\big)$ with an infinite
family $\big(Y_i^\gamma, k_i^\gamma\big)_{\g \in \Gamma_i}$ of copies of~$(Y_i, k_i)$, each $Y_i^\gamma = \bar f_i^{-1} (\g \tilde y_i) $ being glued
to~$\widetilde N_i$ inside
the ball $B_{\widetilde N_i} \big(\g \tilde y_i, \frac{1}{i} \big)$.
Recall that $(\widetilde N_i, \tilde h_i)$ is $\delta_0$-hyperbolic, with
$\delta_0 = \frac{10 (n-1)}{K_0} \ln 3 $, by~\cite[Proposition~1.4.3, p.~12]{CDP},
for its sectional curvature is $\le - \frac{K_0^2}{100 (n-1)^2}$, and that its entropy is bounded above by $K_0$, by the Bishop--Gromov Theorem~\ref{BishopGromov}. From this and the fact that $\big(\overline M_i, \bar g_i\big)$
 is quasi-isometric to $\big(\widetilde N_i, \tilde h_i\big)$, we~deduce that
there exist constants $\delta$ and $K_1$ (independent of~$i$) such that $\big(\overline M_i, \bar g_i\big)$ is $\delta$-hyperbolic and has entropy bounded
above by~$K_1$.

As $\big(\widetilde N_i, \tilde h_i\big)$ has Ricci curvature $ \ge - \frac{K_0^2}{n-1}$, it verifies $\Vol \big( B_{\widetilde N_i} (\tilde x,r)\big) \le C_0\cdot {\rm e}^{K_0 r}$ for every
$r > 0$. From this and from the fact that $d_{\tilde h_i} (\tilde x_i, \g \tilde y_i) \ge D_i = d_{\tilde h_i} (\tilde x_i, \tilde y_i)$ for every $\g \in \Gamma_i$, setting $r_i := \frac{3}{4} D_i $, we~deduce
\begin{gather}\label{nonBishop1}
\dfrac{\Vol\big(B_{\overline M_i}(\tilde x_i,2 r_i)\big)}{\Vol \big(B_{\overline M_i}(\tilde x_i,r_i)\big)} = \dfrac{\Vol \big(B_{\overline M_i} \big(\tilde x_i,\frac{3}{2} D_i\big)\big)} {\Vol\big(B_{\overline M_i} \big(\tilde x_i,\frac{3}{4} D_i\big)\big)}
> \dfrac{C_0 \,{\rm e}^{3i^2 K_0}}{\Vol \big(B_{\widetilde N_i} \big(\tilde x_i,\frac{3}{4} D_i\big) \big)} \ge \dfrac{C_0 \,{\rm e}^{3 i^2 K_0}}{C_0 \,{\rm e}^{\frac{3}{4} K_0 \, D_i }}.
\end{gather}

\begin{itemize}\itemsep=0pt
{\sloppy\item
If we choose $D_i =i$, then~\eqref{nonBishop1} gives $\frac{\Vol
(B_{\overline M_i} (\tilde x_i,2 r_i))}{\Vol (B_{\overline M_i}
(\tilde x_i, r_i))} > {\rm e}^{i (3 i-1) K_0}$; hence, for every choice of~the scale $r_0$, of the factor $C> 0$ and of the exponent $K$, and for every
$i > \Max \big(\frac{4}{3} r_0, \frac{\ln C}{K_0}, \frac{K}{K_0}\big)$, one obtains
that $\frac{\Vol ( B_{\overline M_i} (\tilde x_i,2 r_i))}{\Vol (B_{\overline M_i} (\tilde x_i, r_i))} > C {\rm e}^{K r_i}$, despite the fact that $r_i \ge r_0$. It~thus follows that almost all the $\big(\overline M_i, \bar g_i\big)$'s do not verify any weak Bishop--Gromov inequality.
This ends the proof of Example~\ref{nonBishopbis}.

}

\item
If we choose $D_i = D$, we~have $r_i = r_1 = \frac{3}{4} D$ and, for every choice of~$C$ and $K$ and for every
$i > \frac{1}{\sqrt{3 \,K_0}}\big(\ln C + (K + K_0) D \big)^{1/2}$, inequality~\eqref{nonBishop1} gives
$\frac{\Vol(B_{\overline M_i} (\tilde x_i,2 r_1))}{\Vol(B_{\overline M_i} (\tilde x_i, r_1))} \ge {\rm e}^{3 K_0 (i^2 - D/4)} > C{\rm e}^{K r_1}$.
This proves that almost all the $\big(\overline M_i, \bar g_i\big)$'s do not verify any weak Bishop--Gromov inequality at scale $ r_1 =\frac{3}{4} D$.
On the contrary, applying Theorem~\ref{cocompact2}$(i)$, we~know that every $\Gamma_i$-invariant measure $\mu_i$ satisfies the
weak Bishop--Gromov inequality at scale \mbox{$\frac{5}{2} (7 D + 4 \delta)$}, with factor $ C_0 (K_1, D) = 3\cdot 2^{25/4} {\rm e}^{(1+ 6 \ln 2) K_1 D}$ and
exponent $\frac{36}{5} K_1$. This ends the proof of Example~\ref{nonBishop}.
\hfill \qed
\end{itemize}\renewcommand{\qed}{}
\end{proof}

\begin{exa}\label{infinitetopologies} An infinite family $\big(\overline M_i, \bar g_i\big)_{i \in I}$ of complete, $n$-dimensional, $\delta$-hyperbolic Riemannian
manifolds, with distinct topologies and distinct local topologies, each of them admitting a~co-compact proper action, by isometries, of the same group
$\Gamma$ such that the diameter of~$ \Gamma \backslash \overline M_i$ and
the entropy of~$\big(\overline M_i, \bar g_i\big)$ are respectively bounded by constants
$D$ and $K$ independent on~$i$. However, on each $\big(\overline M_i, \bar g_i\big)$, every $\Gamma$-invariant measure $\mu_i$ verifies a weak Bishop--Gromov
inequa\-lity at scale $\frac{5}{2} (7 D + 4 \delta)$, with exponent $\frac{36}{5} K$ and factor $C_0 (K, D)$ (defined in Example~\ref{nonBishop}).
\end{exa}

In this example, the set of indices $I$ may be chosen as the set $\{Y_i\}_{i \in I}$ of all closed, $n$-dimensional manifolds modulo homeomorphisms. A
 consequence is that the local topology of the quotient spaces $M_i := \Gamma_i \backslash \overline M_i$, which is the connected sum of~$\mathbb B^n$ with $Y_i$, is arbitrary.

\begin{proof}[Construction and proof]
We start from any fixed closed $n$-dimensional Riemannian manifold $(M_0,
g_0)$ with Ricci curvature $\ge - \frac{K_0^2}{n-1}$ and
sectional curvature $\le - \frac{K_0^2}{100 (n-1)^2}$. Let $D$ be an upper bound of~$\diam (M_0, g_0) + 1$. We~still denote by $\{Y_i\}_{i \in I}$ the set of all closed, $n$-dimensional manifolds modulo homeomorphisms. For any choice of~$\e \in{} ]0, 1]$, we~can endow
each mani\-fold~$Y_i$ with a Riemannian metric $h_i$ such that $\diam (Y_i,
h_i) \le \e$.

To each choice of such a manifold $Y_i$, we~associate a new closed $n$-dimensional Riemannian manifold $(M_i, g_i)$, by connected sum of~$(M_0, g_0)$
and $(Y_i,h_i)$, constructed as follows: defining $\e_i := 10^{- 6}\cdot \min \big( \inj (Y_i, h_i), \inj (M_0, g_0)\big) \le 10^{- 6} \e$, we~choose points $y_i \in Y_i$ and
$x_0 \in M_0$, we~excise balls $B_{Y_i} (y_i, \e_i)$ and $B_{M_0} (x_0, \e_i)$ from $Y_i$ and $M_0$ (respectively) and glue $M_0 \setminus B_{M_0} (x_0, \e_i)$ and $Y_i \setminus B_{Y_i} (y_i, \e_i)$ by identification of their boundaries.\footnote{To be isometric, this identification may
require to slighly modify the two metrics $g_0$ and $h_i$ on $B_{M_0} (x_0, 2 \e_i)$ and $B_{Y_i} (y_i, 2 \e_i)$ (respectively), in~order that they become flat
on these two balls.} It~is classical that, as~$Y_i$ and $Y_j$ are not homeomorphic when $i\ne j$, then $M_i = M_0 \# Y_i$ and $M_j = M_0 \# Y_j$ are not homeomorphic. Moreover, when $4\e$ is smaller than the injectivity radius of~$(M_0, g_0)$, for every $x$ on the boundary of~$M_0 \setminus B_{M_0} (x_0, \e_i)$ the balls $B_{M_i} (x, 2 \e)$ and $B_{M_j} (x, 2 \e)$ have different topologies.

One can construct a continuous map $f_i \colon M_i \f M_0$ such that $f_i^{-1} (x_0) = Y_i \setminus B_{Y_i} (y_i, \e_i) $ and which is the identity map on the subset
$M_0 \setminus B_{M_0} (x_0, 2 \e_i)$ (included in both spaces $M_0$ and $M_i = M_0 \# Y$), the metric $g_i$ coinciding with $g_0$ on
$M_0 \setminus B_{M_0} (x_0, 2 \e_i)$. We~call $(f_i)_* $ the induced morphism from $\pi_1 (M_i, x_1)$ to $\pi_1 (M_0, x_1)$, where
$x_1$ is any point of~$M_0 \setminus B_{M_0} (x_0, 2 \e_i)$.

Let now $\pi \colon \big(\widetilde M_0, \tilde g_0\big) \f (M_0, g_0)$ be the Riemannian universal covering of~$(M_0, g_0)$ and choose $\tilde x_0 \in \pi^{-1} (x_0)$. Let
$ p_i \colon \big(\overline M_i, \bar g_i\big) \f (M_i, g_i)$ be the Riemannian covering of~$(M_i, g_i)$ such that $(p_i)_* \big(\pi_1 \big(\overline M_i, \bar x_1\big) \big) = \Ker (f_i)_* $ for some point $\bar x_1 \in p_i^{-1} (x_1)$ (see~\cite[Corollary~(6.9)]{GrH}). As~$\Ker (f_i)_*$ is normal in $\pi_1 (M_i, x_1)$, $p_i$ is a Galois covering, and
the group ${\rm Aut} (p_i)$ of automorphisms of this cove\-ring is $\pi_1 (M_i, x_1) / \Ker (f_i)_* $, hence $(f_i)_* $ induces an isomorphism $f_i^*$ from $\pi_1 (M_i, x_1) / \Ker (f_i)_* $ onto the group $\Gamma$ of
automorphisms of the universal covering $\pi$ of~$M_0$.

Classically (see~\cite[Theorem~(6.1)]{GrH}), there then exists a map $\bar f_i \colon \overline M_i \f \widetilde M_0$ such that $\pi \circ \bar f_i =
f_i \circ p_i$
and $ \bar f_i (\g\cdot \bar x) = f_i^* (\g) \cdot \bar f_i (\bar x)$ for every $\bar x \in \overline M_i$ and every $\g \in {\rm Aut} (p_i)$.

Arguing as in the proof of Example~\ref{nonBishopbis}, we~prove that $p_i^{-1} \big( f_i^{-1} (x_0)\big)$ is the union of the left-translates
$ \g \bar f_i^{-1} ( \tilde x_0) $ of~$\bar f_i^{-1} ( \tilde x_0) $, all
isometric to $ f_i^{-1} (x_0) = Y_i \setminus B_{Y_i} (y_i, \e_i) $ via the map $p_i$.
Using the same arguments as in the proof of Example~\ref{nonBishopbis}, we~also obtain that the restriction of~$\bar f_i$ to
$p_i^{-1} \big( M_0 \setminus B_{M_0} (x_0, 2 \e_i)\big) \subset \overline M_i$ is isometric onto $\widetilde M_0 \setminus \bigcup_{\g \in \Gamma} B_{\widetilde M_0} (\g \tilde x_0, 2 \e_i)$.
Revisiting these arguments, one can say that, from a geometric point of view, $\big(\overline M_i, \bar g_i\big)$ is a connected sum of~$\big(\widetilde M_0, \tilde g_0\big)$ with an infinite
number of copies $\big(Y_i^\gamma, h_i^\gamma\big)_{\g \in \Gamma}$ of~$(Y_i, h_i)$, each copy $Y_i^\gamma = \bar f_i^{-1} (\g \tilde x_0) $ of~$Y_i$ being glued
to $\widetilde M_0$ inside
the ball $B_{\widetilde M_0} (\g \tilde x_0, 2 \e_i)$.

Recall that $\big(\widetilde M_0, \tilde g_0\big)$, whose sectional curvature is $\le - \frac{K_0^2}{100 (n-1)^2}$, is $\delta_0$-hyperbolic, with $\delta_0 =
\frac{10 (n-1)}{K_0} \ln 3 $, by~\cite[Proposition~1.4.3, p.~12]{CDP}, and that its entropy is bounded above by $K_0$, by the Bishop--Gromov Theorem~\ref{BishopGromov}. From this and from the fact that $\big(\overline M_i, \bar g_i\big)$ is quasi-isometric to $\big(\widetilde M_0, \tilde g_0\big)$, we~deduce that
there exist constants $\delta$ and $K$ (independent of~$i$) such that $\big(\overline M_i, \bar g_i\big)$ is $\delta$-hyperbolic and has entropy bounded
above by $K$. As~$\Gamma$ acts on $\overline M_i$ via the isomorphic representation of~$\Gamma$ as ${\rm Aut} (p_i)$, one has
$\diam \big(\Gamma \backslash \overline M_i\big) = \diam (M_i,g_i) \le \diam (M_0, g_0) + 1 \le D $.

Now, applying Theorem~\ref{cocompact2}$(i)$, we~know that every $\Gamma$-invariant measure $\mu_i$ satisfies the
weak Bishop--Gromov inequality at scale $\frac{5}{2} (7 D + 4 \delta)$, with factor $ C_0 (K, D) = 3\cdot 2^{25/4} {\rm e}^{(1+ 6 \ln 2) K D}$ and~exponent
$\frac{36}{5} K$. This ends the proof.
\end{proof}

\begin{exa}\label{infiniteBetti}
For every $n \ge 4$ and every integer $k$ such that $ 2 \le k \le n-2$, there exists a~sequence $(M_i, g_i)_{i \in \N}$
of closed, connected, $n$-dimensional Riemannian manifolds, with diameter
$\le D$, whose universal covers $\big(\widetilde M_i, \tilde g_i\big)$ are $\delta$-hyperbolic with entropy $ \le K$ (where the constants~$D$, $\delta$, $K$ are independent on $i$), and such that $\dim \big( H_k (M_i, \R) \big) \f
+ \infty$ when
$i \f +\infty$.
\end{exa}

\begin{proof}[Construction and proof]
We start from any fixed closed $n$-dimensional Riemannian manifold $(M_0,
g_0)$ with Ricci curvature $\ge - \frac{K_0^2}{n-1}$ and
sectional curvature $\le - \frac{K_0^2}{100 (n-1)^2}$. Let $D$ be an upper bound of~$\diam (M_0, g_0) + 2$.

Let $Y$ be any closed, simply connected, $n$-dimensional manifold. For any choice of~$\e \in{} ]0, 1]$, we~can endow $Y$ with a Riemannian metric
$h_{\e}$ such that
$\diam (Y, h_{\e}) \le \e$. To every $i \in \N$, we~associate a new closed $n$-dimensional Riemannian manifold $(M_i, g_i)$, by connected sum of~$(M_0, g_0)$ with $i$ copies of~$(Y, h_{\e})$, constructed as follows: defining $\e_i := \frac{1}{100}\, {\rm e}^{- 2 i } \min \big( \inj (Y, h_{\e})$;
$(\inj (M_0, g_0))^2; 1 \big) \le {\rm e}^{-i } \frac{\e}{10}$, we~choose points $x_1, \dots, x_i \in M_0$ such that, for any pair of distinct indices $j,
l \le i$, one has
$d_{M_0} (x_j, x_l) \ge 2 \sqrt{\e_i} \, (1 + \sqrt{\e_i})^2 \gg 2 \e_i$. To $ M_0 \setminus \cup_{1 \le j \le i} B_{M_0} (x_j, \e_i)$, we glue~$i$ copies
$\breve Y_1, \dots, \breve Y_i$ of~$(Y, h_{\e}) \setminus B_{Y} (y_0, \e_i)$ identifying the boundary of each $\breve Y_j$ with $\partial B_{M_0} (x_j, \e_i)$ (see the proofs of Examples~\ref{nonBishopbis} and~\ref{infinitetopologies}).

Let us define $M_0^{i {\rm \, holes}}:= M_0 \setminus \cup_{1 \le j \le
i} B_{M_0} (x_j, 2 \e_i)$, which is considered as included in both spaces $M_0$
and $M_i$, then $M_i \setminus M_0^{i {\rm \, holes}} $ has $i$ connected
components, the connected component which contains $\breve Y_j $ (denoted
by $\widehat Y_j$) being the result of the connected sum (constructed above) of~$ B_{M_0} (x_j, 2 \e_i)$ with the $j$-th copy of~$Y$.
One can construct a continuous map $f_i \colon M_i \f M_0$ such that $f_i^{-1} (x_j) = \breve Y_j $, which maps each connected component $\widehat Y_j$
of~$M_i \setminus M_0^{i {\rm \, holes}}$ onto $B_{M_0} (x_j, 2 \e_i)$,
and which is the identity map on the subset $M_0^{i {\rm \, holes}}$, the
metric $g_i$ coinciding with $g_0$ on $M_0^{i {\rm \, holes}}$.

As $Y$ is simply connected, the converse of the induced morphism $\pi_1 (M_i, x_0) \f \pi_1 (M_0, x_0)$, where $x_0 \in M_0^{i {\rm \, holes}}$,
induces an isomorphism, denoted by $\varrho_i $,
between the groups $\Gamma$ and $\Gamma_i$ of deck-transformations of the
universal coverings $\pi_0 \colon \big(\widetilde M_0, \tilde g_0\big)
\f (M_0, g_0)$ and $\pi_i \colon \big(\widetilde M_i, \tilde g_i\big) \f (M_i, g_i)$. As~in the proofs of Examples~\ref{nonBishopbis} and~\ref{infinitetopologies}, there exists a map $\tilde f_i \colon \widetilde M_i \f \widetilde M_0$ such
that $\pi_0 \circ \tilde f_i = f_i \circ \pi_i$ and $ \tilde f_i \big(\varrho_i (\g) \cdot \tilde x\big) = \g \cdot \tilde f_i (\tilde x)$ for every $\tilde x \in \widetilde M_i$ and every $\g \in \Gamma$.

Arguing as in the proof of Examples~\ref{nonBishopbis} and~\ref{infinitetopologies}, we~prove that, for each $j$ ($1 \le j \le i$), denoting by $\tilde x_j$ any
point of~$\pi_0^{-1} (x_j)$, the set $\pi_i^{-1} \big( f_i^{-1} (x_j)\big)$ is the union of the left-translates $ \varrho_i (\g) \tilde f_i^{-1} (
\tilde x_j) $ of~$\tilde f_i^{-1} ( \tilde x_j) $ for all the
$\g \in \Gamma$, each of these translates being isometric to $ f_i^{-1} (x_j) = \breve Y_j $ via the map $\pi_i$. As~in the proof of Examples~\ref{nonBishopbis} and~\ref{infinitetopologies}, we~prove that the restriction of~$\tilde f_i$ to
$\pi_i^{-1} \big( M_0^{i {\rm \, holes}}\big) \subset \widetilde M_i$ is one to one isometric onto $\pi_0^{-1} \big( M_0^{i {\rm \, holes}}\big) =
\widetilde M_0 \setminus \bigcup_{\g \in \Gamma} \big(\cup_{1 \le j \le i} B_{\widetilde M_0} (\g \tilde x_j, 2 \e_i)\big)$ that we shall denote by $\widetilde M_0^{i {\rm \, holes}}$ for the sake of simplicity.

Revisiting these arguments, one can say that, from a geometric point of view, $\big(\widetilde M_i, \tilde g_i\big)$ is a~connected sum of~$(\widetilde M_0, \tilde g_0)$ with an infinite
number of copies $\big(Y_j^\gamma, h_j^\gamma\big)_{\g \in \Gamma}$ of~$(Y, h_\varepsilon)$, each copy $Y_j^\gamma = \tilde f_i^{-1} (\g \tilde x_j) =
\varrho_i (\g) \tilde f_i^{-1} (\tilde x_j) $ being glued to $\widetilde M_0$ inside the ball $B_{\widetilde M_0} \big(\g \tilde x_j, \e_i \big)$.

Considering that $\widetilde M_0^{i {\rm \, holes}}$ is also included in $\widetilde M_i$ via
the inverse of the isometric map $\tilde f_i$: $ \pi_i^{-1} \big( M_0^{i {\rm \, holes}}\big) \f \widetilde M_0^{i {\rm \, holes}}$, the two inclusion maps from $\widetilde M_0^{i {\rm \, holes}}$ into $\widetilde M_0 $ and $\widetilde M_i $ preserve the lengths of paths. On
$\widetilde M_0^{i {\rm \, holes}}$, we~compare the length-distance $d_{\widetilde M_0^{i {\rm \, holes}}}$ $\big($associated to paths lying in $\widetilde M_0^{i {\rm \, holes}}\big)$ with the restriction
of the two distances $d_{\widetilde M_0 }$ and $d_{\widetilde M_i }$ of~$\widetilde M_0$ and $\widetilde M_i$ $\big($which are length-distances associated to
paths lying in $\widetilde M_0 $ and $\widetilde M_i $ respectively$\big)$; it is clear that $ d_{\widetilde M_0 }$, $d_{\widetilde M_i } \le d_{\widetilde M_0^{i {\rm \, holes}}}$.

Between two points $\tilde x, \tilde y$ of~$\widetilde M_0^{i {\rm \, holes}}$, let $c$ be a minimizing geodesic of~$\big(\widetilde M_0, \tilde g_0\big)$, let $\bar c$ be
the disconnected union of geodesic segments which is the intersection of~$c$ with $\widetilde M_0^{i {\rm \, holes}}$ and $\tilde c$ the path of minimal length
among all the continuous paths which coincide with $\bar c$ on the interior of~$\widetilde M_0^{i {\rm \, holes}}$ and which connect the different
connected
components of~$\bar c$ by means of several arcs of circle, each of these arcs lying in one of the $i$ geodesic spheres which are the connected components of
the boundary of~$\widetilde M_0^{i {\rm \, holes}}$; as the length of each of these arcs of circles is at most $2 \pi \e_i (1 + \e_i)$ (when~$i$ is large enough),
and as the length of each connected component of~$\bar c$ (except the first and the last ones) is at least $2 \sqrt{\e_i} (1+\e_i)$, we~get that
${\rm length} (\tilde c) \le (1 + \pi \sqrt{\e_i})\, {\rm length} (\bar
c) + 2 \pi \e_i (1 + \e_i)$ and it follows that
$d_{\widetilde M_0^{i {\rm \, holes}}}(\tilde x,\tilde y) \le (1 + \pi \sqrt{\e_i})\, d_{\widetilde M_0 } (\tilde x, \tilde y) + 2 \pi \e_i (1 + \e_i)$. We~similarly prove that
$d_{\widetilde M_0^{i {\rm \, holes}}}(\tilde x,\tilde y) \le (1 + \pi \sqrt{\e_i})\, d_{\widetilde M_i } (\tilde x, \tilde y) + 2 \pi \e_i (1 + \e_i)$. A first consequence is that
\begin{gather*}
\Ent \big(\widetilde M_i, \tilde g_i\big) \le \big(1 + \pi \sqrt{\e_i}\big)
\Ent \big(\widetilde M_0^{i {\rm \, holes}}, d_{\widetilde M_0^{i {\rm \, holes}}}\big)
\le \big(1 + \pi \sqrt{\e_i}\big) \Ent \big(\widetilde M_0, \tilde g_0\big).
 \end{gather*}
Let now $\tilde x$, $\tilde y$ be any pair of points of~$\widetilde M_i$, $c$ be a minimizing geodesic of~$\big(\widetilde M_i, \tilde g_i\big)$, and $t_0$ (resp.~$t_1$) be the infimum (resp.\ the supremum) of the values of~$t$ such that $c(t) \in \widetilde M_0^{i {\rm \, holes}}$. We~then have
\begin{gather*}
d_{\widetilde M_i} (\tilde x, \tilde y) \le d_{\widetilde M_0^{i {\rm
\, holes}}} (c(t_0), c(t_1)) + 2 (\e + 4 \e_i)
\\ \hphantom{d_{\widetilde M_i} (\tilde x, \tilde y)}
{}\le (1 + \pi \sqrt{\e_i}) \, d_{\widetilde M_0}
( c(t_0), c(t_1)) + 2 \pi \e_i (1 + \e_i) + 2 \e + 8 \e_i
\\ \hphantom{d_{\widetilde M_i } (\tilde x, \tilde y)}
{}\le (1 + \pi \sqrt{\e_i}) \big(d_{\widetilde M_0 } \big(\tilde f_i (\tilde x), \tilde f_i (\tilde y) \big) + 8 \e_i \big) + 2 \pi \e_i (1 + \e_i) + 2 \e + 8 \e_i,
\end{gather*}
where the last inequality uses the fact that $c(t_0)$ and $\tilde x$ (resp.\ $c(t_1)$ and $\tilde y$), if they do not coincide, belong to the same
connected
component of~$ \widetilde M_i \setminus \widetilde M_0^{i {\rm \, holes}}$, thus that $\tilde f_i (c(t_0))$ and $\tilde f_i (\tilde x)$ $\big($resp.\ $\tilde f_i (c(t_1))$
and $\tilde f_i (\tilde y)\big)$ belong to the same connected component of~$ \widetilde M_0 \setminus \widetilde M_0^{i {\rm \, holes}}$, thus to the same ball
of radius $2 \e_i$. Analogous arguments yield
\begin{gather*}
d_{\widetilde M_0}\big(\tilde f_i(\tilde x),\tilde f_i (\tilde y)\big) \le \big(1+ \pi\sqrt{\e_i}\big) \big(d_{\widetilde M_i } (\tilde x,\tilde y) + 2 \e + 8 \e_i \big) + 2 \pi \e_i (1 + \e_i) + 8 \e_i.
\end{gather*}
This proves that $\tilde f_i $ is a quasi-isometry from $\big(\widetilde M_i, \tilde g_i\big)$ onto $\big(\widetilde M_0, \tilde g_0\big)$. Recall that $\big(\widetilde M_0, \tilde g_0\big)$ is
$\delta_0$-hyperbolic with $\delta_0 = \frac{10 (n-1)}{K_0} \ln 3 $, because its sectional curvature is $\le - \frac{K_0^2}{100 (n-1)^2}$ (by~\cite[Proposition~1.4.3, p.~12]{CDP}), and that it has entropy bounded
above by $K_0$ (by the Bishop--Gromov Theorem~\ref{BishopGromov}). This quasi-isometric comparison implies that
there exist constants $\delta$ and $K$ (independent on $i$) such that $\big(\widetilde M_i, \tilde g_i\big)$ is $\delta$-hyperbolic and has entropy bounded
above by $K$. As~$\Gamma$ acts on $\widetilde M_i$ via the isomorphic representation $\varrho_i$, one has
$\diam \big(\Gamma \backslash \widetilde M_i\big) = \diam (M_i,g_i) \le \diam (M_0, g_0) + 2 \le D $. This ends the first part of the proof.

\medskip

In the sequel the homology groups are with coefficients in $\R$. Let us prove that $\dim \big( H_k (M_i) \big)$ $\f + \infty$ with $i$ when $ 2 \le k \le n-2$.
From the topological viewpoint, as~$M_{i+1} = M_i \# Y$, it is the union of two open subsets $M_i^*$ and $Y^*$, where $M_i^*$ (resp.\ $Y^*$) is
obtained
by removing from $M_i$ (resp.\ from $Y$) a small ball contained in $B_{M_i} (x_{i+1}, \e_i)$ (resp.\ in $B_Y (y_0, \e_i)$); then, using the fact that
$M_i^* \cap Y^*$ is homeomorphic to $]{-} \e_i, \e_i[{} \times \mathbb S^{n-1}$, the Mayer--Vietoris exact
sequence associated to the covering $M_i^*\cup Y^*$ of~$M_{i+1}$
reads
\begin{gather*}
\cdots \f H_{k+1} (M_i \# Y) \f H_k \big(\mathbb S^{n-1}\big) \f H_k ( M_i^*) \oplus H_k (Y^*) \f H_k (M_i \# Y)
\\ \hphantom{\cdots}
\f H_{k-1} (\mathbb S^{n-1})\f \cdots.
\end{gather*}
As $ n \ge 4$ and $2 \le k \le n-2$, one has $H_{k} \big(\mathbb S^{n-1}\big) =
\{0\} = H_{k-1} \big(\mathbb S^{n-1}\big)$ and $H_{k} (\mathbb B^{n}) = \{0\}$, hence
\begin{gather*}
H_k (M_{i+1}) = H_k (M_i \# Y) \simeq H_k ( M_i^*) \oplus H_k (Y^*) \simeq H_k ( M_i) \oplus H_k (Y).
\end{gather*}
Let us consider three examples of manifolds $Y$:

\begin{itemize}\itemsep=0pt
\item when $n$ is even ($n = 2 d \ge 4$), if we choose $Y = \C P^d$, then $H_k (Y) \simeq \R$ for every even integer $k$ such that $ 2 \le k \le n-2$ and
we have, in~this case,
$\dim \big( H_k (M_{i+1})\big) = \dim \big( H_k (M_{i})\big) + 1$; this ends the proof in this case,

\item when $n$ is odd ($n = 2 d + 1\ge 5$), if we choose $Y = \C P^{d-1} \times \mathbb S^3$, we~have $\dim H_k (Y) \ge 1 $ for every
$k$ such that $ 2 \le k \le n-2$ and thus
$\dim \big( H_k (M_{i+1})\big) \ge \dim \big( H_k (M_{i})\big) + 1$,
and this ends the proof in this case,

\item when $n$ is even ($n = 2 d \ge 8$), if we choose $Y = \C P^{d-3} \times \mathbb S^3 \times \mathbb S^3$, then $\dim H_k (Y) \ge 1 $ for
every $k$
such that $ 2 \le k \le n-2$ (except when $d = 4$ and $ k = 4$) and then $\dim \big( H_k (M_{i+1})\big) \ge \dim \big( H_k (M_{i})\big) + 1$, and this ends the proof in this case.
\hfill \qed
\end{itemize}\renewcommand{\qed}{}
\end{proof}

\pdfbookmark[1]{References}{ref}
\LastPageEnding

\end{document}